\documentclass[12pt,oneside]{article}
\RequirePackage{epsfig}
\RequirePackage{amsmath}
\RequirePackage{amsfonts,amssymb,amsfonts}
\RequirePackage{amsthm}

\usepackage{color}
\usepackage{graphicx}
\usepackage{graphics}
\usepackage{subfigure}
\usepackage{footnote}

\newtheorem{theorem}{Theorem}
\newtheorem{proposition}{Proposition}

\newtheorem{prb}[theorem]{Problem}
\newcommand{\be}{\begin{equation}}
\newcommand{\ee}{\end{equation}}
\newcommand{\bi}{\begin{itemize}}
\newcommand{\ei}{\end{itemize}}
\newcommand{\ben}{\begin{enumerate}}
\newcommand{\een}{\end{enumerate}}
\newcommand{\ba}{\begin{eqnarray}}
\newcommand{\ea}{\end{eqnarray}}
\newcommand{\ZZ}{\mathbb{Z}}
\newcommand\bp{\mathbf{p}}
\newcommand\bq{\mathbf{q}}
\newcommand\bt{\mathbf{t}}
\newcommand\bl{\mathbf{l}}
\newcommand\bm{\mathbf{m}}
\newcommand\bs{\mathbf{s}}
\newcommand\bg{\mathbf{g}}

\def\pp{{\bf p}}
\def\uu{{\bf u}}
\def\mm{{\bf m}}
\def\ll{{\bf l}}
\def\tt{{\bf t}}

\def\eee{{\bf e}}

\catcode`\à=\active \def à{\` a} \catcode`\â=\active \def â{\^ a}
\catcode`\é=\active \def é{\' e} \catcode`\è=\active \def è{\` e}
\catcode`\ê=\active \def ê{\^ e} \catcode`\ù=\active \def ù{\` u}
\catcode`\ô=\active \def ô{\^ o} \catcode`\û=\active \def û{\^ u}
\catcode`\ç=\active \def ç{\c c} \catcode`\î=\active \def î{\^ \i}
\catcode`\ï=\active \def ï{\" \i} \catcode`\À=\active \def À{\` A}

\begin{document}

\newpage

\begin{center}
\textbf{CONVEXITY PRESERVING INTERPOLATORY SUBDIVISION WITH CONIC PRECISION}

\bigskip

\textsc{Gudrun Albrecht}\\
\footnotesize
\textrm{Univ. Lille Nord de France, UVHC, LAMAV-CGAO\\
FR no. 2956, F-59313 Valenciennes, France}\\
\textit{gudrun.albrecht@univ-valenciennes.fr}

\bigskip
\normalsize
\textsc{Lucia Romani}\\
\footnotesize
\textrm{Univ. of Milano-Bicocca, Dept. of Mathematics and Applications\\
Via R. Cozzi 53, 20125 Milano, Italy}\\
\textit{lucia.romani@unimib.it}

\normalsize

\bigskip

\begin{abstract}
The paper is concerned with the problem of shape preserving
interpolatory subdivision. For arbitrarily spaced, planar input data
an efficient non-linear subdivision algorithm is presented that
results in $G^1$ limit curves, reproduces conic sections and
respects the convexity properties of the initial data. Significant
numerical examples illustrate the effectiveness of the proposed
method.



\end{abstract}

\end{center}

\bigskip

\normalsize

{\renewcommand{\thefootnote}{}

\footnotetext{Date: January 26, 2010}
}

\section{Introduction and state of the art}

Subdivision schemes constitute a powerful alternative for the design
of curves and surfaces over the widely studied parametric and
implicit forms. In fact, they offer a really versatile tool that is,
at the same time, very intuitive and easy to use and implement. This
is due to the fact that subdivision schemes are defined via
iterative algorithms which exploit simple refinement rules to
generate denser and denser point sequences that, under appropriate
hypotheses, converge to a continuous, and potentially smooth,
function.

In the univariate case, the iteration starts with a sequence of
points denoted by $\bp^0=(\bp_i^0 \ : \ i\in \ZZ)$, attached to the
integer grid, and then for any $k\geq0$ one subsequently computes a
sequence $\bp^{k+1}=S\bp^k$, where $S:\ell(\ZZ) \rightarrow
\ell(\ZZ)$ identifies the so-called subdivision operator.

Subdivision operators can be broadly classified into two main
categories: {\sl interpolating} and {\sl approximating}
\cite{DL02,ZS00}. Interpolating schemes are required to generate
limit curves passing through all the vertices of the given polyline
$\bp^0$. Thus they are featured by refinement rules maintaining the
points generated at each step of the recursion in all the successive
iterations. Approximating schemes, instead, are not required to
match the original position of vertices on the assigned polyline and
thus they adjust their positions aiming at very smooth and visually
pleasing limit shapes. As a consequence, while in the case of {\sl
approximating} subdivision the refinement rules rely on a recursive
corner cutting process applied to the starting polygon $\bp^0$, in
the case of {\sl interpolatory} subdivision, in every iteration a
finer data set $\bp^{k+1}$ is obtained by taking the old data values
$\bp^k$ and inserting new points in between them. Every such new
point is calculated using a finite number of existing, usually
neighboring points. In particular, if the computation of the new
points is carried out through a linear combination of the existing
points, the scheme is said to be {\sl linear}, otherwise {\sl
non-linear}.

Then, inside the above identified categories, the schemes can also be further classified.
More specifically, they can be distinguished between {\sl
stationary} (when the refinement rules do not depend on the
recursion level) and {\sl non-stationary}; between {\sl uniform}
(when the refinement rules do not vary from point to point) and {\sl
non-uniform}; between {\sl binary} (when the number of points is
doubled at each iteration) and {\sl $N$-ary}, namely of arity $N>2$.

Most of the univariate subdivision schemes studied in the literature
are binary, uniform, stationary and linear. These characteristics,
in fact, ease to study the mathematical properties of the limit
curve, but seriously limit the applications of the scheme.
Exceptions from a binary, or a uniform, or a stationary, or a linear
approach, have recently appeared (see for example \cite{S05} and
references therein), but none of the proposed methods provides an
interpolatory algorithm that can fulfill the list of {\sl all}
fundamental features considered essential in applications. These
features can briefly be summarized as:

\begin{itemize}
\item[(i)] generating a visually-pleasing limit curve which faithfully mimics the behaviour
of the underlying polyline without creating unwanted oscillations;
\item[(ii)] preserving the shape, i.e., the convexity properties of the given data;
\item[(iii)] identifying geometric primitives like circles and more generally conic sections,
the starting polyline had been sampled from, and reproducing them.
\end{itemize}

\noindent Requirement (i) derives from the fact that, despite
interpolating schemes are considered very well-suited for handling
practical models to meet industrial needs (due to their evident link
with the initial configuration of points representing the object to
be designed), compared to their approximating counterparts, they are
more difficult to control and tend to produce bulges and unwanted
folds when the initial data are not uniformly spaced. Recently this
problem has been addressed by using non-uniform refinement rules
\cite{BCR09b} opportunely designed to take into account the
irregular distribution of the data. But, although their established
merit of providing visually pleasing results, there is no guarantee
that such methods are convexity-preserving, i.e., that if a convex
data set is given, a convex interpolating curve can be obtained.
 This is due to the fact that, such non-uniform schemes
are linear and, as it is well-known \cite{K98}, linear refinement
operators that are $C^1$ cannot preserve convexity in general.

The property (ii) of convexity preservation is of great practical
importance in modelling curves and surfaces tailored to industrial
design (e.g. related to car, aeroplane or ship modelling where
convexity is imposed by technical and physical conditions as well as
by aesthetic requirements). In fact, if shape information as
convexity is not enforced, interpolatory curves, though smooth, may
not be satisfactory as they may contain redundant wiggles and bumps
rather than those suggested by the data points, i.e., they feature
unacceptable visual artifacts. Preserving convexity, while a curve
is interpolated to a given data set, is far from trivial. But much
progress has been made in this field, evidence of which is given by
the recent burgeoning literature. In most publications, the
introduction of subdivision schemes fulfilling requirement (ii) has
been achieved through the definition of non-linear refinement rules.
In fact, despite linear subdivision schemes turn out to be simple to
implement, easy to analyze and affine invariant, they have many
difficulties to control the shape of the limit curve and avoid
artifacts and undesired inflexions that usually occur when the
starting polygon $\bp^0$ is made of highly non-uniform edges.
Non-linear schemes, instead, offer effective algorithms to be used
in shape-preserving data interpolation
\cite{cinesi,DLL92,D05,K98,KVD00,MDD05}.

On the basis of the well--known, linear Dubuc-Deslauriers
interpolatory 4-point scheme \cite{DD89} for example, several
non-linear analogues have been presented in order to solve at least
one of the three above listed properties. On the one hand,
non-linear modifications of the classical 4-point scheme have been
introduced to reduce the oscillations that usually occur in the
limit curve when applying the refinement algorithm to polylines with
short and long adjacent edges. These have been presented in
\cite{CDM03} and \cite{DFH09}, and as concerns the case of
convexity-preserving strategies (which are the ones capable of
completely eliminating the artifacts arising during the subdivision
process), we find the papers \cite{FM98} and \cite{MDD05}. On the
other hand, for the purpose of enriching the Dubuc-Deslauriers
4-point scheme with the property (iii) of geometric primitives
preservation, a non-linear 4-point scheme reproducing circles and
reducing curvature variation for data off the circle, has been
defined \cite{SD05}. With the same intent, another modification of
the classical 4-point scheme in a non-linear fashion, had been given
in \cite{AEV03}.

With these papers, the theoretical investigation of non-linear
interpolatory subdivision has only begun. A lot is still to be done,
in particular as concerns the use of non-linear rules for
reproducing salient curves other than circles, considered of
fundamental importance in several applications. So far, it has been
shown that non-linear updating formulas can be used in the
definition of non-stationary subdivision schemes aimed at
reproducing polynomials and some common transcendental functions. In
particular, \cite{R09b} respectively \cite{BCR07,BCR09a,R09a}
present subdivision algorithms that turn out to be circle-preserving
respectively able to exactly represent any conic section. While the
first is able to guarantee reproduction starting from given samples
with any arbitrary spacing, for the latter ones the property of
conics precision is confined to the case of equally-spaced samples.
Most recently a shape and circle preserving scheme has been
presented in \cite{cinesi}.

Therefore, an outstanding issue that should be considered is the
possibility of defining an interpolatory subdivision scheme that is
at the same time shape-preserving and artifact free, as well as
capable of generating conic sections starting from any
arbitrarily-spaced samples coming from a conic. This is exactly the
purpose of this paper. Based on an approximation order four strategy
presented in \cite{ABFH08} for estimating tangents to planar convex
data sequences, we propose a convexity-preserving interpolatory
subdivision scheme with conic precision. This turns out to be a new
kind of non-linear and geometry-driven subdivision method for curve
interpolation. In Section \ref{def_scheme} we start by describing the refinement strategy
which relies on a classical cross-ratio property for conic sections and uses the
tangent estimator from \cite{ABFH08}, for the case of globally convex data.
In Section \ref{nonconvex} we adapt the scheme to general, not necessarily convex data by
segmenting the given polygon into convex segments and by presenting
new refinement rules next to the junction points. In Section
\ref{algorithm} we summarize the whole subdivision algorithm in all its steps.
Then, in Section \ref{secadapt} we present an adaptive variant of the scheme which is
aimed at producing regularly spaced points in every round of subdivision.
Section \ref{properties} contains proofs for the scheme's shape preservation and conic
reproducing properties, as well as for the $G^1$ continuity of its
limit curve. Section \ref{examples} is devoted to illustrating the
scheme by several significant application examples, and we conclude
in Section \ref{concl}.

\section{Definition of the scheme for globally convex data}\label{def_scheme}

In this section we define a convexity preserving subdivision scheme
for globally convex data which will then be the basis for the final
shape preserving subdivision algorithm for general data.

Curve subdivision schemes iteratively apply a subdivision operator
$S$ to a starting point sequence $\bp^0=(\bp_i^0 \ : \ i\in \ZZ)$
yielding a new sequence $\bp^{k+1}=S\bp^k$ for any level $k \geq 0$. Our scheme being
interpolatory, has refinement rules of the following form:

\begin{equation} \label{def_phi}
\begin{array}{ll}
&\bp_{2i}^{k+1}=\bp_i^k, \\
&\bp_{2i+1}^{k+1}=\varphi(\bp_{i-\nu}^k,...,\bp_i^k, \bp_{i+1}^k,
...,\bp_{i+\nu}^k, \bp_{i+\nu+1}^k ; \; \bp)
\end{array}
\end{equation}

\smallskip\noindent
where $\nu =2$ is the number of points taken into account in the
left and right hand neighborhoods of the segment $\bp^k_i
\bp_{i+1}^k$ in order to define the newly inserted vertex
$\bp_{2i+1}^{k+1}$, and $\bp$ is a parameter point specified later.
$\varphi$ is a non linear function, which we will define in form of
an algorithm.

In order to detail the idea of this convexity--preserving scheme, we
thus consider the following problem where, for simplicity, we omit
the upper indices $k$.

\smallskip
\begin{prb} \label{prb1}
Given $n$ points $\pp_i((p_i)_x,(p_i)_y)$, $i=1,\ldots, n$ $(n \ge
5)$, in convex position in the affine plane, we wish to obtain one
new point $\uu_i$ related to the $i$-th edge $\pp_i \, \pp_{i+1}$.
\end{prb}

\smallskip
We will carry out the construction of the new points in the
projectively extended affine plane. To this end we denote the
projective counterparts of the affine points
$\pp_i((p_i)_x,(p_i)_y)$, $i=1,\ldots, n$ by $P_i(p_{i,0},p_{i,1},
p_{i,2})$ where $p_{i,0}=1, p_{i,1} = (p_i)_x, p_{i,2} = (p_i)_y$.

By projective geometry's principle of duality, a line $L$ may be
represented either by a linear equation
\[
l_0 x_0 + l_1 x_1 + l_2 x_2 = 0
\]
in variable {\it point coordinates} $(x_0, x_1, x_2)$ or by a triple
$(l_0, l_1, l_2)$ of constant {\it line coordinates}. The line
coordinates of the line $L(l_0, l_1, l_2)$ joining two points
$X_1(x_{1,0}, x_{1,1},$ $x_{1,2})$ and $X_2(x_{2,0}, x_{2,1},
x_{2,2})$ may simply be calculated by the vector product $X_1 \wedge
X_2 = L$. In the same way, the point coordinates of the intersection
point $P(p_0, p_1, p_2)$ of two lines $L_1(l_{1,0}, l_{1,1},
l_{1,2})$ and $L_2(l_{2,0}, l_{2,1}, l_{2,2})$ is obtained as $L_1
\wedge L_2 = P$. Without loss of generality we apply a normalization
to the homogeneous point coordinates such that $x_0 \in \{0, 1\}$
for all calculated points.

In order to preserve global convexity of the points we apply the
following preprocessing procedure to the point set. In every given
point $P_i$ we estimate a tangent from a subset of five points
(including the point $P_i$) by the conic tangent estimator presented
in \cite{ABFH08}. If the given points represent a closed polygon,
then the five-point subset is composed of the point $P_i$, its
preceding two points $P_{i-1}, P_{i-2}$ as well as its successive
two points $P_{i+1}, P_{i+2}$ (by considering $P_0 = P_n$, $P_{-1} =
P_{n-1}$, $P_{n+1} = P_1$, $P_{n+2} = P_2$). If the given points
represent an open polygon, then for $i=3, \ldots, n-2$ we proceed as
above and for $i \in \{1,2\}$ (respectively $i \in \{n-1,n\}$) the
points $\{P_1, P_2, P_3, P_4, P_5\}$ (respectively $\{P_{n-4},
P_{n-3}, P_{n-2}, P_{n-1}, P_{n}\}$) are taken.

\medskip
\begin{figure}[ht!]
\centering
\resizebox{5.5cm}{!}{\includegraphics{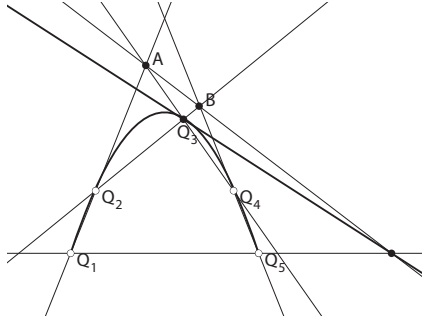}}
\caption{Illustration of the tangent construction.}
\label{fig1}
\end{figure}
\medskip

\begin{figure}[ht!]
\centering
\resizebox{4.0cm}{!}{\includegraphics{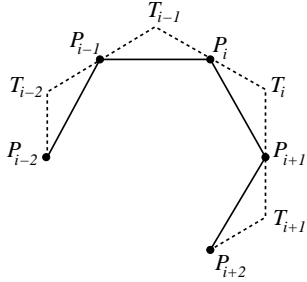}}
\caption{Convex delimiting polygon (dashed line) for a given convex data set (solid line).}
\label{fig1b}
\end{figure}

In order to apply the conic tangent estimator from \cite{ABFH08} we locally rename the
five points around $P_i$ by $Q_3 = P_i$, and by arbitrarily mapping
the remaining four points to $Q_1, Q_2, Q_4, Q_5$ by a one--to--one
map. The desired tangent in the point $Q_3$ is then calculated by
the formula (see \cite{ABFH08}) \be \label{est_tang} M_{33} := Q_3
\wedge (M_{15} \wedge (A \wedge B))\, , \ee where $M_{ij} = Q_i
\wedge Q_j$ for $i \ne j$, $A = M_{12} \wedge M_{34}$, and $B =
M_{54} \wedge M_{32}$. See Figure \ref{fig1} for an illustration.\\
We then denote the obtained line in the point $P_i$ by $L_i$ and
intersect every two consecutive lines generating the intersection
points
\[
T_i = L_i \wedge L_{i+1}\, ,
\]
see Figure \ref{fig1b}.

\medskip
\medskip

The dashed lines in Figure \ref{fig1b} constitute a convex
delimiting polygon for the new points generated in the next
subdivision level. If the initial points come from a conic section,
the constructed lines $L_i$ are the tangents to this conic in the
respective points $P_i$. Otherwise the lines $L_i$ approximate the
tangents with approximation order $4$, see \cite{ABFH08}.

After this preprocessing step we now get back to the initial
subdivision problem \ref{prb1}, i.e., between every two points $P_i$
and $P_{i+1}$ insert a new point $U_i$ by applying a classical
result from projective geometry which, for the readers convenience, we recall in the following
theorem (see e.g., \cite{A08, C65, P10}).

\medskip
\begin{theorem} \label{thrmcr}

\begin{enumerate}
\item[a)]
Let $X, E, E_0, E_1$ be four points of a projective line ${\cal
P}_1$, where the points $E, E_0, E_1$ are mutually distinct, and let
$(x_0, x_1)$ be the projective coordinates of the point $X$ with
respect to the projective coordinate system $\{E_0, E_1; E\}$ of
${\cal P}_1$. Then, the cross ratio $cr(X, E, E_0, E_1)$ of the four
points $X, E, E_0, E_1$ in this order is defined by
\[
cr(X, E, E_0, E_1)= \frac{x_1}{x_0}.
\]
\item[b)]
Let $P_1, P_2$ be two points on a conic section $c$, and $t_1, t_2$
the tangents of $c$ in $P_1, P_2$ respectively, and let $T$ be the
intersection point of $t_1$ and $t_2$ $(T=t_1 \cap t_2)$. Then, the
point $T$ and the line $P_1 P_2$ are pole and polar with respect to
the conic $c$. Let's further denote the intersection points of any
line $l_T$ through $T$ with the conic $c$ by $P$ and $U$, and the
intersection point of $l_T$ and $T$'s polar $P_1 P_2$ by $X$ $(l_T
\cap c = \{P,U\}, \, l_T \cap P_1 P_2 = X)$. Then
\[
cr(U,P,X,T) = -1,
\]
and the four points $(U,P,X,T)$ are said to be in harmonic position.
\end{enumerate}
\end{theorem}

\begin{figure}[ht!]
\centering
\resizebox{3.5cm}{!}{\includegraphics{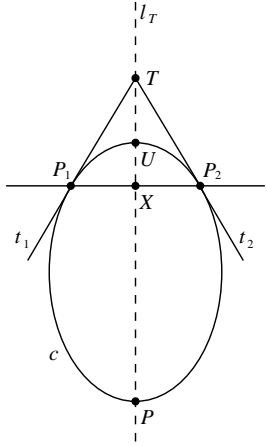}}
\caption{Harmonic cross ratio condition for conic sections.}
\label{fig2}
\end{figure}

\medskip
Theorem \ref{thrmcr} gives us the means of constructing a point $U$
from three known collinear points $P,X,T$ such that the harmonic
cross ratio condition for conic sections is satisfied. For an illustration see
Figure \ref{fig2}.

\medskip
\medskip

\begin{figure}[ht!]
\centering
\resizebox{4.0cm}{!}{\includegraphics{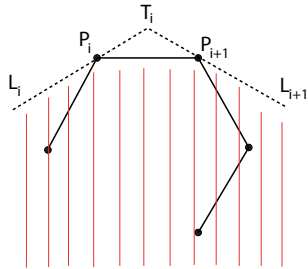}}
\caption{Suitable region for the choice of the parameter point $P$.}
\label{fig_region}
\end{figure}

For the choice of the parameter point $P$ which is needed for the
construction of the new point $U_i$ inside the triangle $\Delta P_i
T_i P_{i+1}$, the whole region bounded by the lines $L_i$, $L_{i+1}$
and $P_i P_{i+1}$ containing the given convex polygon (see Figure
\ref{fig_region}) is suitable. In particular, every point $P_{j_i}$
of the given convex polygon can be taken as parameter point (with
exception of $P_i$ and $P_{i+1}$), and since such a choice
guarantees the reproduction of conic sections we opt for it.

\medskip
\medskip

Let's thus denote by $X_i$ the intersection point of the lines
$P_iP_{i+1}$ and $P_{j_i} T_i$ for a chosen $j_i \in \{ 1, \ldots,
n\} \setminus \{i, i+1\}$, i.e., $X_i = P_iP_{i+1} \cap P_{j_i}
T_i$.

In order to guarantee a regular distribution of the inserted points
we propose to choose the index $j_i$ by the following angle criterion.
To this end we temporarily come back to the
Euclidean plane and introduce the midpoint $\mm_i$ for each segment
$\pp_i \pp_{i+1}$. Let $g_i = \tt_i \mm_i$, $h_{ij} = \tt_i
\pp_j$ be the connecting lines of the points $\tt_i$ and $\mm_i$
respectively $\tt_i$ and $\pp_j$, and $\alpha_j^i = \angle(g_i,
h_{ij})$ the angle between these two lines\footnote{The smaller one
of the two complementary angles is taken in each case.} for $j =
i+2, \ldots, n+i-1$ by considering $\pp_{n+r} = \pp_r$ for $r \ge
1$. For every $i \in \{1, \ldots, n\}$ in the case of a closed
polygon and for every $i \in \{1, \ldots, n-1\}$ in the case of an
open polygon, we then obtain a value $j_i$ from the condition \be
\label{cond_ki_angle} \alpha_{j_i}^i = \min_{j = i+2, \ldots, n+i-1}
\alpha_j^i. \ee

Once the point $P_{j_i}$ has been selected by exploiting the illustrated criterion,
we then establish the projective coordinate system $\{X_i, T_i;
P_{j_i}\}$ on the straight line $P_{j_i} T_i$ by calculating the
projective representatives of $X_i$ and $T_i$ by solving
\[
\gamma_i X_i + \mu_i T_i = P_{j_i}
\]
for $\gamma_i$ and $\mu_i$. We obtain
\[
\gamma_i = D_{i,1}/D_i \, , \ \ \mu_i = D_{i,2}/D_i \, ,
\]
where \be \label{Ds} D_i=\det \left( \begin{array}{cc} x_{i,l} &
t_{i,l}
\\ x_{i,m} & t_{i,m}
\end{array}\right) \, ,\;
D_{i,1}=\det \left( \begin{array}{cc} p_{{j_i},l} & t_{i,l} \\
p_{{j_i},m} & t_{i,m}
\end{array}\right) \, ,\;
D_{i,2}=\det \left( \begin{array}{cc} x_{i,l} & p_{{j_i},l} \\
x_{i,m} & p_{{j_i},m}
\end{array}\right) \, ,\;
\ee
$l \ne m \in \{0,1,2\}$.\\
By Theorem \ref{thrmcr} the point
$U_{i}$ is thus obtained as
\[
U_{i} = D_{i,1} X_i - D_{i,2} T_i.
\]

\section{Modification of the scheme for non-convex data}
\label{nonconvex}

If the input data are not convex, we segment them according to the
following criteria in order to obtain piecewise convex segments.\\
First, consecutive collinear points are identified in the
following way. \be\label{critcol} \begin{minipage}{11cm}
By applying a dominant points selection algorithm like, for instance, the one in \cite{MS03},
we can easily detect the end points of a sequence containing at least 3 collinear vertices.
Let us denote them by $\bp_j^0$ and $\bp_l^0$.
\end{minipage}  \ee
Since in CAD applications as we have in mind linear features are
usually intentional, we do not smooth the angle between a subpolygon
consisting of (at least 3) collinear points and its neighbors. The
insertion rule for these straight line subpolygons between the
points $\bp_j^k$ and $\bp_l^k$ ($k=0, 1, 2, 3, \ldots$) simply
reads as: \be \label{insertline} \bp_{2i+1}^{k+1} = \frac{1}{2}
(\bp_{i}^{k} + \bp_{i+1}^{k}) \, , \; i=j, \ldots, l-1. \ee
The remaining subpolygons do not anymore contain collinear segments. For
these remaining subpolygons, inflection edges are identified by the
following criterion, see also \cite{cinesi,vicky}.
\be\label{inflcrit} \begin{minipage}{11cm}
An edge $\bp_i^0
\bp_{i+1}^0$ is identified as {\it inflection edge} if the points
$\bp_{i-1}^0$ and $\bp_{i+2}^0$ lie in different half planes with
respect to it.
\end{minipage}  \ee
On an inflection edge we insert a new point, e.g., as midpoint of
the edge corners, and we thus have a sequence of subpolygons without
inflections. Next, we check each of these subpolygons $\bp_{j}^0
\ldots \bp_{l}^0$ for total convexity by the following criterion:
\be\label{convcrit} \begin{minipage}{11cm} If for every edge of the
subpolygon all the points of the subpolygon lie either on the edge
or in the same half plane with respect to the edge, then the
subpolygon is totally convex, otherwise it is only locally convex.
\end{minipage}  \ee
If a subpolygon is only locally convex we divide it into two new
subpolygons $\bp_{j}^0 \ldots \bp_{i}^0$ and $\bp_{i}^0 \ldots
\bp_{l}^0$, where $i = j + \llcorner\frac{l-j+1}{2} \lrcorner$. We repeat
this procedure until we have a sequence of totally convex
subpolygons, each of which we suppose to be composed of at least
five points.

Whereas in the interior of every convex subpolygon we apply the
algorithm detailed in the previous section, we now define the method
next to
\begin{enumerate}
\item
an inserted junction point on an inflection edge, which we refer to
as {\it inflection point},
\item
a junction point between two convex subpolygons, which we refer to
as {\it convex junction point}.
\end{enumerate}

In the case of an {\it inflection point} $\bp_i^0$ let us denote the
line of the inflection edge $\bp_{i-1}^0 \bp_i^0 \bp_{i+1}^0$ by
$\eee_i$, and the convex subpolygons meeting in $\bp_i^0$ by
$\bs_{l,i}^0$ and $\bs_{r,i}^0$. We estimate a left and a right
tangent in $\bp_i^0$, $\bl_{l,i}^0$ and $\bl_{r,i}^0$ respectively,
by applying the tangent estimation method from \cite{ABFH08} to the
five points $\bp_{i-4}^0, \ldots, \bp_{i-1}^0, \bp_{i}^0$ of polygon
$\bs_{l,i}^0$, respectively $\bp_{i}^0, \ldots, \bp_{i+3}^0,
\bp_{i+4}^0$ of polygon $\bs_{r,i}^0$.

\medskip
\begin{figure}[ht!]
\centering
\resizebox{5.0cm}{!}{\includegraphics{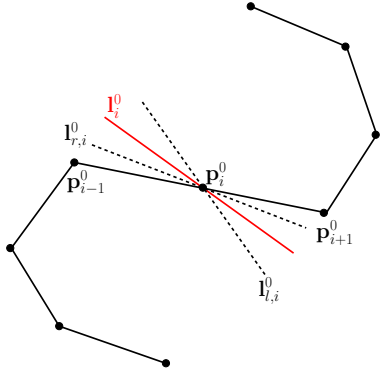}}
\caption{Definition of an initial tangent $\bl_i^0$ in an inflection
point $\bp_i^0$.} \label{figl0_inflpoint}
\end{figure}
\medskip

\begin{figure}[ht!]
\centering
\resizebox{6.5cm}{!}{\includegraphics{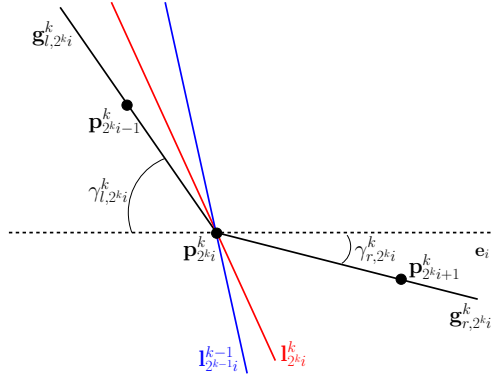}}
\caption{Tangent definition in an inflection point.}
\label{figl2ki_inflpoint}
\end{figure}

\noindent
We then combine these two lines $\bl_{l,i}^0$ and $\bl_{r,i}^0$ for defining an initial
tangent $\bl_i^0$ in $\bp_i^0$ (see Figure \ref{figl0_inflpoint}):
\be \label{li0} \bl_i^0 = \lambda_i^0 \bl_{l,i}^0 + \mu_i^0
\bl_{r,i}^0\, , \ee where $\lambda_i^0 + \mu_i^0 = 1$, $\lambda_i^0,
\mu_i^0 > 0$. The pair $(\lambda_i^0, \mu_i^0)$ thus plays the role
of a shape parameter. The tangents in the other vertices of the
polygons $\bs_{l,i}^0$ and $\bs_{r,i}^0$, and thus their new vertices,
are calculated as described in the previous section by treating
$\bs_{l,i}^0$ and $\bs_{r,i}^0$ separately. We obtain the new
polygons $\bs_{l,i}^1$ and $\bs_{r,i}^1$. Let us now describe how to
obtain the tangents $\bl_{2^k i}^k$ in the point $\bp_i^0= \bp_{2^k
i}^k$ in the following iterations ($k=1, 2, 3, \ldots$), see Figure
\ref{figl2ki_inflpoint} for an illustration.

\medskip
\medskip

\noindent
Let \be \label{glrk}
\bg_{l,2^k i}^k = \bp_{2^k i-1}^k \bp_{2^k i}^k\; \, \mbox{and} \;
\; \bg_{r,2^k i}^k = \bp_{2^k i}^k \bp_{2^k i+1}^k \ee be the edges
that are incident in $\bp_i^0= \bp_{2^k i}^k$, and \be
\label{gammalrk} \gamma_{r,2^k i}^k = \angle(\bg_{r,2^k i}^k,
\eee_i) \; \, \mbox{and} \; \; \gamma_{l,2^k i}^k =
\angle(\bg_{l,2^k i}^k, \eee_i)\ee their respective angles with the
initial inflection edge $\eee_i$.
We then define the line $\bg_{2^k i}^k$ by choosing that one of the
lines $\bg_{l,2^k i}^k$ and $\bg_{r,2^k i}^k$ from (\ref{glrk})
yielding the maximum angle \be \label{gammamax} \gamma_{2^k i}^k =
\max\{\gamma_{l,2^k i}^k, \gamma_{r,2^k i}^k\}\,. \ee The tangent
$\bl_{2^k i}^k$ in the point $\bp_i^0= \bp_{2^k i}^k$ is then
defined as \be \label{lnew} \bl_{2^k i}^k = \lambda_{2^k i}^k
\bl_{2^{k-1} i}^{k-1} + \mu_{2^k i}^k \bg_{2^k i}^k\, , \ee where
$\lambda_{2^k i}^k + \mu_{2^k i}^k = 1$ and $\lambda_{2^k i}^k,
\mu_{2^k i}^k > 0$. In the other vertices of $\bs_{l,i}^k$ and
$\bs_{r,i}^k$ we proceed as in section \ref{def_scheme} for
estimating the tangents; this allows us to calculate the new
polygons $\bs_{l,i}^{k+1}$ and $\bs_{r,i}^{k+1}$ by separately
applying the ``convex'' procedure from the previous section.

\begin{figure}[ht!]
\centering
\resizebox{7.0cm}{!}{\includegraphics{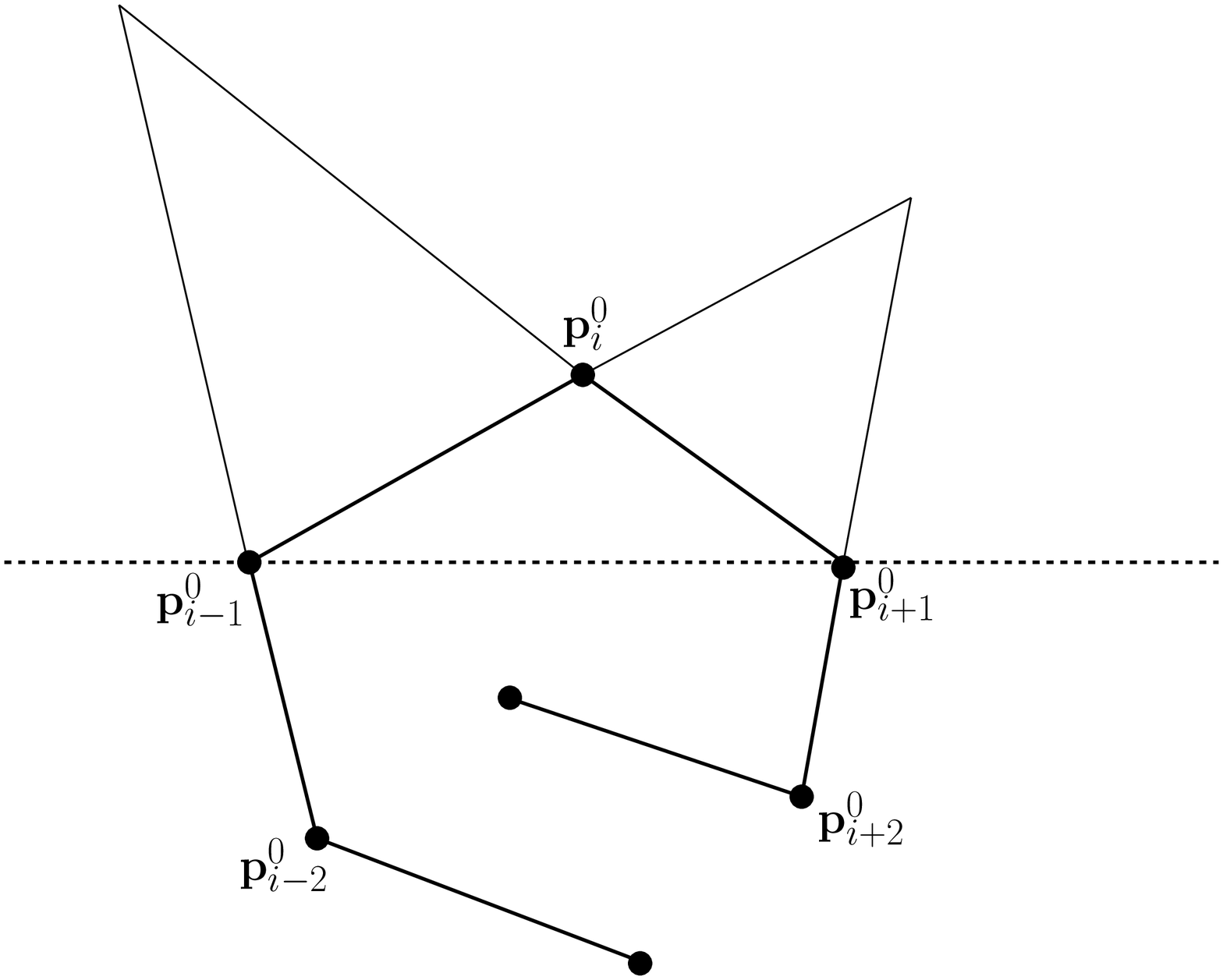}}
\caption{Situation in a convex junction point.}
\label{fig_condconvjoint}
\end{figure}

In the case of a {\it convex junction point} $\bp_i^0$  we suppose
our data to be such that the following condition holds\footnote{This
condition is guaranteed by a sufficiently dense sampling of the
initial data points.} (see Figure \ref{fig_condconvjoint}): \be
\label{condconvjoint} \mbox{\begin{minipage}{10cm} the intersection
points $\bp_{i-2}^0 \bp_{i-1}^0 \cap \bp_{i}^0 \bp_{i+1}^0$ and
$\bp_{i-1}^0 \bp_{i}^0 \cap \bp_{i+1}^0 \bp_{i+2}^0$ lie in the same
half plane with respect to the line $\bp_{i-1}^0 \bp_{i+1}^0$ as the
point $\bp_{i}^0$. \end{minipage}} \ee

\medskip
\medskip

\begin{figure}[ht!]
\centering
\resizebox{6.4cm}{!}{\includegraphics{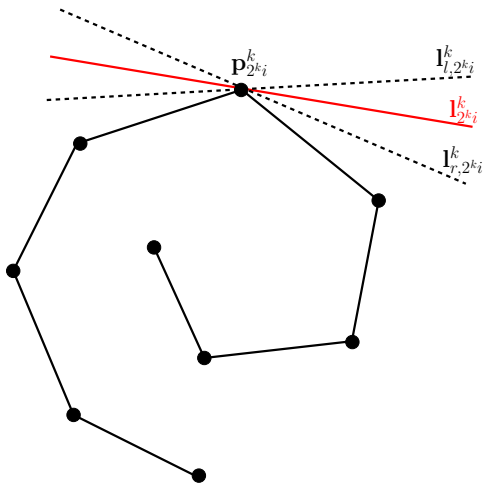}}
\caption{Tangent definition in a convex junction point.}
\label{figl2ki_convjoint}
\end{figure}

Let us denote the convex subpolygons meeting in $\bp_{i}^0$ by $\bs_{l,i}^0$ and
$\bs_{r,i}^0$. As in the case of an inflection point we estimate a
left and a right tangent in $\bp_i^0$, $\bl_{l,i}^0$ and
$\bl_{r,i}^0$ respectively, by applying the tangent estimation
method from \cite{ABFH08} to the five points $\bp_{i-4}^0, \ldots,
\bp_{i-1}^0, \bp_{i}^0$ of polygon $\bs_{l,i}^0$, respectively
$\bp_{i}^0, \ldots, \bp_{i+3}^0, \bp_{i+4}^0$ of polygon
$\bs_{r,i}^0$. If the points $\bp_{i-1}^0$ and $\bp_{i+1}^0$ lie in
different half planes with respect to $\bl_{l,i}^0$ ($\bl_{r,i}^0$
respectively) we replace $\bl_{l,i}^0$ ($\bl_{r,i}^0$ respectively)
by the line $\bp_i^0 \bp_{i+1}^0$ ($\bp_{i-1}^0 \bp_{i}^0$
respectively). We then combine these two lines $\bl_{l,i}^0$ and
$\bl_{r,i}^0$ as in (\ref{li0}) for defining an initial tangent
$\bl_i^0$ in $\bp_i^0$. The tangents in the other vertices of the
polygons $\bs_{l,i}^0$ and $\bs_{r,i}^0$, and thus their new vertices,
are calculated as described in the previous section by treating
$\bs_{l,i}^0$ and $\bs_{r,i}^0$ separately. We obtain the new
polygons $\bs_{l,i}^1$ and $\bs_{r,i}^1$. We then iterate this
procedure and obtain the tangents $\bl_{2^k i}^k$ in the point
$\bp_i^0 = \bp_{2^k i}^k$ in the following iterations ($k=1,2,3,
\ldots$) as \be \label{lconv} \bl_{2^k i}^k = \lambda_{2^k i}^k
\bl_{l,2^k i}^k + \mu_{2^k i}^k \bl_{r,2^k i}^k\, , \ee where
$\bl_{l,2^k i}^k $ and $\bl_{r,2^k i}^k$ is the respective left and
right tangent in $\bp_{2^k i}^k$ and $\lambda_{2^k i}^k + \mu_{2^k
i}^k = 1$, $\lambda_{2^k i}^k, \mu_{2^k i}^k > 0$ (see Figure
\ref{figl2ki_convjoint} for an illustration).

\medskip
\medskip

\begin{figure}[ht!]
\centering
\resizebox{6.0cm}{!}{\includegraphics{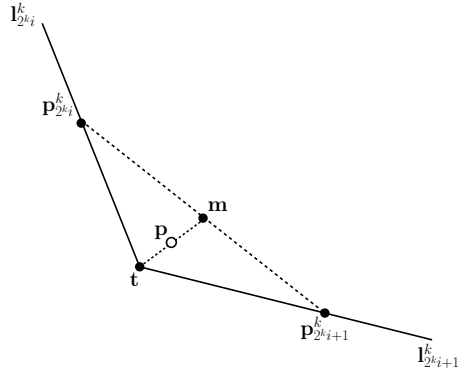}}
\caption{Modified point insertion rule next to a junction point.}
\label{fig_endpointrule}
\end{figure}

\begin{figure}[ht!]
\centering
\hspace{-0.5cm}
\subfigure[$k=1$]{\includegraphics[trim = 25mm 5mm 25mm 5mm, clip, width=4.0cm]{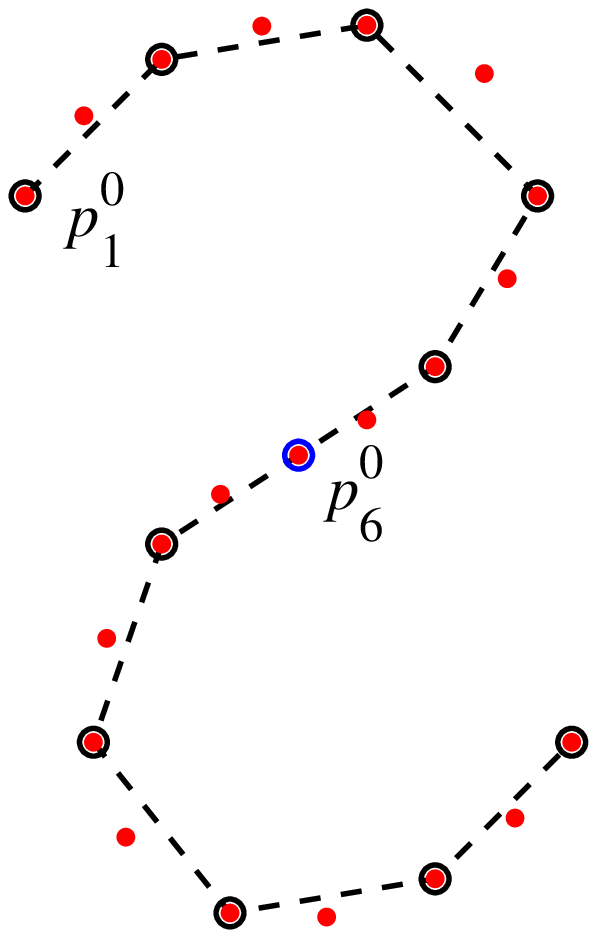}}\hspace{0.3cm}
\subfigure[$k=2$]{\includegraphics[trim = 25mm 5mm 25mm 5mm, clip, width=4.0cm]{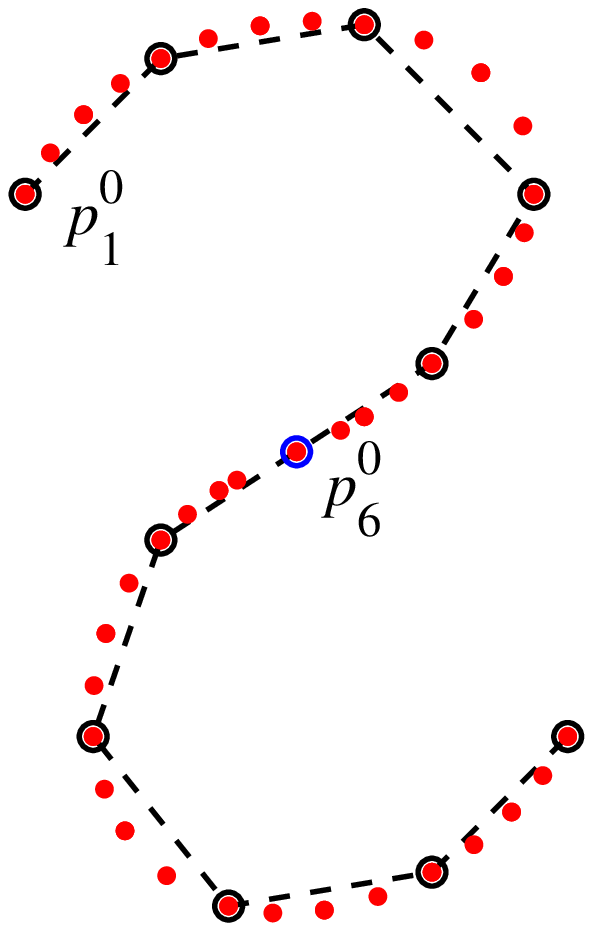}}\hspace{0.3cm}
\subfigure[$k=3$]{\includegraphics[trim = 25mm 5mm 25mm 5mm, clip, width=4.0cm]{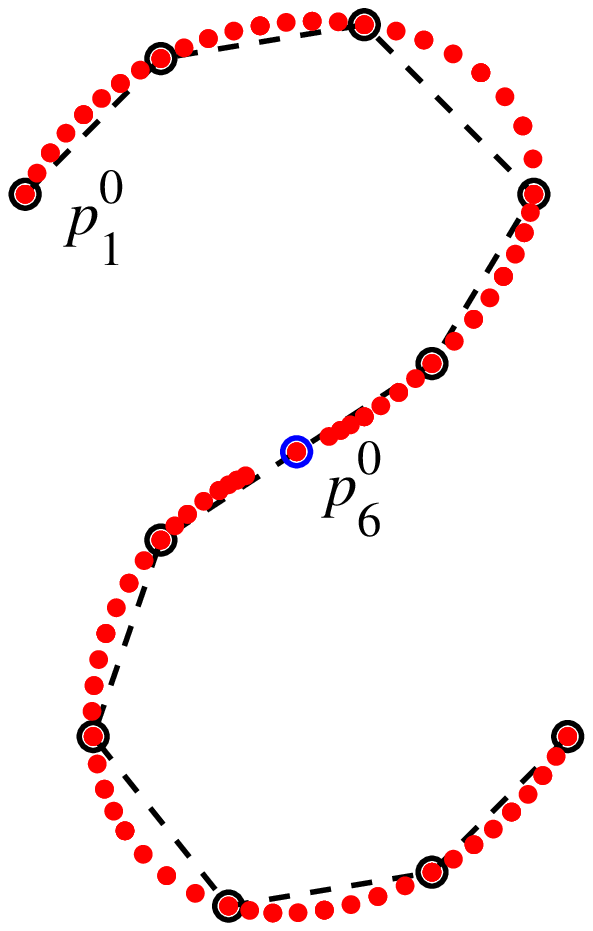}}
\caption{Example of slow convergence of the newly inserted
points towards the inflection point $p_6^0$.}
\label{slow_conv_example}
\end{figure}

In order to avoid a rather slow convergence of the newly inserted
points towards the junction points (inflection points or convex
junction points), see Figure \ref{slow_conv_example},
we modify the point insertion rule for the points $\bp_{2^{k+1} i-1}^{k+1}$ and
$\bp_{2^{k+1} i+1}^{k+1}$ in the tangent triangles adjacent to the junction point $\bp_{2^k
i}^k$ formed by the lines $\bl_{2^k i-1}^k$, $\bl_{2^k i}^k$ and the
edge $\bp_{2^{k} i-1}^{k} \bp_{2^k i}^k$, and $\bl_{2^k i}^k$,
$\bl_{2^k i+1}^k$ and the edge $\bp_{2^{k} i}^{k} \bp_{2^k i+1}^k$
respectively.
The modified rule for the first and the last new point
of a subpolygon reads as follows \be \label{newpointmod} \bp = \rho
\bt + \sigma \bm\, , \ \; \rho + \sigma =1\,, \ \; \rho, \sigma > 0\, ,
\ee where $\bp$ stands for the new point $\bp_{2^{k+1} i-1}^{k+1}$
respectively $\bp_{2^{k+1} i+1}^{k+1}$, $\bt$ designates the
intersection point of the two tangents in the corner points of the
corresponding edge and $\bm$ is the mid-point of this edge, see
Figure \ref{fig_endpointrule}.

\medskip
\medskip

\begin{figure}[ht!]
\centering
\hspace{-0.5cm}
\subfigure[$k=1$]{\includegraphics[trim = 25mm 5mm 25mm 5mm, clip, width=4.0cm]{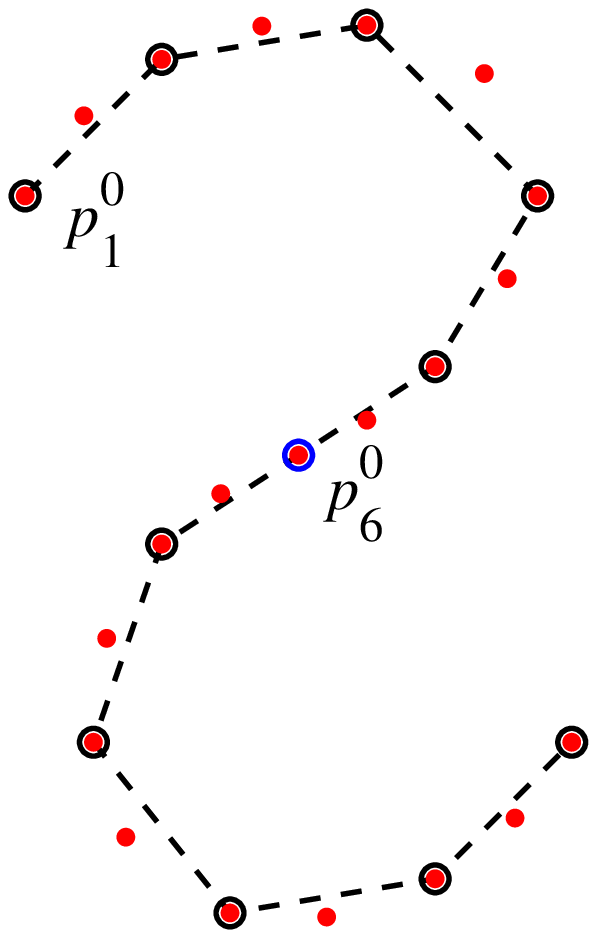}}\hspace{0.3cm}
\subfigure[$k=2$]{\includegraphics[trim = 25mm 5mm 25mm 5mm, clip, width=4.0cm]{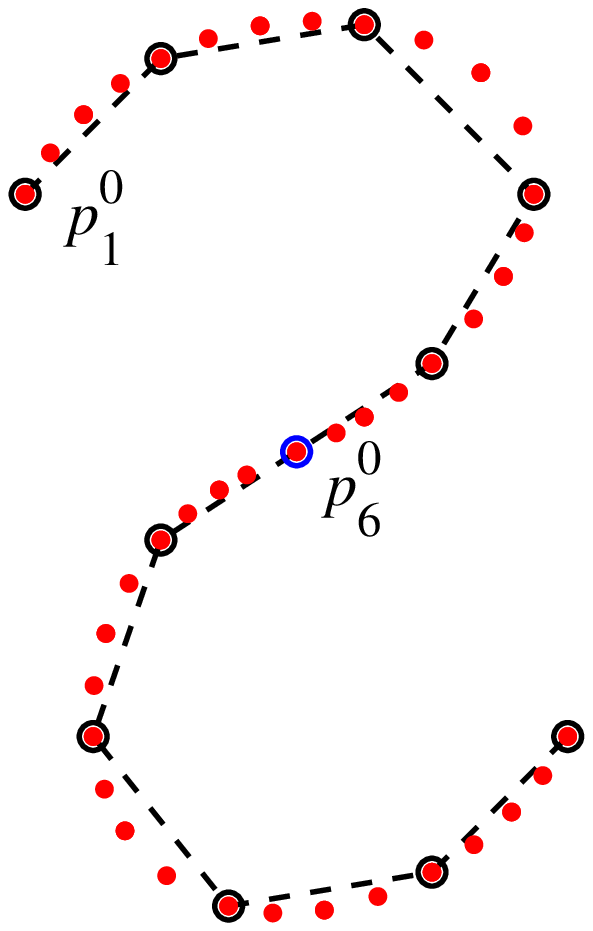}}\hspace{0.3cm}
\subfigure[$k=3$]{\includegraphics[trim = 25mm 5mm 25mm 5mm, clip, width=4.0cm]{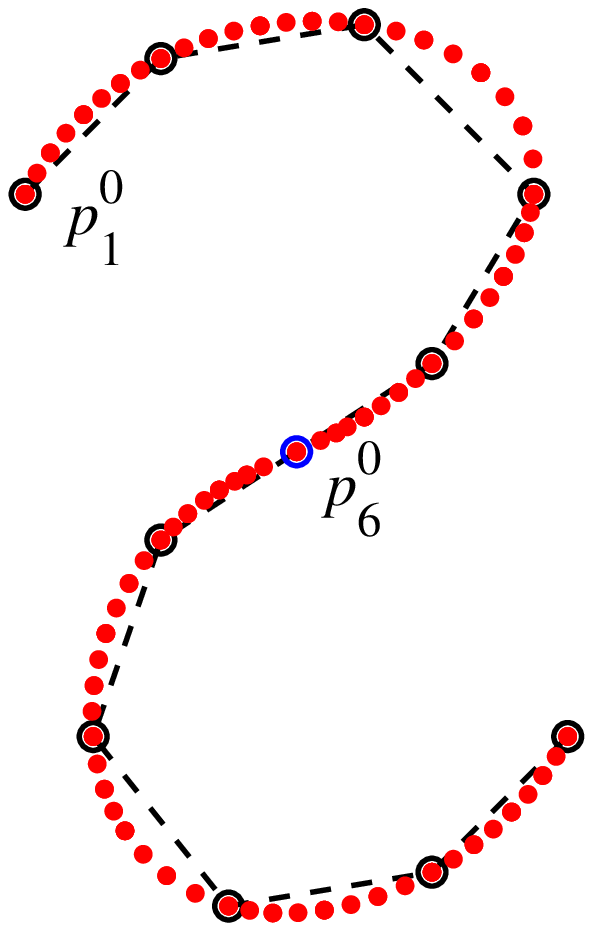}}
\caption{Application example of the new ``end point rule'' in equation (\ref{newpointmod}) to accelerate convergence of the newly inserted
points towards the inflection point $p_6^0$.}
\label{fast_conv_example}
\end{figure}

\noindent
This new ``end point rule'' avoids holes around the junction points (see Figure \ref{fast_conv_example}).

\medskip
\medskip
\medskip

\section{The subdivision algorithm}\label{algorithm}

The above detailed procedure may be summarized in the following
algorithm, which we apply to the given point sequence
$\bp^0=(\bp_i^0 \ : \ i \in \ZZ)$.
 For the sake of clarity we will first
describe the procedure to be used in the case of totally convex data
(Subsection \ref{total_conv_algo}), which constitutes an essential
ingredient of the final algorithm for general non-convex data
presented in Subsection \ref{non_conv_algo}.

\medskip

\subsection{Algorithm for totally-convex data}
\label{total_conv_algo}

Let $\bp^k=(\bp_i^k \ : \ i=1, \ldots,n_k)$ be the vertices of the
totally convex polygon at the $k$-level refinement. Hereafter we
denote by $P^{k}_{i}$ the projective counterparts of the affine
points $\bp^{k}_{i}$, see section \ref{def_scheme}.\\ The algorithm
that implements the function $\varphi$ from (\ref{def_phi}) proceeds
as follows in order to calculate the points
$\bp_{2i+1}^{k+1}=\varphi(\bp_{i-2}^k,\bp_{i-1}^k,\bp_i^k,\bp_{i+1}^k,\bp_{i+2}^k,
\bp_{i+3}^k \, ; \, \bp)$.

\smallskip
\begin{enumerate}
\item[\sl{Step 1:}] Preprocessing

\begin{enumerate}
\item[a)] For a closed polygon we set $P_0^k = P_{n_k}^k$, $P_{-1}^k = P_{n_k-1}^k$,
$P_{n_k+1}^k = P_1^k$, $P_{n_k+2}^k = P_2^k$.

For an open polygon we set $P_0^k = P_5^k$, $P_{-1}^k = P_{4}^k$,
$P_{n_k+1}^k = P_{n_k-4}^k$, $P_{n_k+2}^k = P_{n_k-3}^k$.

For $i=1, \ldots, n_k$ we then assign $Q_1=P_{i-2}^k$,
$Q_2=P_{i-1}^k$, $Q_3=P_{i}^k$, $Q_4=P_{i+1}^k$, $Q_5=P_{i+2}^k$,
and apply formula (\ref{est_tang}) for obtaining the tangent $L_i^k$
in $P_i^k$.

\item[b)]
According to the angle criterion, we calculate the angles $\alpha_j^i$ for $j=i+2,\ldots,n_k+i-1$ and
for $i = 1, \ldots, n_k$ in the case of a closed polygon, while for $i = 1,
\ldots, n_k-1$ in the case of an open polygon. We then select the value $\alpha_{j_i}^i$ satisfying condition (\ref{cond_ki_angle}).
\end{enumerate}

\medskip

\item[\sl{Step 2:}]

For a closed polygon we set $L_{n_k+1}^k = L_1^k$ and we consider $i
= 1, \ldots, n_k$, whereas for an open polygon we consider $i = 1,
\ldots, n_k-1$. We thus calculate the intersection points
\[
T_i^k = L_i^k \wedge L_{i+1}^k\, .
\]

\item[\sl{Step 3:}]

Calculate the lines (for $i = 1, \ldots, n_k$ for a closed polygon,
and for $i = 1, \ldots, n_k-1$ for an open polygon)
\[
N_i^k = P_i^k \wedge P_{i+1}^k\, , \; \Lambda_i^k = P_{j_i}^k \wedge
T_i^k\, ,
\]
as well as their intersection points
\[
X_i^k = N_i^k \wedge \Lambda_i^k.
\]

\item[\sl{Step 4:}]

Calculate the points $P_{2i+1}^{k+1}$ as
\[
P_{2i+1}^{k+1} = D_{i,1}^k X_i^k - D_{i,2}^k T_i^k
\]
with $D_{i,1}^k, D_{i,2}^k$ according to (\ref{Ds}).

\end{enumerate}

\bigskip

\subsection{Algorithm for non-convex data}
\label{non_conv_algo}

\begin{enumerate}
\item[\sl{Step 1:}] Preprocessing

We preprocess the data according to the criteria (\ref{critcol}),
(\ref{inflcrit}) and (\ref{convcrit}) thus identifying subpolygons
contained in a straight line and introducing and /or identifying
inflection vertices and convex junction vertices within the initial
data points of the remaining subpolygons yielding a sequence
composed of totally convex subpolygons and ``straight line''
subpolygons.

\medskip

\item[\sl{Step 2:}]

\begin{itemize}
\item
For a ``straight line'' subpolygon consisting of collinear segments we fix the tangents
at its end points equal to the straight line passing through them and we
apply everywhere the insertion rule (\ref{insertline}).
\item
For a totally convex subpolygon $\bp_j^0 \ldots \bp_l^0 = \ldots =
\bp_{2^k j}^k \ldots \bp_{2^k l}^k$ between the junction points
$\bp_j^0$ and $\bp_l^0$ (where $l \ge j+4$) we apply the algorithm
of Subsection \ref{total_conv_algo} in order to calculate the
tangents in the points $\bp_{2^k j+1}^k, \ldots, \bp_{2^k l-1}^k$
and the new points (with upper index $k+1$) between $\bp_{2^k
j+1}^k$ and $\bp_{2^k l-1}^k$. According to the type of the
subpolygon end points $\bp_{2^k j}^k$ and $\bp_{2^k l}^k$, i.e.,
inflection point or convex junction point, we apply the rule
(\ref{lnew}) for an inflection point and (\ref{lconv}) for a convex
junction point in order to calculate the respective tangent, and
apply rule (\ref{newpointmod}) for calculating the first
respectively last new point, i.e., $\bp_{2^{k+1} j+1}^{k+1}$
respectively $\bp_{2^{k+1} l-1}^{k+1}$.
If $\bp_{2^k j}^k$ or $\bp_{2^k l}^k$ is a junction point with a
straight line subpolygon, we fix the tangent at such an end point
equal to the straight line passing through it and we define the new point
in the adjacent triangle using the standard rule.
Finally, if $\bp_{2^k j}^k$ or $\bp_{2^k l}^k$ coincides with an end point of the whole sequence, we proceed
by computing the tangent at that location and the new point closest to it exactly
following the same procedure we described in the open polygon case of Subsection \ref{total_conv_algo}.
\end{itemize}
\end{enumerate}

\section{Adaptive version of the scheme}\label{secadapt}

\begin{figure}[ht!]
\centering
\hspace{-0.2cm}
\subfigure[]{\resizebox{4.3cm}{!}{\includegraphics{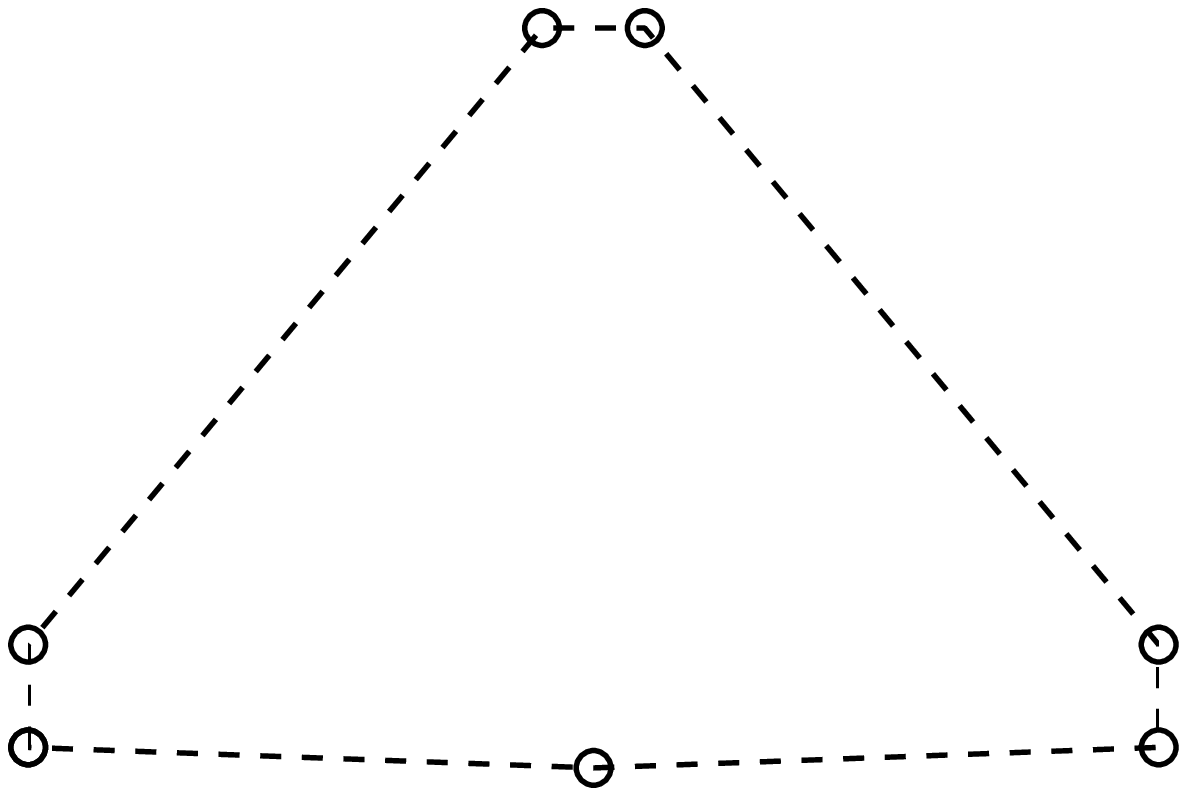}}}\hspace{-0.3cm}
\subfigure[]{\resizebox{4.3cm}{!}{\includegraphics{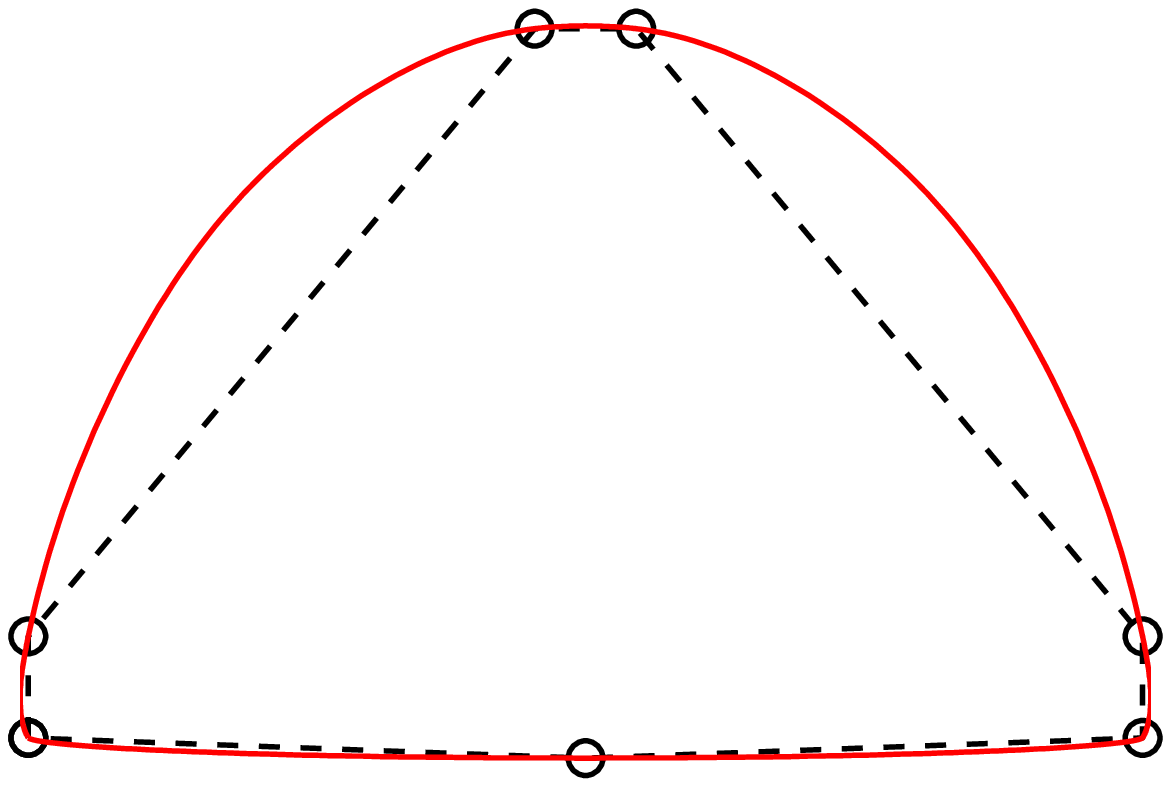}}}\hspace{-0.3cm}
\subfigure[]{\resizebox{4.6cm}{!}{\includegraphics{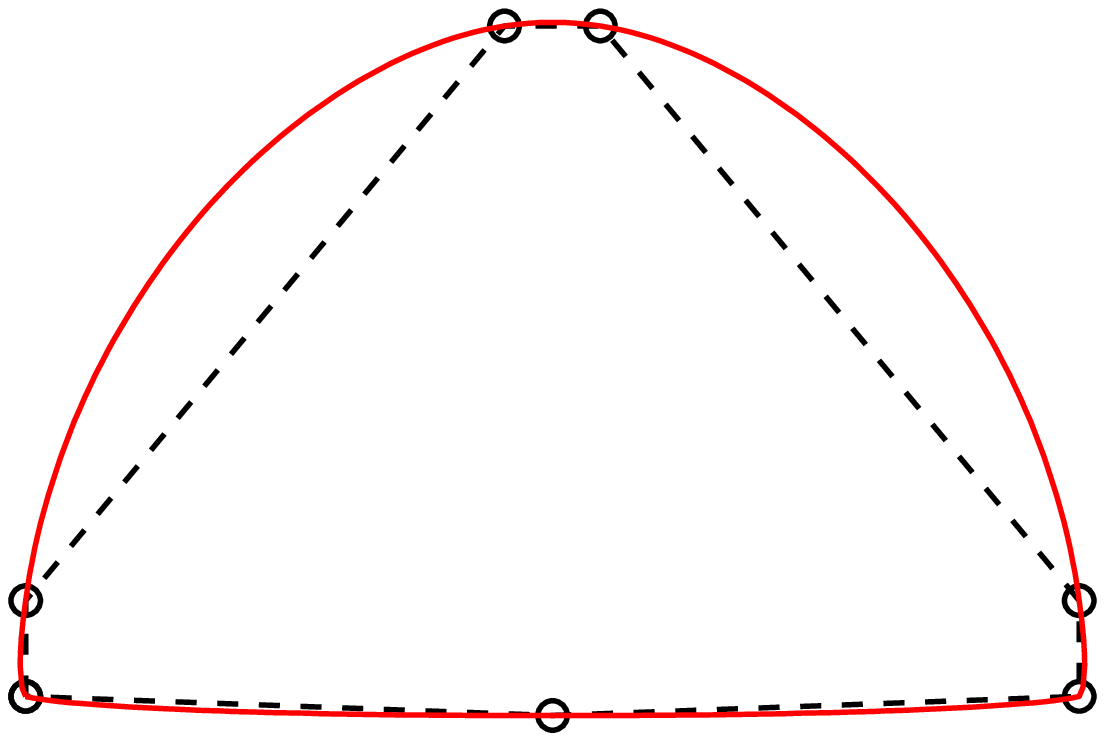}}}
\caption{Application example of the adaptive algorithm: (a) starting polyline; (b) refined polyline after 6 steps of the basic algorithm; (c) refined polyline after 6 steps of the adaptive algorithm.}
\label{adaptive_ref}
\end{figure}

Subdivision curves are visualized by drawing a polyline on a level of refinement which evokes the impression of sufficient approximation of the given data.
For a high quality rendering, the task is therefore to calculate and draw a level of subdivision which is a visually sufficient approximation of the limit shape.
The subdivision algorithm presented in Section \ref{algorithm} provides a process of global refinement at every level.
Therefore, when the starting polyline is highly non-uniform, the required level of refinement is determined by those locations which approximate the limit curve most unfavorably. Obviously, these may cause unnecessary fine subdivisions at other locations of the curve, thus leading to an unreasonable resource demanding algorithm.\\
To overcome this problem we propose an adaptive version of the subdivision algorithm previously described.
Adaptive subdivision is achieved by applying the mechanism of subdivision only locally, i.e. only at those locations of the initial polyline that are not approximated with the desired quality.
The decision where high resolution refinement is needed, strongly depends on the underlying application.
As concerns our algorithm, adaptivity may be controlled either by the user or by an automatic criterion.
In fact, the user may specify which portions of the polygon should be subdivided or the process may be automated by controlling whether the length of an edge is greater or not than a specified threshold.
Only in the positive case we insert a new point in correspondence of the considered edge.
Of course, besides improving the visual quality of the limit curve, the adaptive version of the scheme reduces the computational cost of the algorithm.\\
In the remainder of this section we take highly non-uniform
polylines and compare our adaptive refinement algorithm with the
basic one. Figure \ref{adaptive_ref} shows the comparison between
the refined polylines obtained by applying the two algorithms to a
totally convex set of points, while Figure \ref{adaptive_ref2}
compares the two algorithms on a highly non-uniform polyline with
convex junction points and inflection points.

\medskip
\begin{figure}[ht!]
\centering
\subfigure[]{\resizebox{5.3cm}{!}{\includegraphics{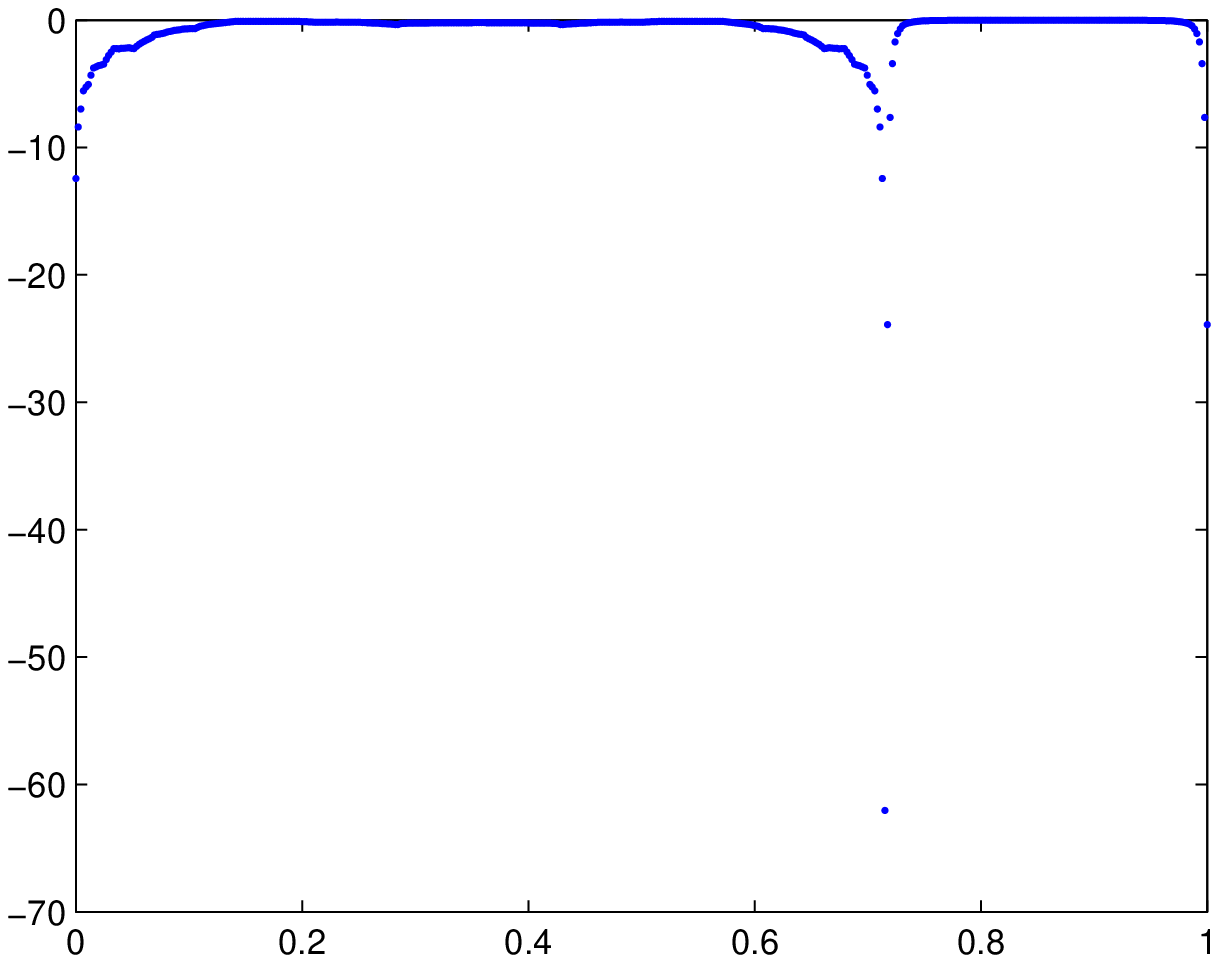}}}\hspace{-0.3cm}
\subfigure[]{\resizebox{5.3cm}{!}{\includegraphics{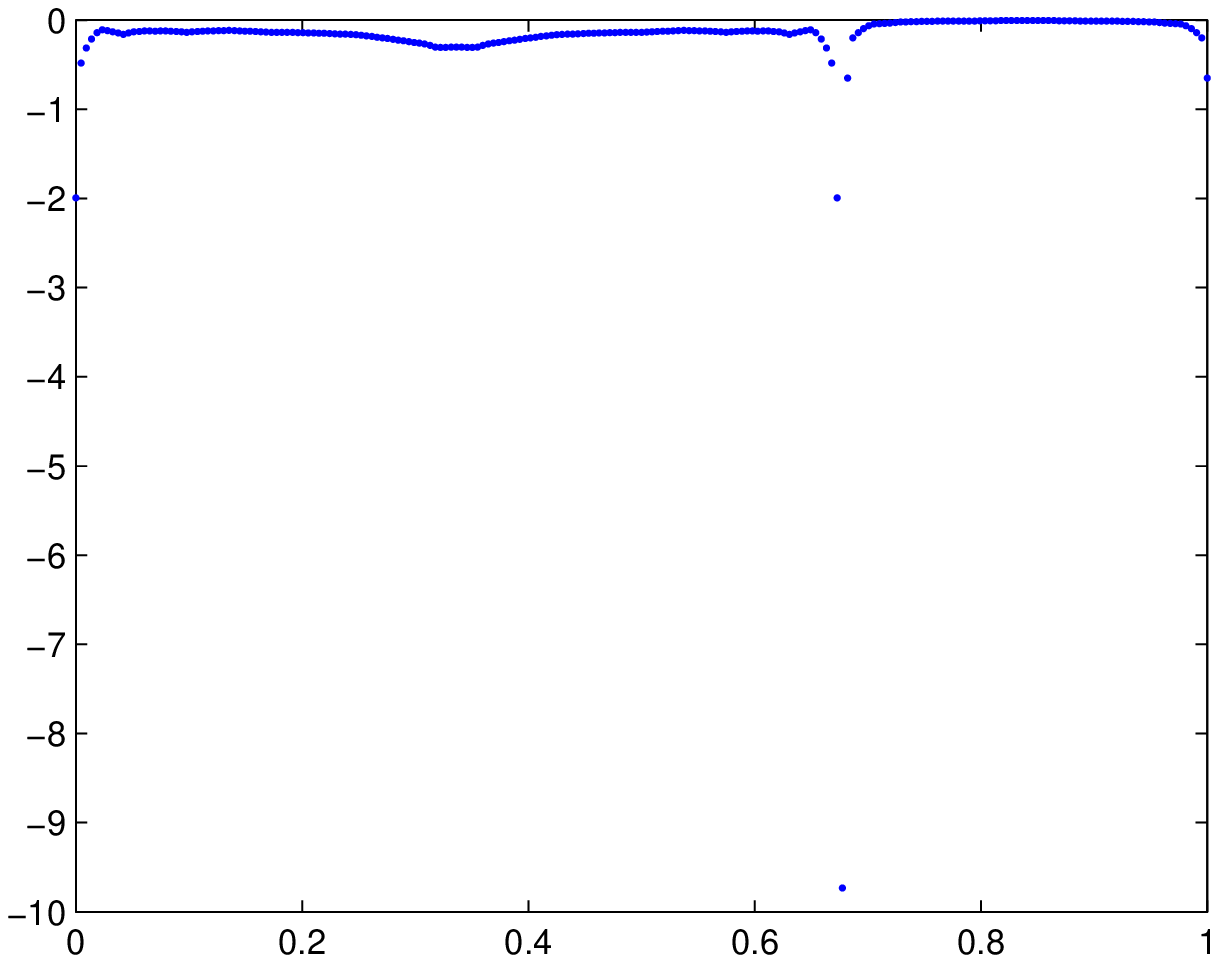}}}
\caption{Comparison between the discrete curvature plots of the refined polylines in Figure \ref{adaptive_ref} (b) and \ref{adaptive_ref} (c).}
\label{adaptive_ref_curv}
\end{figure}
\medskip

Figures \ref{adaptive_ref_curv} and \ref{adaptive_ref_curv2} illustrate the corresponding discrete curvature plots.
It can be easily seen that the use of the adaptive refinement scheme results in a considerable improvement of the curvature behaviour.

\medskip
\begin{figure}[ht!]
\centering
\hspace{-0.2cm}
\subfigure[]{\resizebox{4.3cm}{!}{\includegraphics{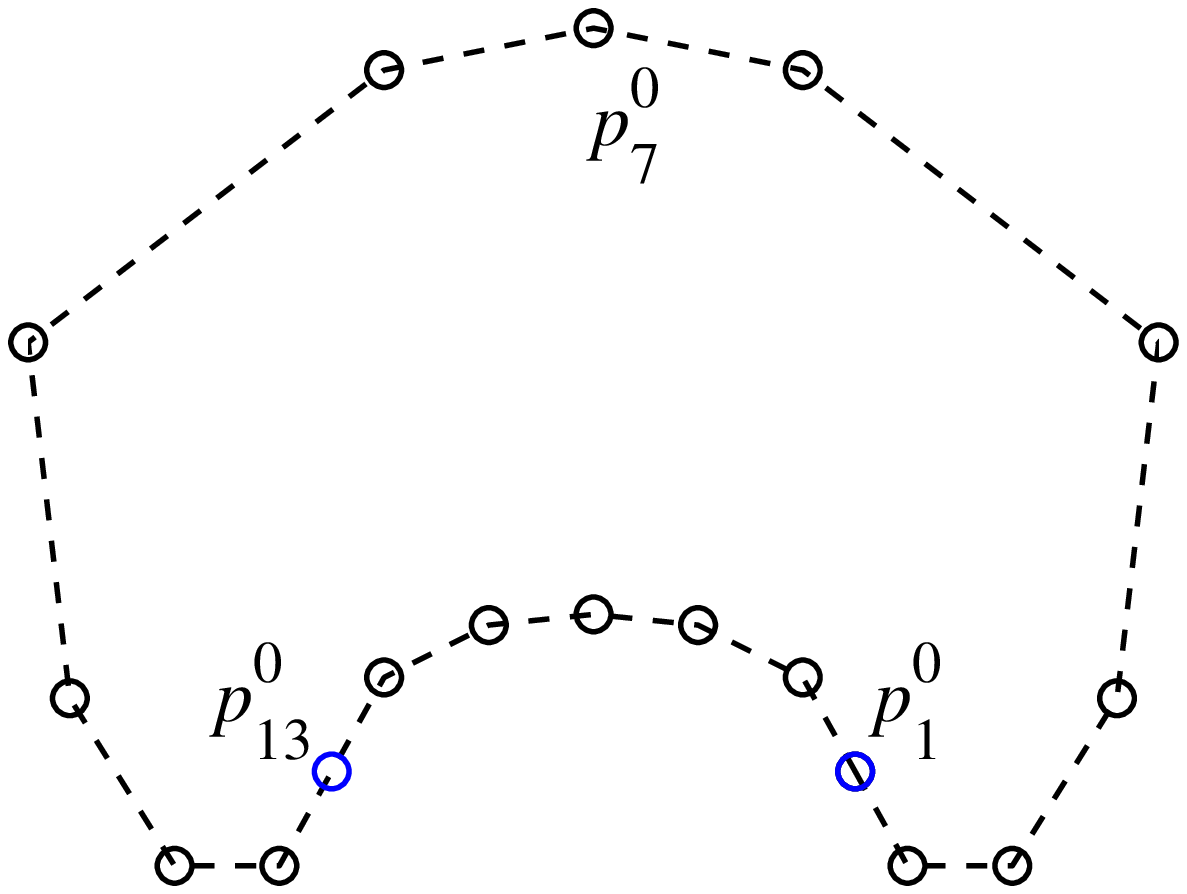}}}\hspace{-0.3cm}
\subfigure[]{\resizebox{4.3cm}{!}{\includegraphics{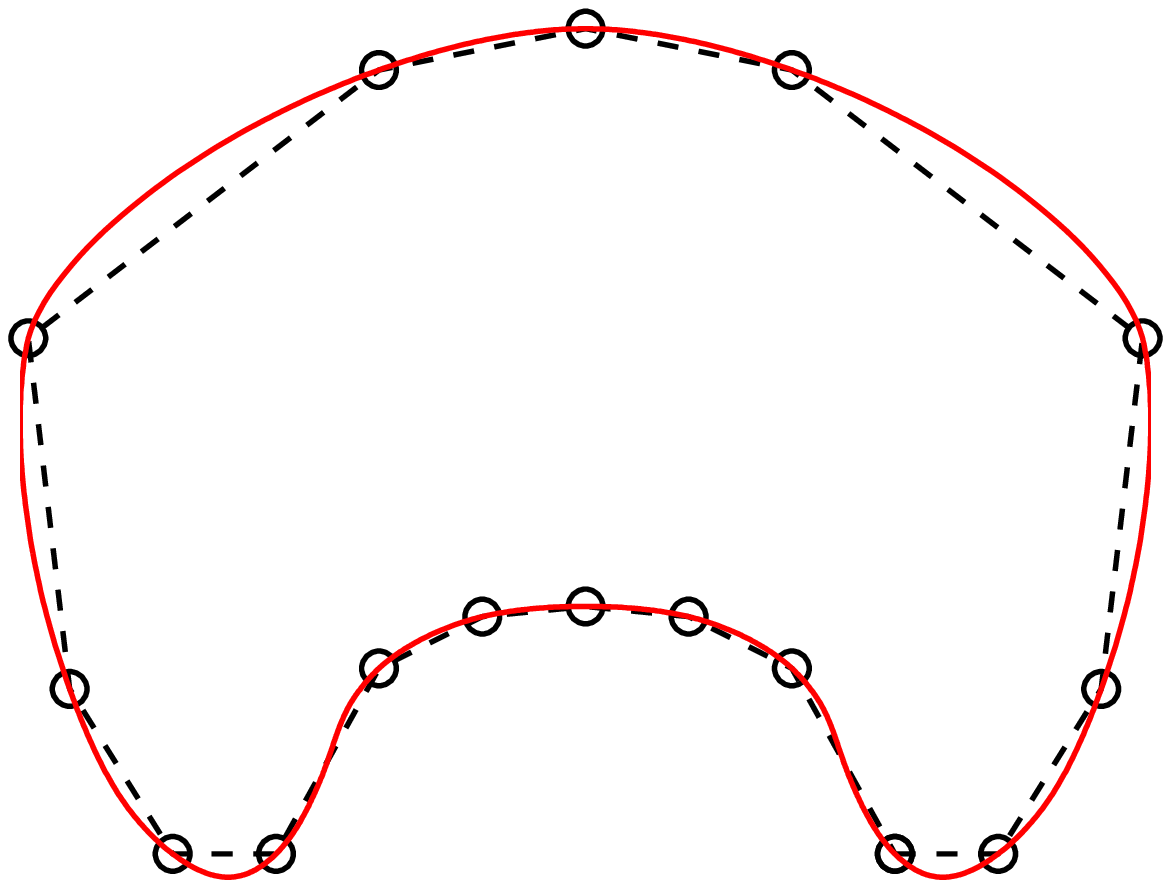}}}\hspace{-0.3cm}
\subfigure[]{\resizebox{4.7cm}{!}{\includegraphics{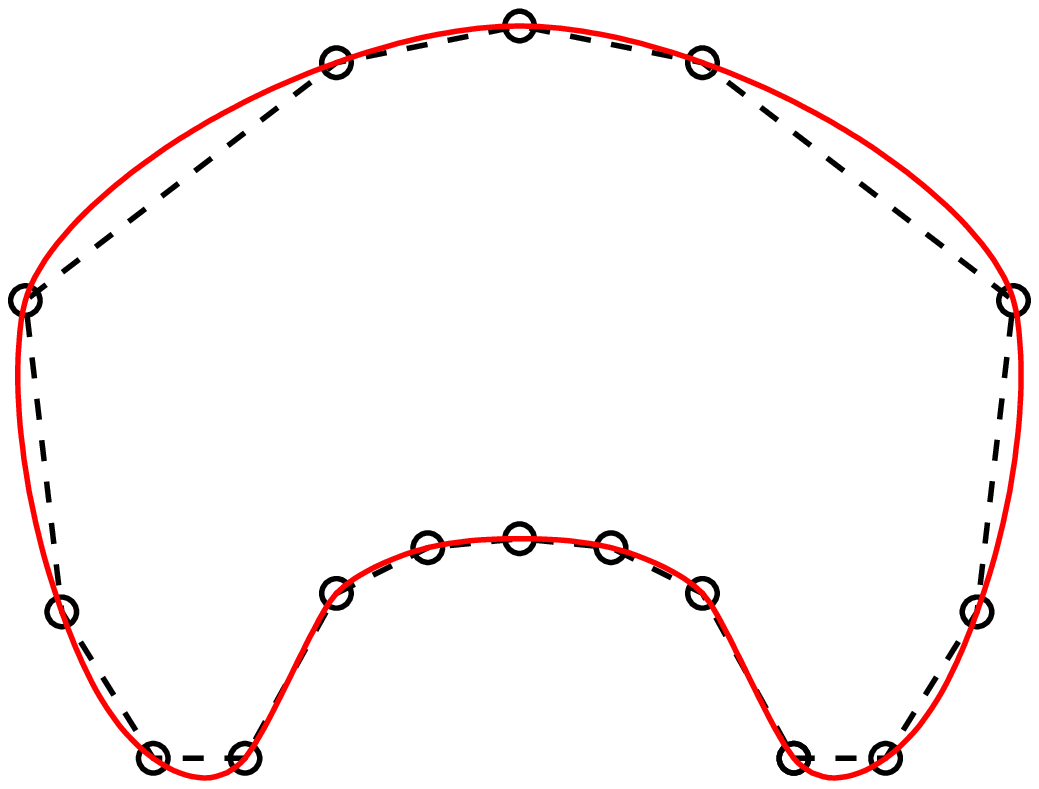}}}
\caption{Application example of the adaptive algorithm: (a) starting polyline; (b) refined polyline after 6 steps of the basic algorithm; (c) refined polyline after 6 steps of the adaptive algorithm. The labels in Subfigure (a) denote the end points of consecutive subpolygons: while $p_7^0$ is a convex junction point, $p_1^0$ and $p_{13}^0$ are inflection points.}
\label{adaptive_ref2}
\end{figure}
\medskip

\medskip
\begin{figure}[ht!]
\centering
\subfigure[]{\resizebox{5.3cm}{!}{\includegraphics{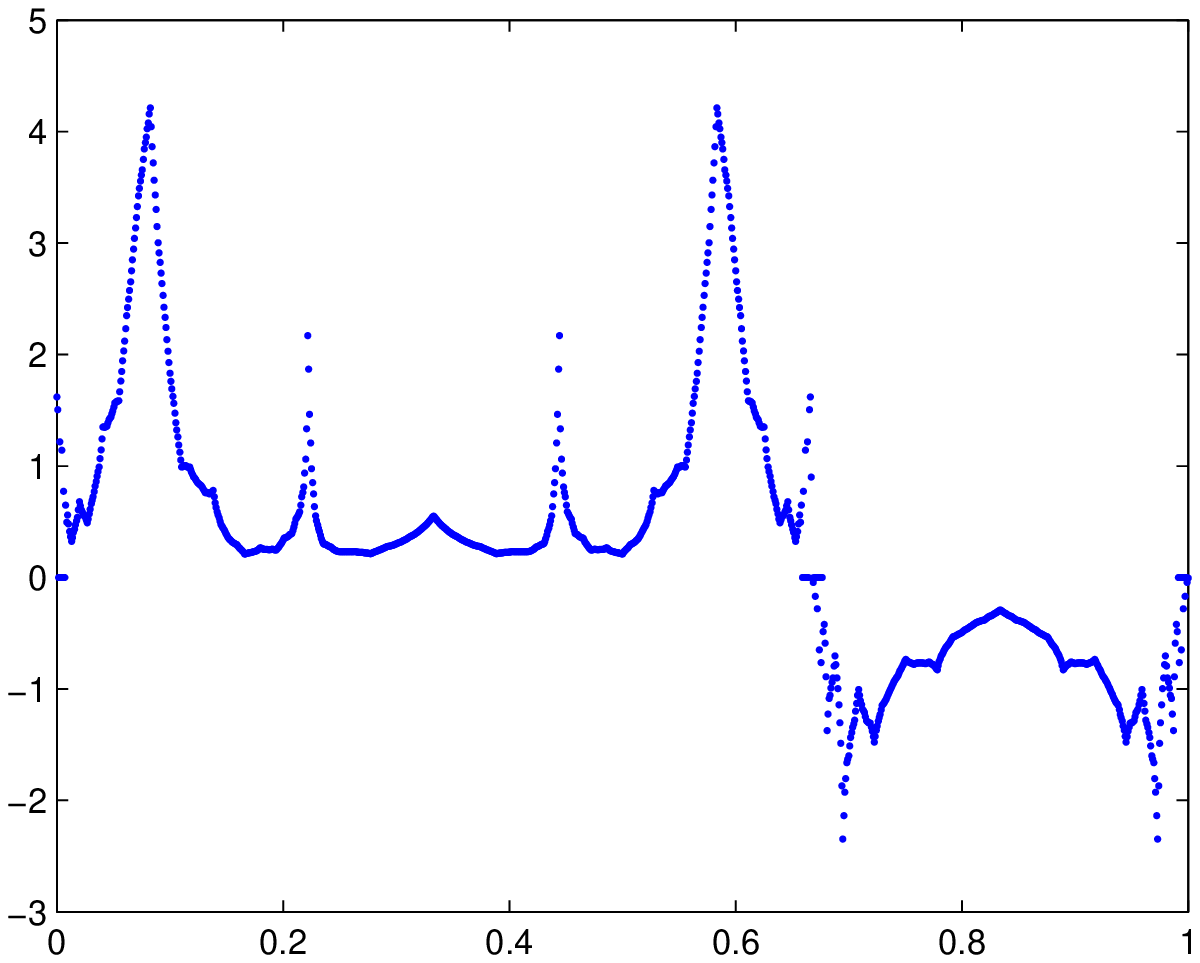}}}\hspace{-0.3cm}
\subfigure[]{\resizebox{5.3cm}{!}{\includegraphics{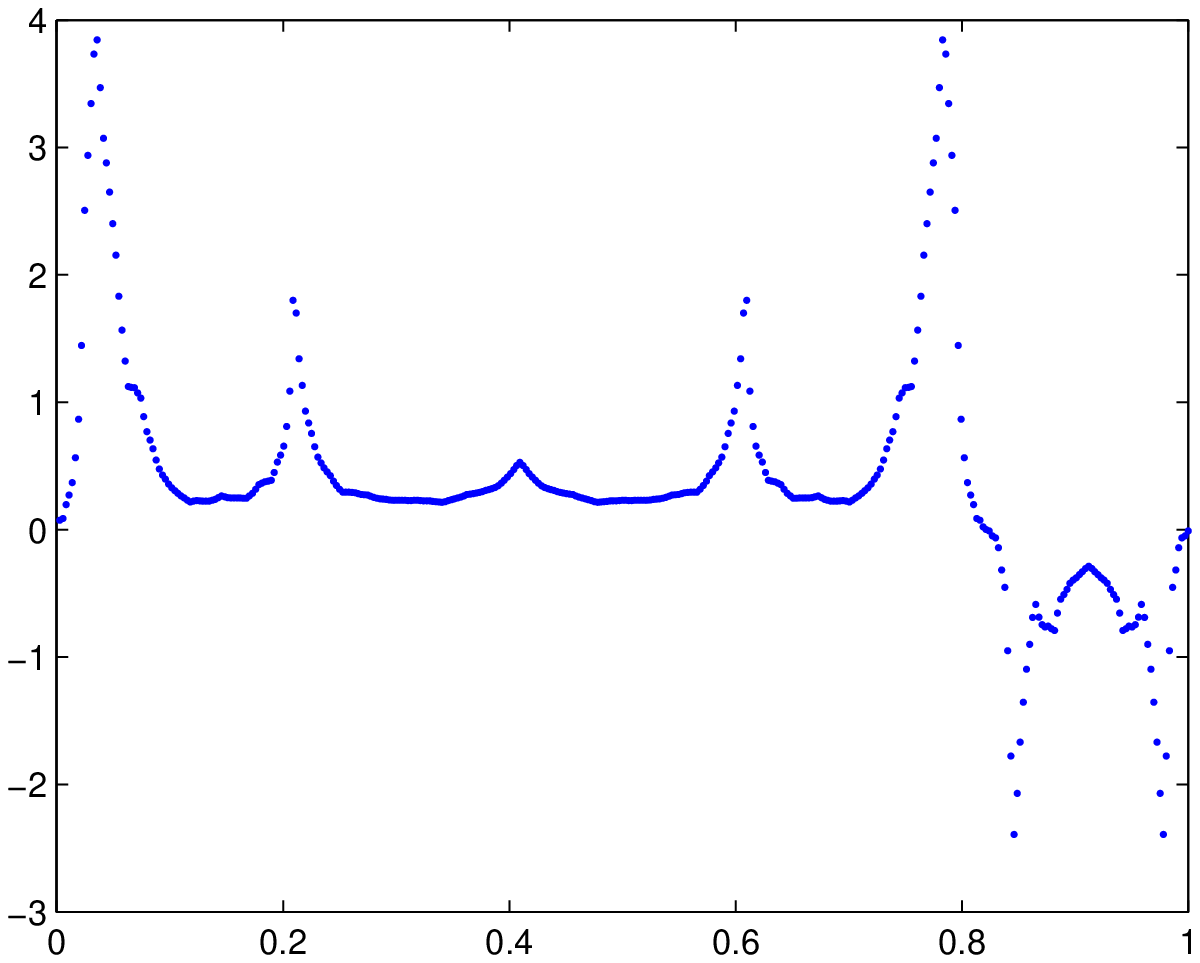}}}
\caption{Comparison between the discrete curvature plots of the refined polylines in Figure \ref{adaptive_ref2} (b) and \ref{adaptive_ref2} (c).}
\label{adaptive_ref_curv2}
\end{figure}
\medskip

\section{Properties of the scheme and smoothness analysis}\label{properties}

By construction the presented scheme enjoys the properties of being
shape preserving and conic reproducing as summarized in the
following proposition.

\smallskip
\begin{proposition} \label{prop}
The presented interpolatory curve subdivision scheme
\begin{enumerate}
\item[a)]
is shape preserving, i.e., if the starting point sequence
$\bp^0=(\bp_i^0 \ : \ i=1, \ldots,n_0)$ consists of convex, straight
line and concave segments, then all generated subsequent point
sequences $\bp^k=(\bp_i^k \ : \ i=1, \ldots,n_k)$, for
$k=1,2,3,\ldots$ respect the same behavior;
\item[b)]
reproduces conic sections, i.e., if the starting point sequence
$\bp^0$ is sampled from a conic section $c$, then the limit curve
coincides with $c$.
\end{enumerate}
\end{proposition}

\smallskip
\proof
\begin{enumerate}
\item[{\it a)}]
In the preprocessing step the initial point sequence is segmented
into straight line segments and totally convex segments, which are
not changed during the subdivision procedure. A straight line
subpolygon is reproduced as such. For a totally convex subpolygon
$\bp^k=(\bp_i^k \ : \ i=j_k, \ldots,l_k)$ the algorithm first
generates a line $\ll_i^k$ in every point $\bp_i^k$. By construction
(see Section \ref{def_scheme} and (\ref{lnew}), (\ref{lconv})) these
lines form a convex delimiting polygon for the points of the
sequence $\bp^{k+1}$ of the next subdivision level, see Figure
\ref{fig1b}. The fact that by construction the point
$\bp_{2i+1}^{k+1}$ is contained in the triangle formed by the points
$\bp_i^k, \tt_i^k, \bp_{i+1}^k$ ($\bp_{2i+1}^{k+1} \in
\Delta(\bp_i^k, \tt_i^k, \bp_{i+1}^k)$) guarantees convexity of the
sequence $\bp^{k+1}$.
\item[{\it b)}]
If the point sequence $\bp^k$ comes from a conic section $c$, then
the lines $\ll_i^k$ generated in the preprocessing step are the
tangents of $c$ in the points $\bp_i^k$. By Theorem \ref{thrmcr} the
constructed points $\bp_{2i+1}^{k+1}$ lie on $c$.
\endproof
\end{enumerate}

\smallskip
The special set-up of the presented subdivision scheme allows to
obtain a result on the smoothness of the limit curve as formulated
in the following proposition for which we suppose our data not to
contain any subpolygons that consist of consecutive collinear
segments. The transitions from and to such subpolygons are purposely
$C^0$ by construction, and these ``straight line'' subpolygons are
exactly reproduced by the subdivision algorithm. We can thus
restrict our attention to polygons not containing any straight line
subpolygons.

\smallskip
\begin{proposition} \label{smoothness}
\ben
\item[a)]
For the presented subdivision scheme the polygon series $\bp^k$
converges to a continuous curve.
\item[b)]
The limit curve of the presented subdivision scheme is of continuity
class $G^1$. \een
\end{proposition}

\smallskip
\proof The proof is formulated for the non--adaptive version of the
subdivision algorithm where every polygon edge is replaced by two
new edges in every step. For the sake of simplicity we omit the
details for the adaptive case, where the indices change but not the
general idea. We proceed by demonstrating the $C^0$ and $G^1$
continuity in two steps. First, we consider (locally) convex
segments (including convex junction points), and then we treat
inflection points. The locally convex segments are composed of
several totally convex
segments joined together at convex junction points.\\

\begin{figure}[ht!]
\centering
\resizebox{9.5cm}{!}{\includegraphics{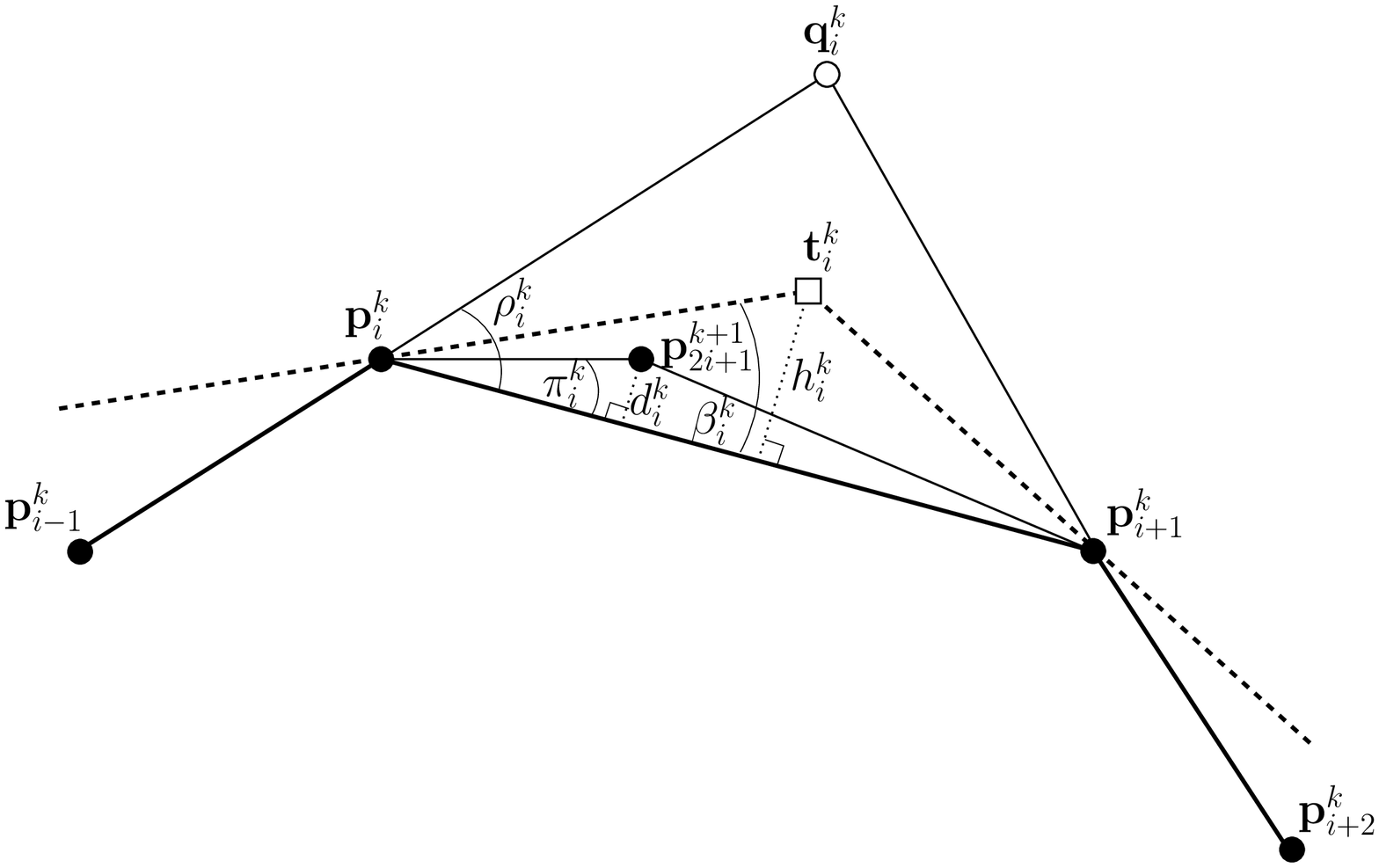}}
\caption{Notation used throughout the proof.} \label{fignot}
\end{figure}

\noindent We introduce the following notation. See Figure
\ref{fignot} for an illustration.

\medskip
\medskip

\begin{tabular}{cl}
$\bq_i^k$ & intersection point of the edges $\bp_{i-1}^k \bp_{i}^k$
and $\bp_{i+1}^k \bp_{i+2}^k$ \\
$\rho_i^k$ & inner angle in $\bp_{i}^k$ of the triangle
$\Delta(\bp_{i+1}^k \bp_{i}^k \bq_{i}^k)$ \\
$\beta_i^k$ & inner angle in $\bp_{i}^k$ of the triangle
$\Delta(\bp_{i+1}^k \bp_{i}^k \bt_{i}^k)$ \\
$\pi_i^k$ & inner angle in $\bp_{i}^k$ of the triangle
$\Delta(\bp_{i+1}^k \bp_{i}^k \bp_{2i+1}^{k+1})$ \\
$h_i^k$ & height of the triangle $\Delta(\bp_{i+1}^k \bp_{i}^k
\bt_{i}^k)$ from $\bt_{i}^k$
onto the edge $\bp_{i}^k \bp_{i+1}^k$\\
$d_i^k$ & height of the triangle $\Delta(\bp_{i+1}^k \bp_{i}^k
\bp_{2i+1}^{k+1})$ from $\bp_{2i+1}^{k+1}$
onto the edge $\bp_{i}^k \bp_{i+1}^k$\\
$c_i^k$ & length of the edge $\bp_{i}^k \bp_{i+1}^k$\\
$l_i^k$ & length of the line segment $\bp_{i}^k \bt_{i}^k$
\end{tabular} \\

\medskip
\noindent {\bf $C^0$ continuity for (locally) convex segments:}\\

\begin{figure}[ht!]
\centering
\resizebox{9.5cm}{!}{\includegraphics{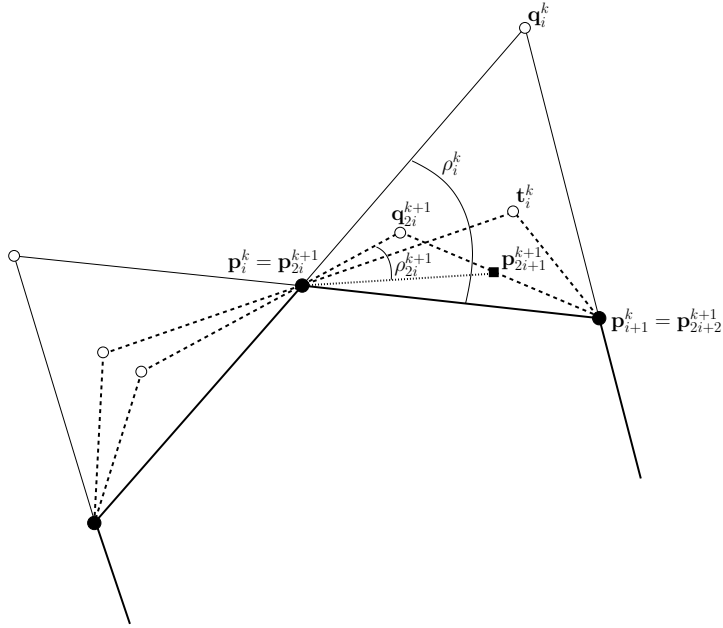}}
\caption{Relation between the angles $\rho_i^k$ and
$\rho_{2i}^{k+1}$.} \label{figrelrhos}
\end{figure}

\noindent
In order to show the continuity of the limit curve we compute the
distance $d^k$ between the polygon $\bp^{k+1}$ and the polygon
$\bp^{k}$. According to the above notation we define $d^k = \max_i
\{d_i^k\}$. By construction\footnote{In a convex junction point
this is guaranteed by condition (\ref{condconvjoint}) and the
respective tangent definition.} the point $\bt_{i}^k$ lies inside
the triangle $\Delta(\bp_{i+1}^k \bp_{i}^k \bq_{i}^k)$ and the point
$\bp_{2i+1}^{k+1}$ lies inside the triangle $\Delta(\bp_{i+1}^k
\bp_{i}^k \bt_{i}^k)$. Thus (see Figures \ref{fignot},
\ref{figrelrhos}) \be \label{pialpharho} 0 < \pi_i^k < \beta_i^k <
\rho_i^k \ee and \be \label{rho} \rho_{2i}^{k+1} < \rho_i^k \,.\ee

\medskip
\medskip

\noindent
Let $\bp_j^{k_0}$ be a point that appears in the $k_0$--th
iteration. Since the presented scheme is interpolatory we have
\[
\bp_j^{k_0} = \bp_{2j}^{k_0+1} = \bp_{4j}^{k_0+2} = \ldots =
\bp_{2^n j}^{k_0+n}\,.
\]
From (\ref{rho}) we obtain in the point $\bp_j^{k_0}$: \be
\label{ineqrho} \rho_{2^n j}^{k_0 +n} \le (\epsilon_j)^n
\rho_j^{k_0}\, ,   \ee and by (\ref{pialpharho}) we have \be
\label{ineqbetapi} \beta_{2^n j}^{k_0 +n}  \le  (\epsilon_j)^n
\rho_j^{k_0}\, , \;  \pi_{2^n j}^{k_0 +n}  \le  (\epsilon_j)^n
\rho_j^{k_0}\, , \ee where the $\epsilon_j$'s are constants such
that $0 < \epsilon_j < 1$. Let $\epsilon := \max_j\{ \epsilon_j\}$
and $\rho^{k_0} := \max_j\{ \rho_j^{k_0}\}$, where $0 < \epsilon <
1$. Then, \ba \rho_{2^n j}^{k_0
+ n} & \le & \epsilon^n \rho^{k_0}\, , \label{limr}\\
\beta_{2^n j}^{k_0
+ n} & \le & \epsilon^n \rho^{k_0}\, , \label{lima}\\
\pi_{2^n j}^{k_0 + n} & \le & \epsilon^n \rho^{k_0}\, . \label{limp}
\ea In the triangle $\Delta(\bp_{i+1}^k \bp_{i}^k \bp_{2i+1}^{k+1})$
we have
\[
\sin(\pi_i^k) = \frac{d_i^k}{c_{2i}^{k+1}}
\]
and thus \be \label{djkn} d_{2^n j}^{k_0+n} =
c_{2^{n+1}j}^{k_0+n+1} \sin(\pi_{2^n j}^{k_0+n}). \ee Since the
angles $\rho_i^k$ become smaller in every step and the longest edge
of a triangle is the one opposite to its biggest inner angle there
exists an index $k'$ such that for all $k \ge k'$ the circle
centered in $\bp_i^k$ with radius $c_i^k$ contains the triangle
$\Delta(\bp_{i+1}^k \bp_{i}^k \bq_{i}^k)$. Thus there exists an
index $n_0$ such that for all $n \ge n_0$ and for all $j$ we have
\[
c_{2^{n+1}j}^{k_0+n+1} \le c_{2^{n_0+1}j}^{k_0+n_0+1} \le
\max_j\{c_{2^{n_0+1}j}^{k_0+n_0+1}\} =: c^{k_0+n_0+1}.
\]
Since $\sin(x) \le x$ for $x \ge 0$, we deduce from (\ref{djkn}) and
(\ref{limp}) for $n \ge n_0$ that: \ba d_{2^n j}^{k_0+n} &
\le & c^{k_0+n_0+1} \pi_{2^n j}^{k_0+n} \nonumber \\
& \le & c^{k_0+n_0+1} \epsilon^n \rho^{k_0}. \nonumber \ea Since this
holds for all $j$ we have \be \label{dk} d^k \le c^{k_0+n_0+1}
\epsilon^{k-k_0} \rho^{k_0}\, , \ee where $k \ge k_0 + n_0$. The
polygons $\{\bp^k\}$ thus form a Cauchy sequence and this sequence
of polygons converges uniformly. Each polygon being piecewise
linear, the limit curve is continuous.\\

\medskip
\noindent {\bf $G^1$ continuity for (locally) convex segments:}\\

\begin{figure}[ht!]
\centering
\resizebox{9.8cm}{!}{\includegraphics{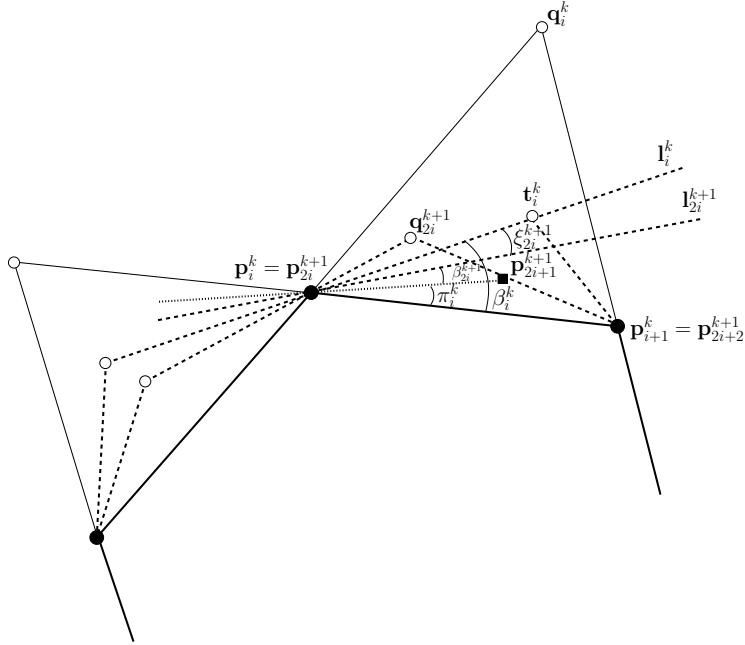}}
\caption{Illustration of the uniform convergence of the tangent
sequence.} \label{figproofG1}
\end{figure}

\noindent
Let
\[
h^k = \max_i\{h_i^k\}\, .
\]
For the triangle $\Delta(\bp_{i+1}^k \bp_{i}^k \bt_{i}^k)$ it holds
\[
\sin(\beta_i^k) = \frac{h_i^k}{l_i^k}\, .
\]
By the analogous reasoning as above for $d^k$ we obtain for $h^k$
(since there exists an index $k'$ such that for all $k \ge k'$: $
l_i^k \le c_i^k$):
\[
h^k \le c^{k_0+\bar{n}} \epsilon^{k-k_0} \rho^{k_0}\, ,
\]
 where $k \ge k_0 + \bar{n}$ for a certain value of $\bar{n}$ and
 \[
c^{k_0+\bar{n}} = \max_j\{c_{2^{\bar{n}} j}^{k_0 + \bar{n}}\}\, .
\]
In every iteration level the presented subdivision scheme constructs
the new points by sampling them from a $G^1$--continuous conic
spline $\Gamma^k$. By its Bézier construction the distance of the
conic segment of $\Gamma^k$ corresponding to the edge $\bp_i^k
\bp_{i+1}^k$ is bounded by $h_k$. The sequence of conic splines
$\{\Gamma^k\}$ thus converges to the same limit curve as the
sequence of polygons $\{\bp^k\}$.\\
It thus remains to be shown that the sequence of tangents
$\{\ll^k\}$, where $\ll^k = (\ll_i^k: i \in \ZZ)$ converges
uniformly. Let $\xi_{2i}^{k+1}$ be the angle between the tangents
$\ll_i^k$ and $\ll_{2i}^{k+1}$ in $\bp_i^k$. Then,
\[
| \xi_{2i}^{k+1} | = | \beta_{2i}^{k+1} + \pi_i^k - \beta_i^k|\, ,
\]
see Figure \ref{figproofG1}.

\medskip
\noindent
Thus by (\ref{lima}) and (\ref{limp})
with $i= 2^nj$ and $k= k_0 + n$: \ba | \xi_{2^{n+1}j}^{k_0+n+1}| &
\le & |\beta_{2^{n+1}j}^{k_0+n+1}| + |\pi_{2^{n}j}^{k_0+n}| +
|\beta_{2^{n}j}^{k_0+n}| \nonumber \\
& \le & \epsilon^n (\epsilon+2) |\rho^{k_0}|. \nonumber \ea The set
of tangents $\{\ll^k\}$ thus forms a Cauchy sequence and this
sequence converges uniformly. The limit curve is therefore of
continuity class $G^1$.\\

\medskip

\noindent {\bf $C^0$ and $G^1$ continuity for inflection points:}\\

\begin{figure}[ht!]
\centering
\resizebox{7.4cm}{!}{\includegraphics{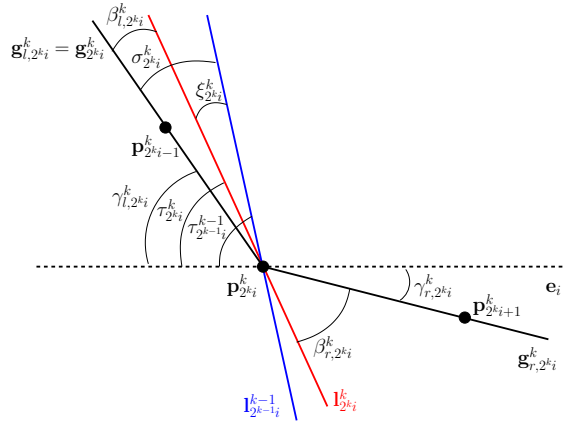}}
\caption{Illustration of $C^0$ and $G^1$ continuity for inflection
points.} \label{figinfl_angles}
\end{figure}

\noindent
In an inflection point $\bp_i^0 = \bp_{2^k i}^k$ we define \ba
\label{defangles} \sigma_{2^k i}^k & = &
\angle(\bl_{2^{k-1}i}^{k-1},
\bg_{2^k i}^k)\,, \nonumber \\
\tau_{2^k i}^k & = & \angle(\bl_{2^{k}i}^{k},
\eee_i)\,, \nonumber \\
\xi_{2^k i}^k & = & \angle(\bl_{2^{k-1}i}^{k-1},
\bl_{2^{k}i}^{k})\,, \nonumber \\
\beta_{l,2^k i}^k & = & \angle(\bl_{2^{k}i}^{k},
\bg_{l,2^k i}^k)\,, \nonumber \\
\beta_{r,2^k i}^k & = & \angle(\bl_{2^{k}i}^{k}, \bg_{r,2^k i}^k)
\nonumber\,, \ea where $\bl_{2^{j}i}^{j}$, $\bg_{2^k i}^k$,
$\bg_{l,2^k i}^k$, and $\bg_{r,2^k i}^k$ are defined as in section
\ref{nonconvex}, see Figure \ref{figinfl_angles}. We have \be
\label{stg} \sigma_{2^k i}^k = \tau_{2^{k-1}i}^{k-1} - \gamma_{2^k
i}^k\,, \ee with $\gamma_{2^k i}^k$ from (\ref{gammamax}).

\medskip
\medskip

\noindent By construction these angles satisfy the following
inequalities (see Figure \ref{figinfl_angles}): \ba 0 <
\xi_{2^k i}^k & < & \sigma_{2^k i}^k\,, \label{xs} \\
0 <
\tau_{2^k i}^k & < & \tau_{2^{k-1} i}^{k-1}\,, \label{tt} \\
0 <
\gamma_{2^k i}^k & < & \gamma_{2^{k+1} i}^{k+1}\,, \label{gg} \\
0 <
\beta_{l,2^{k+1} i}^{k+1} & < & \beta_{l,2^{k} i}^{k}\,, \label{bl} \\
0 < \beta_{r,2^{k+1} i}^{k+1} & < & \beta_{r,2^{k} i}^{k}\,.
\label{br} \ea

\noindent
The relations (\ref{bl}) and (\ref{br}) imply that
the tangent triangles adjacent to the inflection point $\bp_i^0$,
used for the determination of the new points next to $\bp_i^0$,
become continuously flatter, thus yielding $C^0$ continuity in
$\bp_i^0$.\\
Furthermore, we have \be \label{ss} \sigma_{2^k i}^k = \tau_{2^{k-1}
i}^{k-1} - \gamma_{2^k i}^k > \tau_{2^k i}^k - \gamma_{2^{k+1}
i}^{k+1} = \sigma_{2^{k+1} i}^{k+1} > 0\, . \ee By introducing an
$\epsilon$ with $0 < \epsilon < 1$ we thus have
\[ \sigma_{2^k i}^k \le \epsilon^{k-1} \sigma_{2i}^1
\]
and consequently
\[
0 < \xi_{2^k i}^k < \sigma_{2^k i}^k \le \epsilon^{k-1}
\sigma_{2i}^1\,.
\]
The angle between two consecutive tangents in an inflection point
$\bp_i^0$ thus continuously decreases yielding a $G^1$ continuous
inflection joint in the limit.
\endproof

\section{Numerical examples}\label{examples}
In this section we want to illustrate the performance of the basic and adaptive
subdivision algorithms presented in Sections \ref{algorithm} and \ref{secadapt}, respectively.
Generally, when the starting points are uniformly spaced we use the first proposal, while in the case of irregularly
distributed vertices we apply the second one.

\begin{figure}[h!]
\centering
\resizebox{2.4cm}{!}{\includegraphics{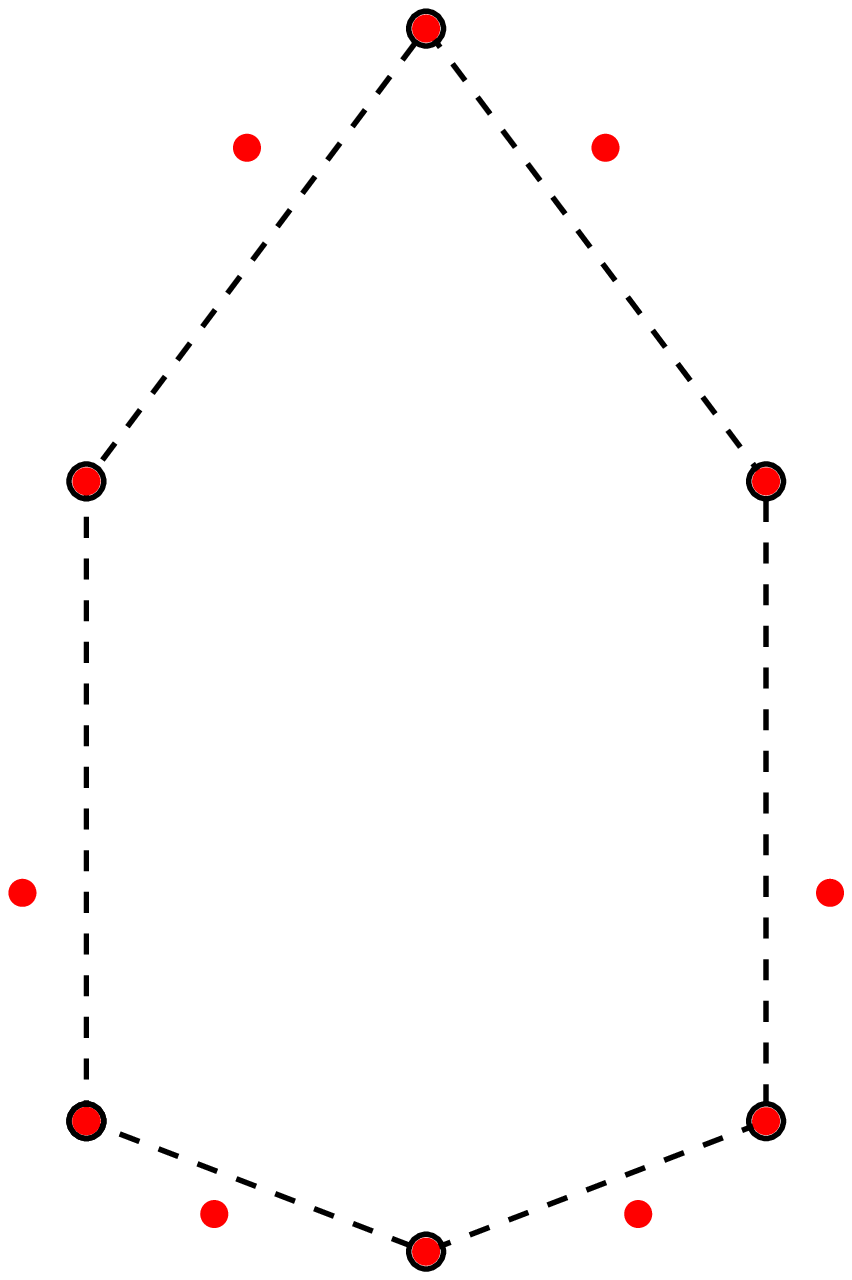}}\hspace{0.5cm}
\resizebox{2.4cm}{!}{\includegraphics{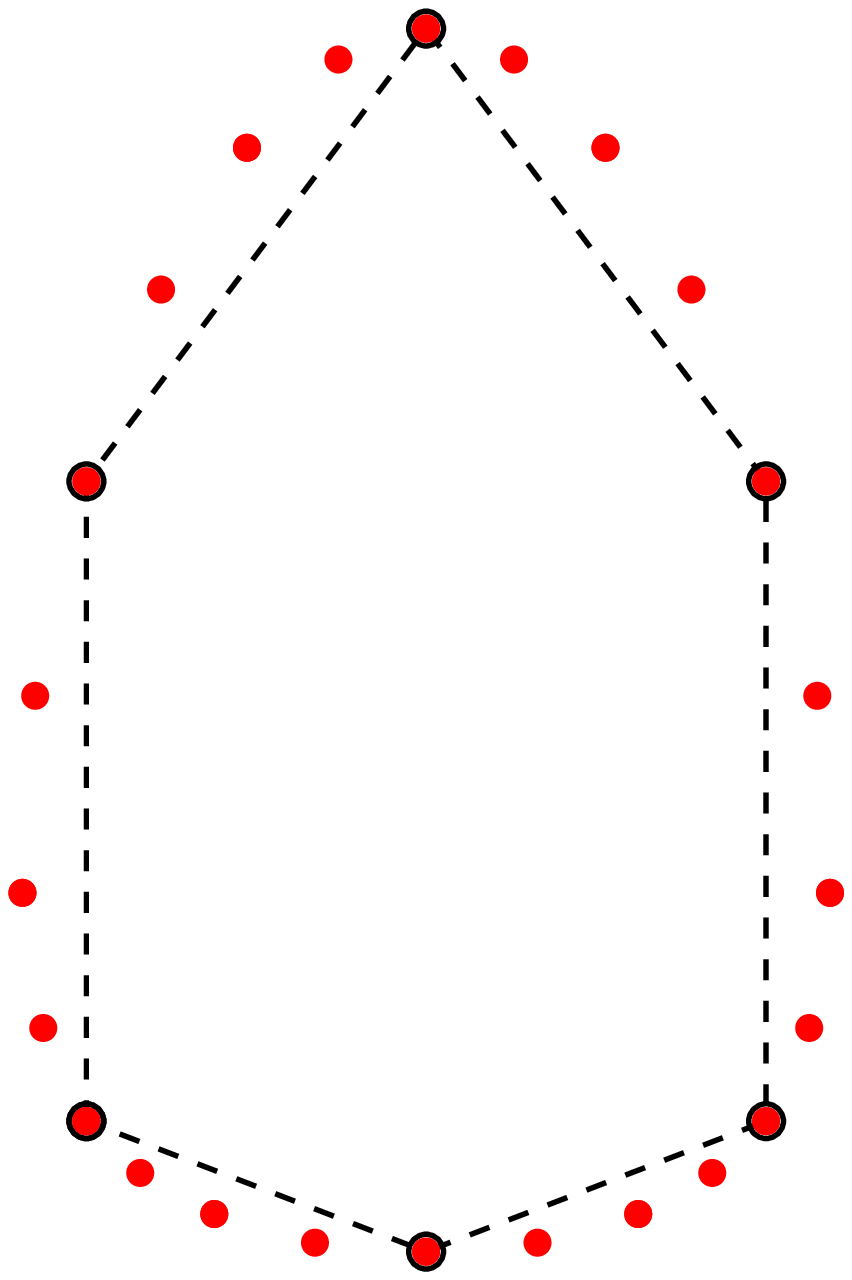}}\hspace{0.5cm}
\resizebox{2.1cm}{!}{\includegraphics{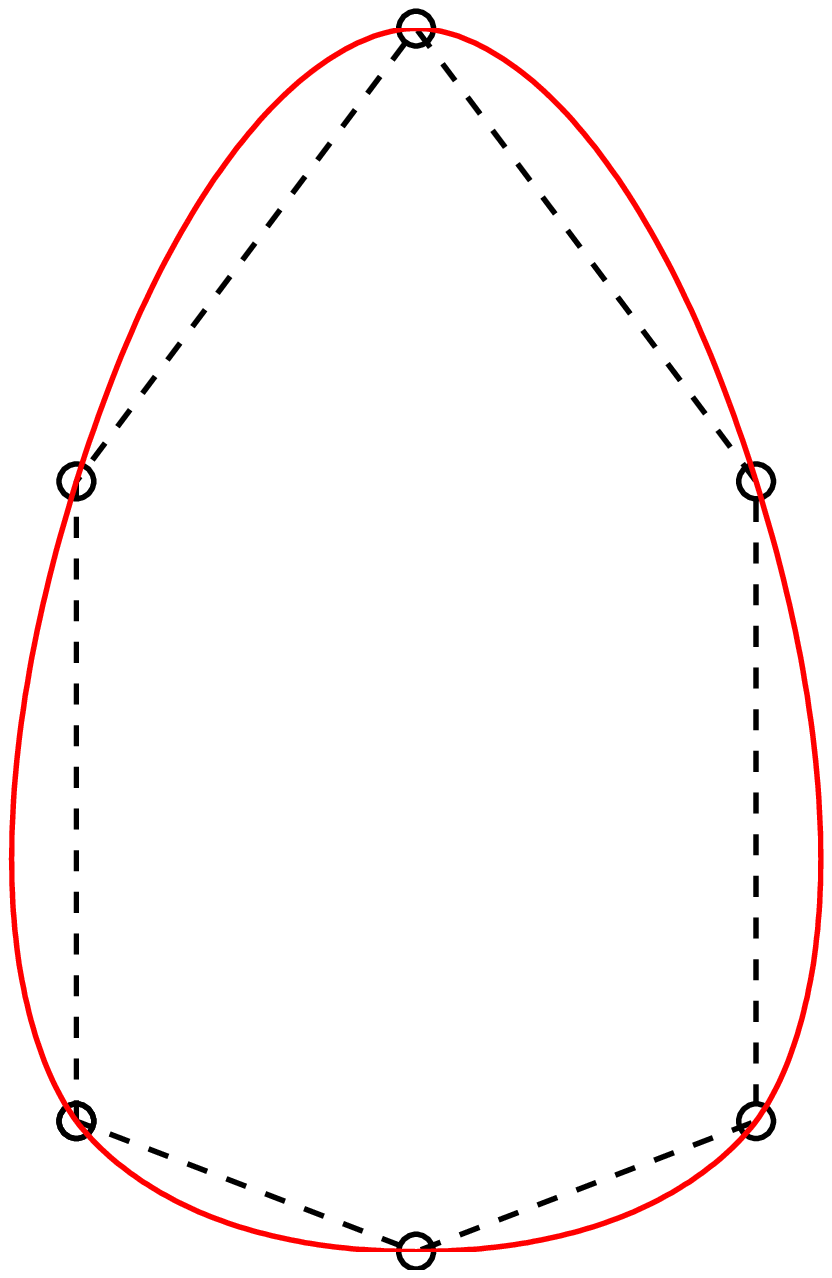}}\vspace{0.3cm}\\
\resizebox{1.8cm}{!}{\includegraphics{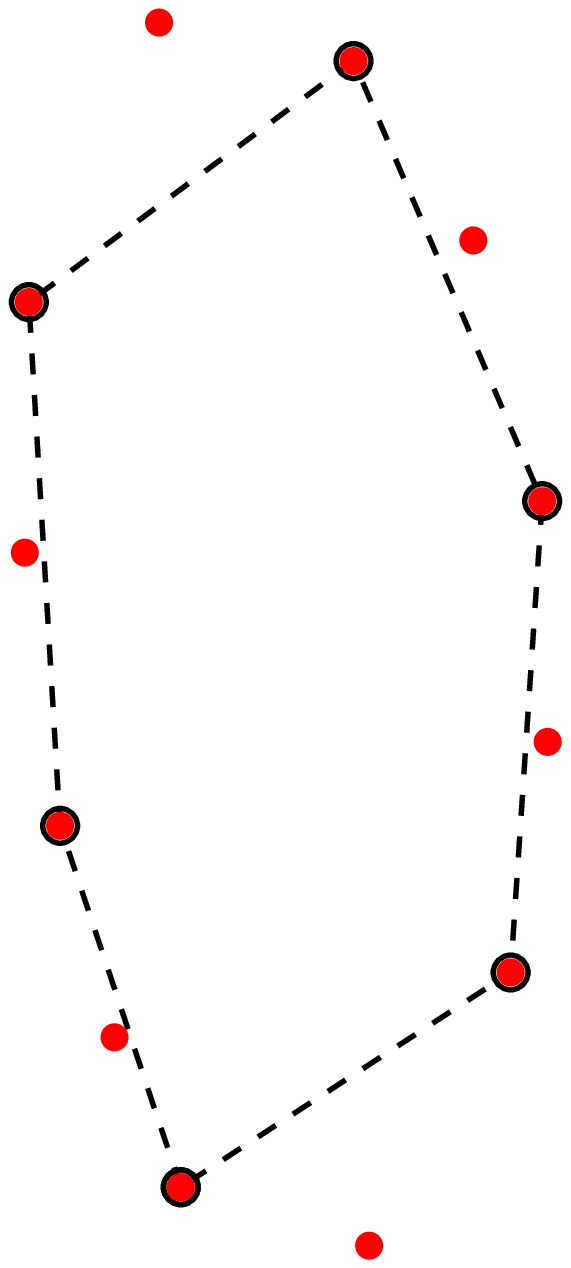}}\hspace{0.8cm}
\resizebox{1.8cm}{!}{\includegraphics{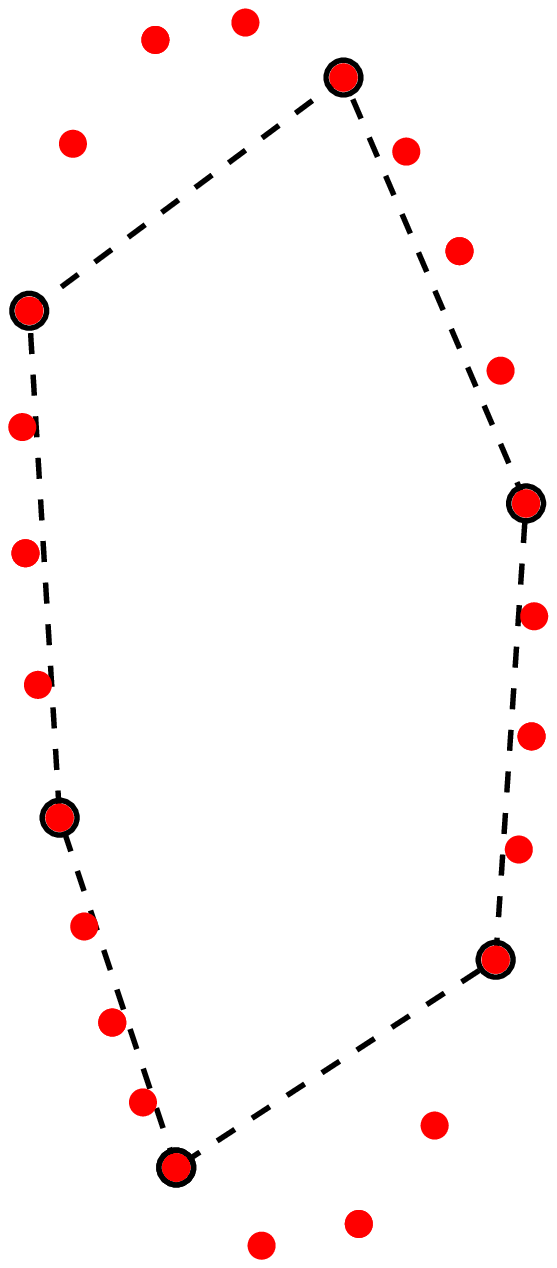}}\hspace{0.8cm}
\resizebox{1.6cm}{!}{\includegraphics{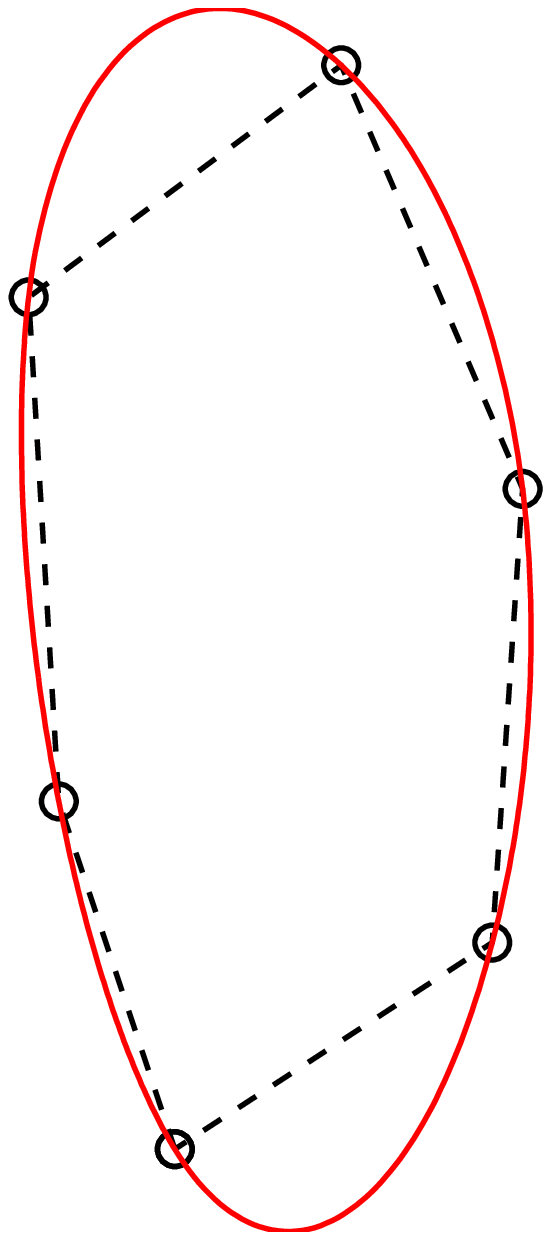}}
\caption{Application examples of the subdivision algorithm of
Subsection \ref{total_conv_algo} to nearly uniform closed polylines.
From left to right: points at 1st and 2nd level of refinement;
refined polyline after 6 steps of the algorithm.} \label{fig7}
\end{figure}

\begin{figure}[h!]
\centering
{\includegraphics[trim= 10mm 6mm 10mm 6mm, clip, width=3.2cm]{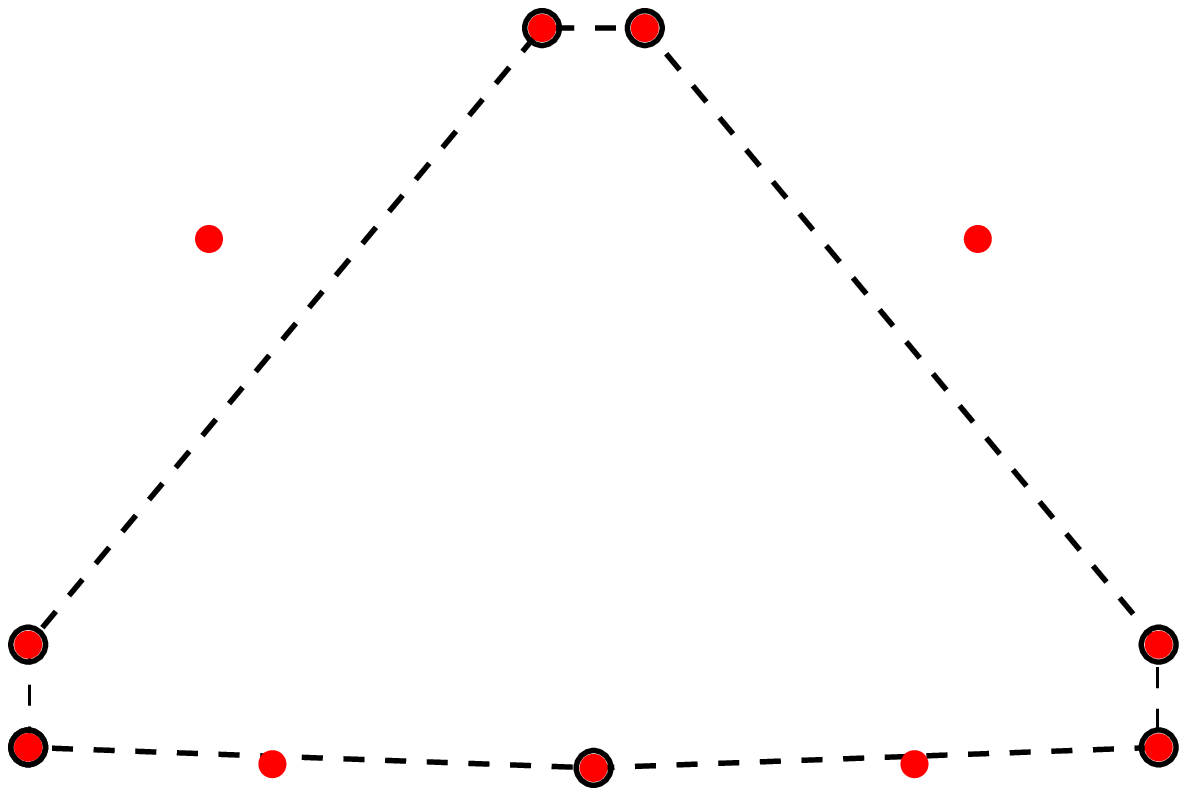}}\hspace{0.5cm}
{\includegraphics[trim= 10mm 6mm 10mm 6mm, clip, width=3.2cm]{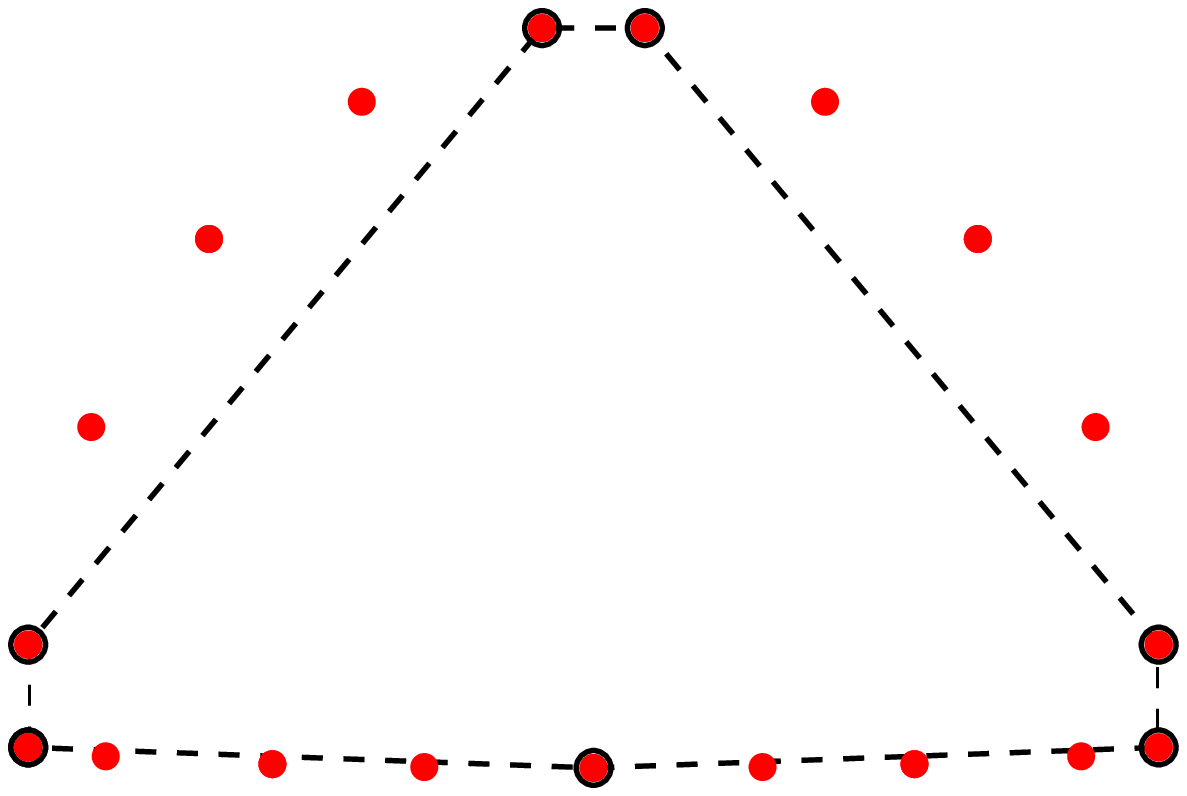}}\hspace{0.5cm}
{\includegraphics[trim= 10mm 6mm 10mm 6mm, clip, width=3.2cm]{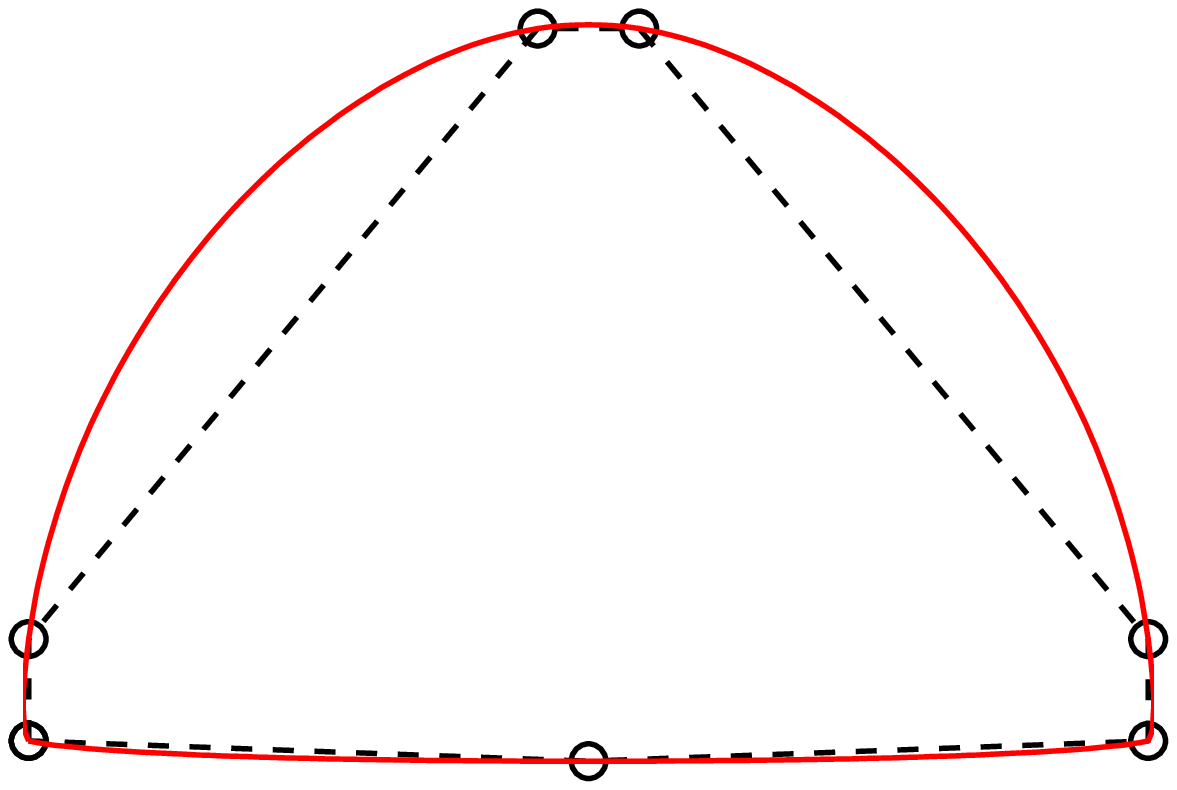}}\vspace{-0.1cm}\\
{\includegraphics[trim= 30mm 6mm 30mm 6mm, clip, width=2.1cm]{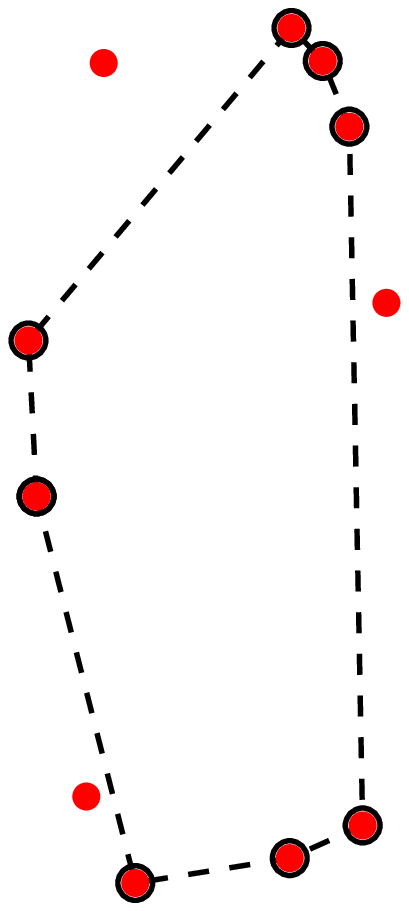}}\hspace{0.8cm}
{\includegraphics[trim= 30mm 6mm 30mm 6mm, clip, width=2.1cm]{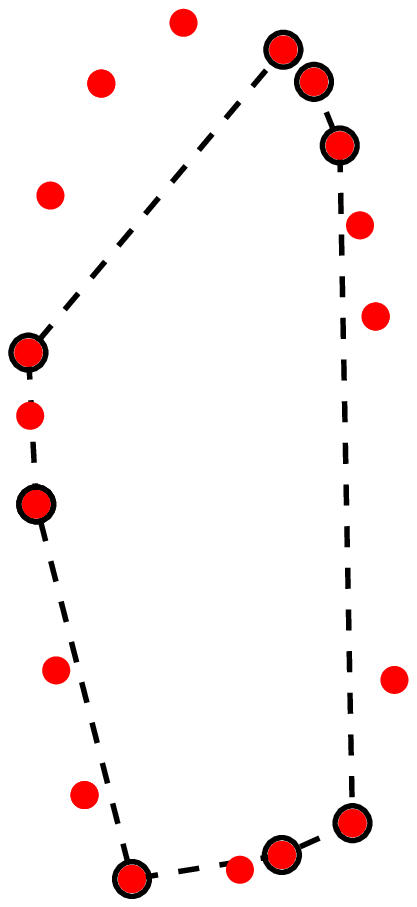}}\hspace{0.8cm}
{\includegraphics[trim= 30mm 6mm 30mm 6mm, clip, width=2.1cm]{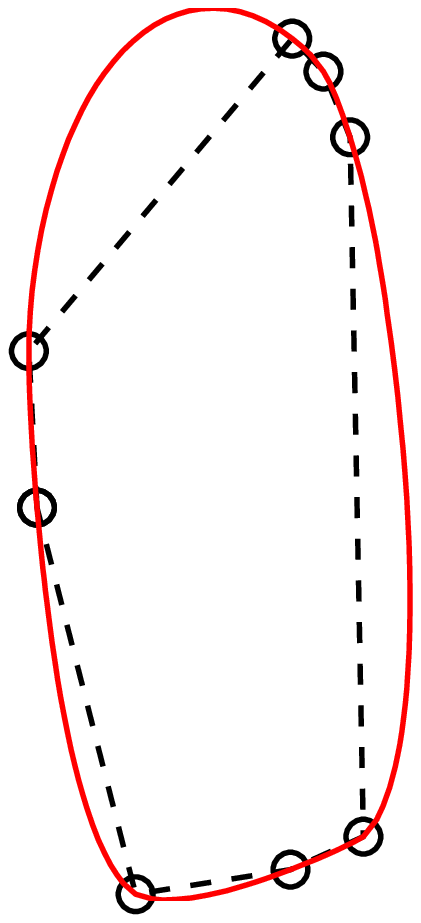}}
\caption{Application examples of the adaptive subdivision algorithm of
Section \ref{secadapt} to highly non-uniform closed
polylines. From left to right: points at 1st and 2nd level of
refinement; refined polyline after 6 steps of the algorithm.}
\label{fig6}
\end{figure}

\begin{figure}[h!]
\centering
\resizebox{4.1cm}{!}{\includegraphics{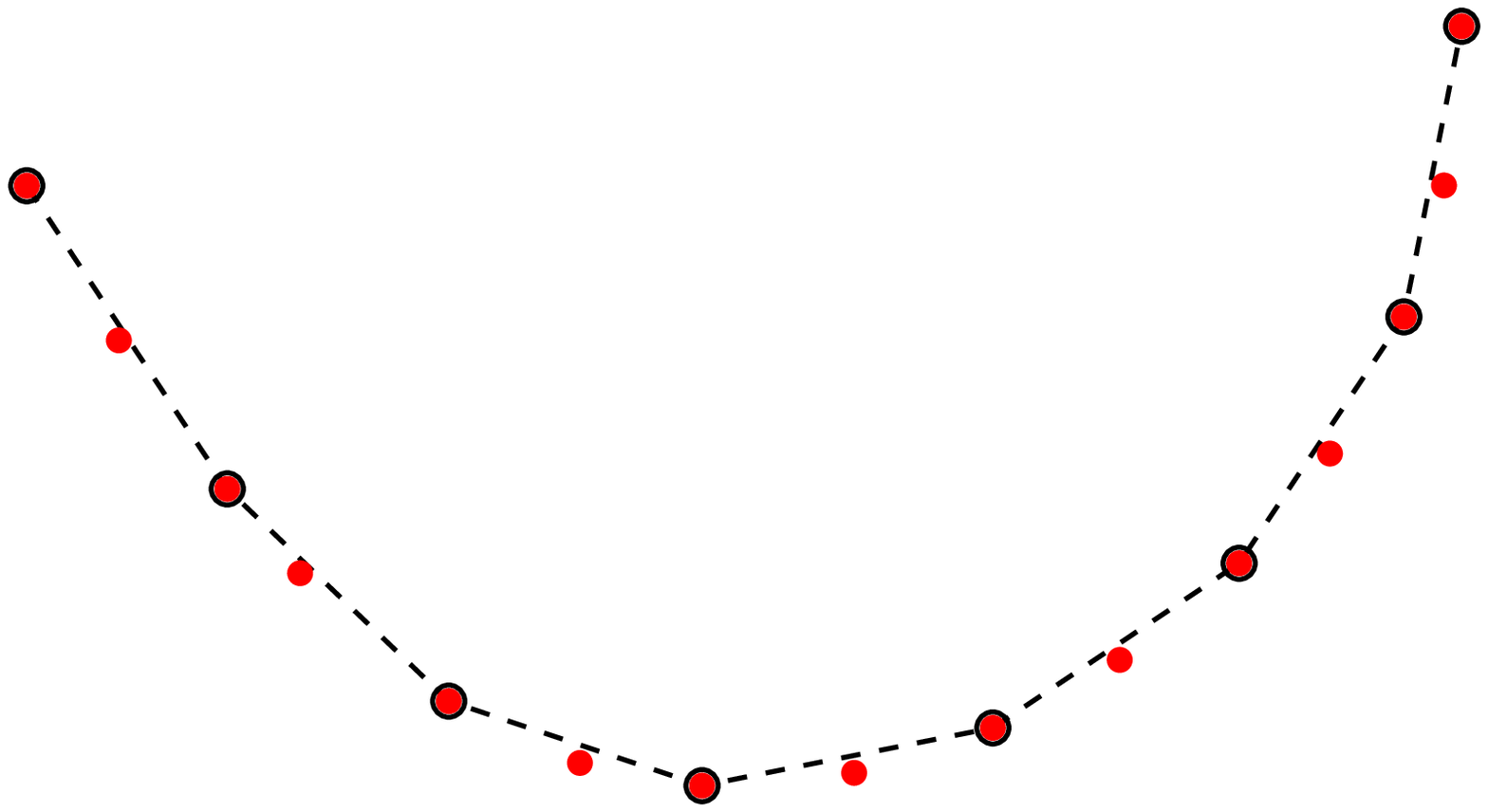}}\hspace{0.2cm}
\resizebox{4.1cm}{!}{\includegraphics{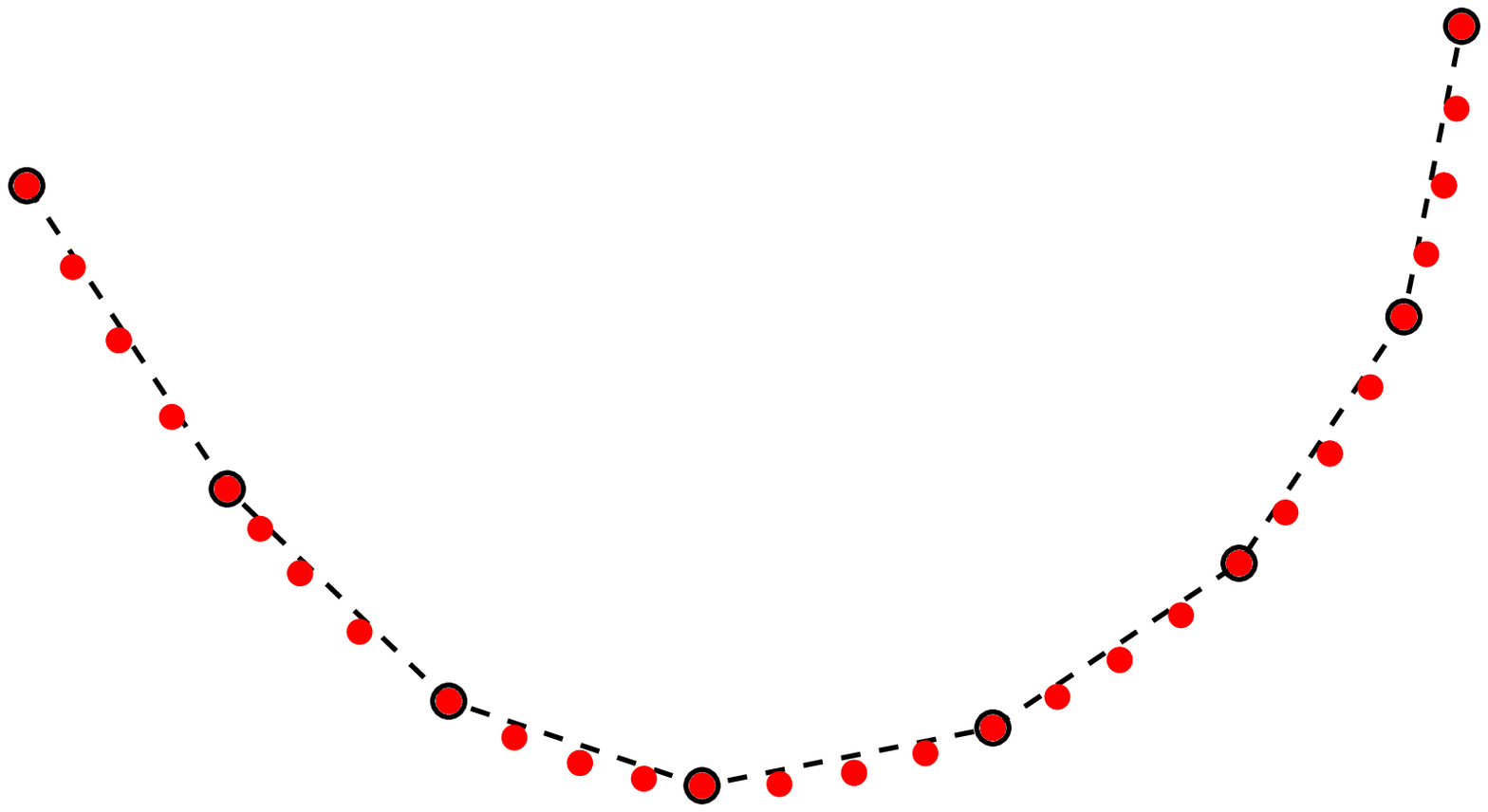}}\hspace{0.2cm}
\resizebox{4.2cm}{!}{\includegraphics{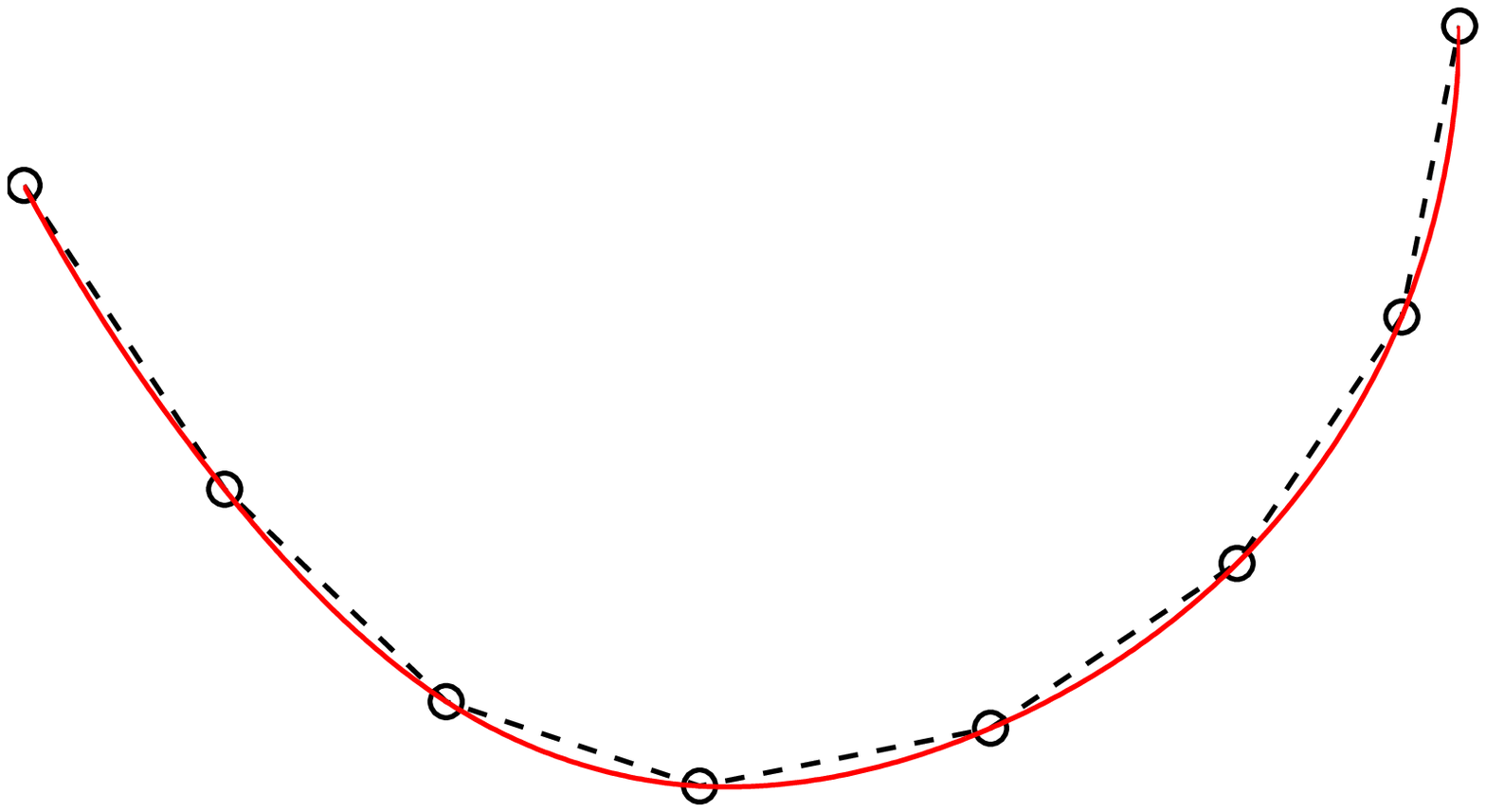}}\vspace{0.1cm}\\
\resizebox{1.5cm}{!}{\includegraphics{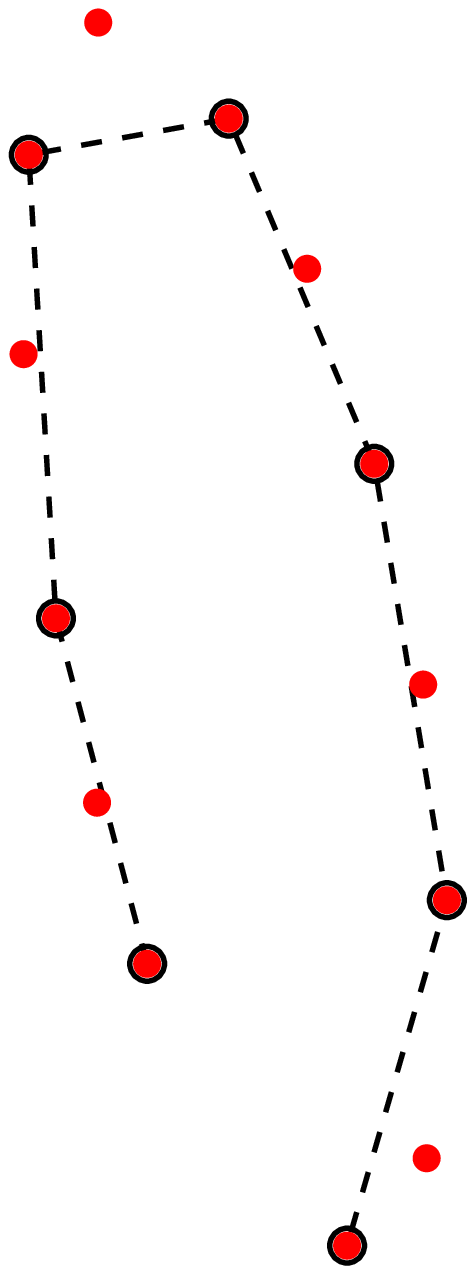}}\hspace{0.6cm}
\resizebox{1.5cm}{!}{\includegraphics{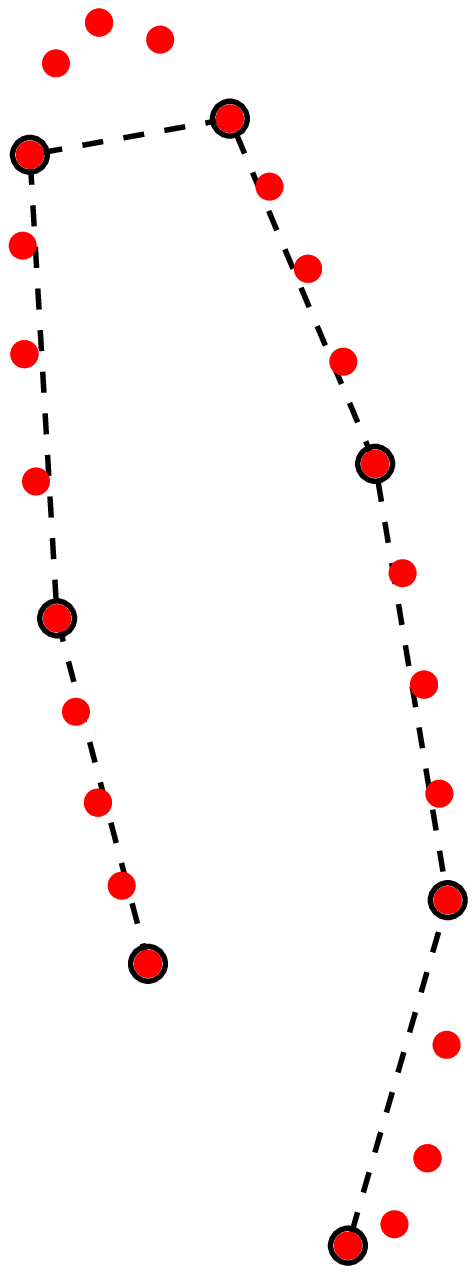}}\hspace{0.6cm}
\resizebox{1.4cm}{!}{\includegraphics{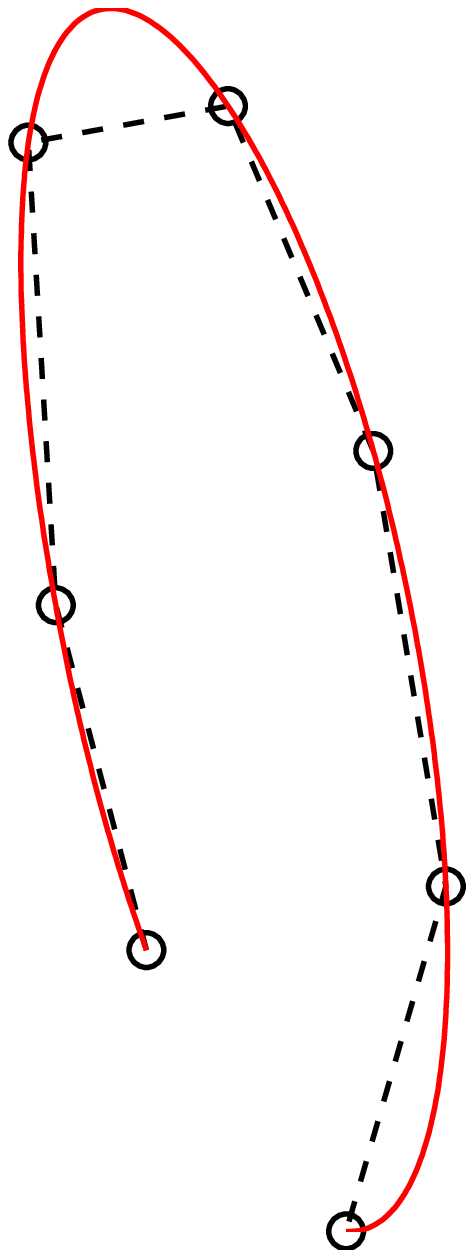}}
\caption{Application examples of the subdivision algorithm of
Subsection \ref{total_conv_algo} to nearly uniform open polylines.
From left to right: points at 1st and 2nd level of refinement;
refined polyline after 6 steps of the algorithm.} \label{fig9}
\end{figure}

\begin{figure}[h!]
\centering
{\includegraphics[trim= 10mm 8mm 10mm 8mm, clip, width=3.2cm]{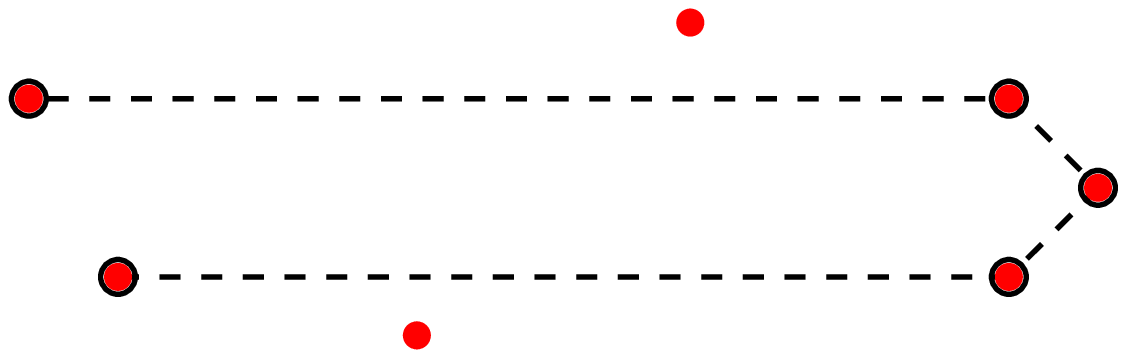}}\hspace{0.3cm}
{\includegraphics[trim= 10mm 8mm 10mm 8mm, clip, width=3.2cm]{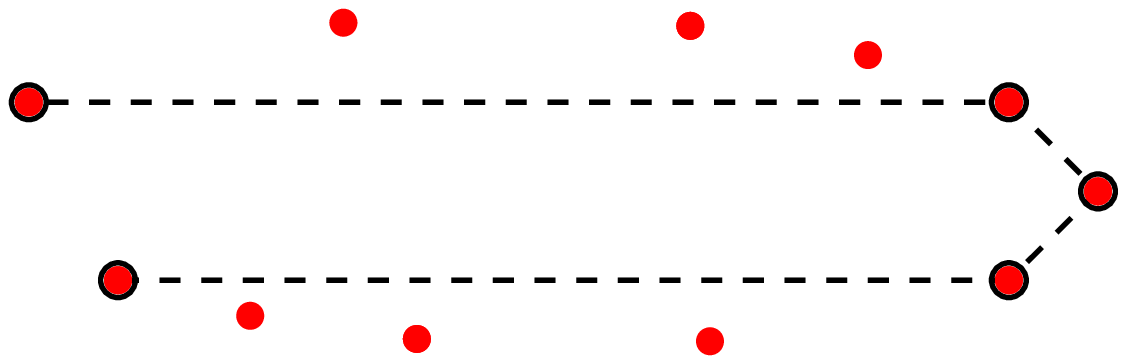}}\hspace{0.3cm}
{\includegraphics[trim= 10mm 8mm 10mm 8mm, clip, width=3.2cm]{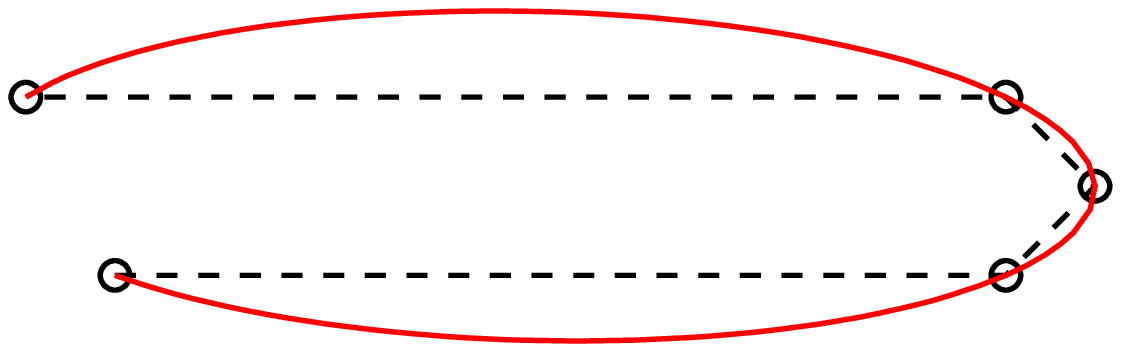}}\vspace{-0.2cm}\\
{\includegraphics[trim= 30mm 6mm 30mm 6mm, clip, width=2.25cm]{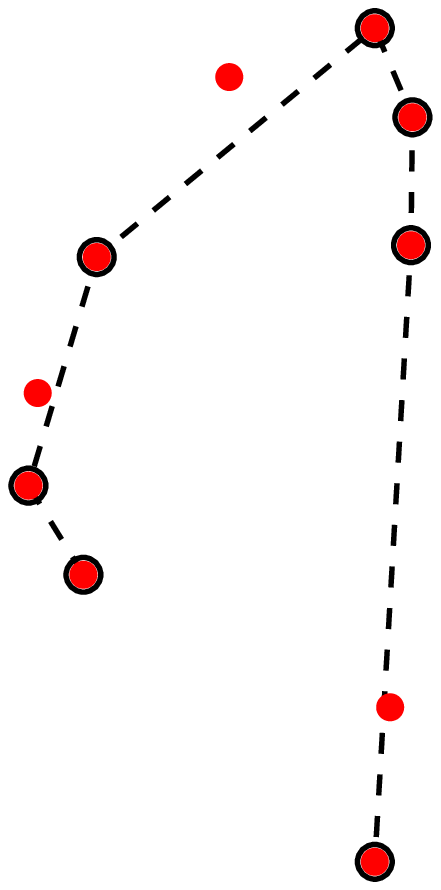}}\hspace{0.6cm}
{\includegraphics[trim= 30mm 6mm 30mm 6mm, clip, width=2.25cm]{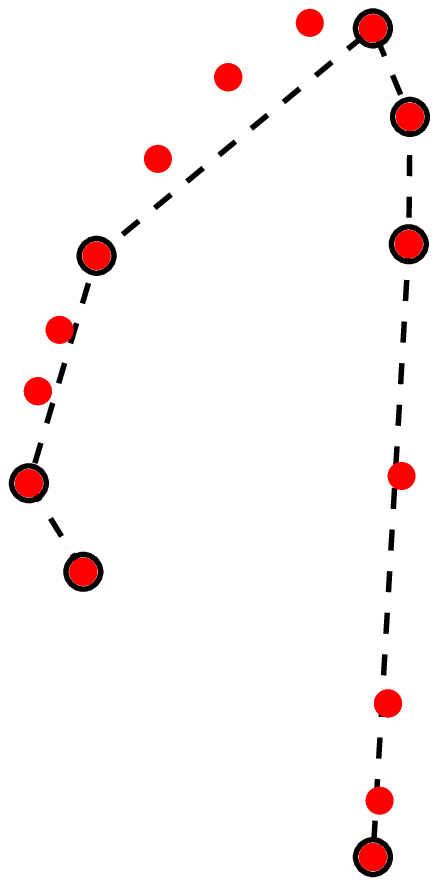}}\hspace{0.6cm}
{\includegraphics[trim= 30mm 6mm 30mm 6mm, clip, width=2.25cm]{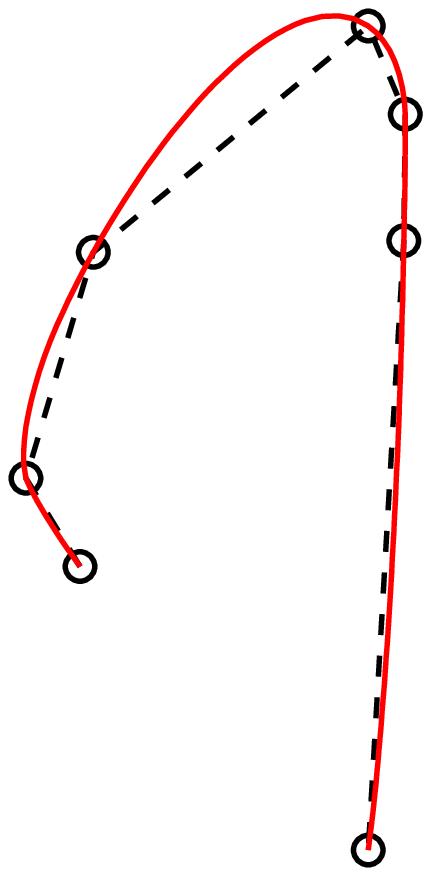}}
\caption{Application examples of the adaptive subdivision algorithm of
Section \ref{secadapt} to highly non-uniform open
polylines. From left to right: points at 1st and 2nd level of
refinement; refined polyline after 6 steps of the algorithm.}
\label{fig8}
\end{figure}

\begin{figure}[h!]
\centering
\hspace{-0.2cm}
{\includegraphics[trim= 1mm 4mm 1mm 4mm, clip, width=3.0cm]{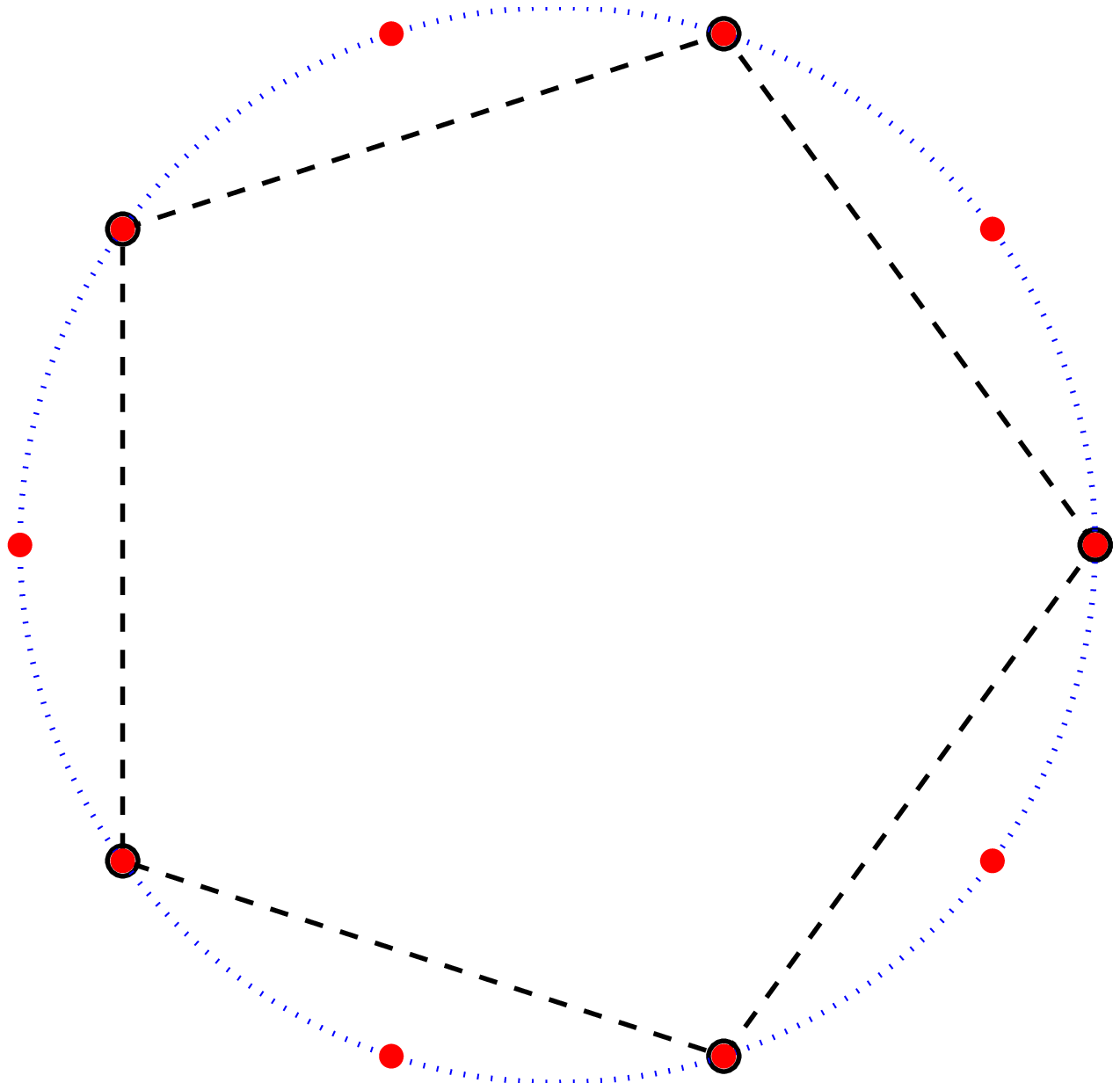}}\hspace{0.8cm}
{\includegraphics[trim= 1mm 4mm 1mm 4mm, clip, width=3.0cm]{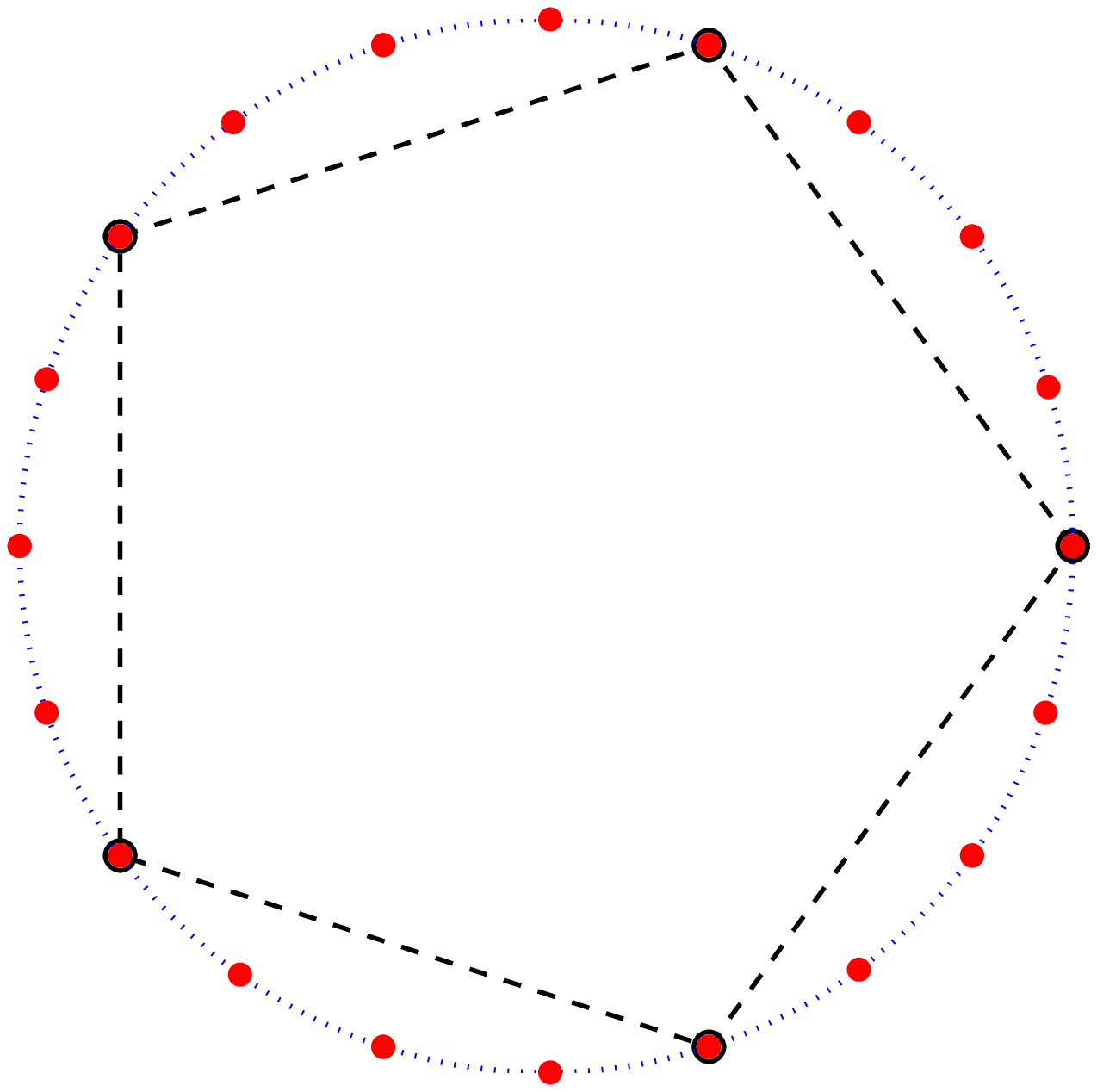}}\hspace{0.8cm}
{\includegraphics[trim= -1mm 1mm -1mm 1mm, clip, width=2.9cm]{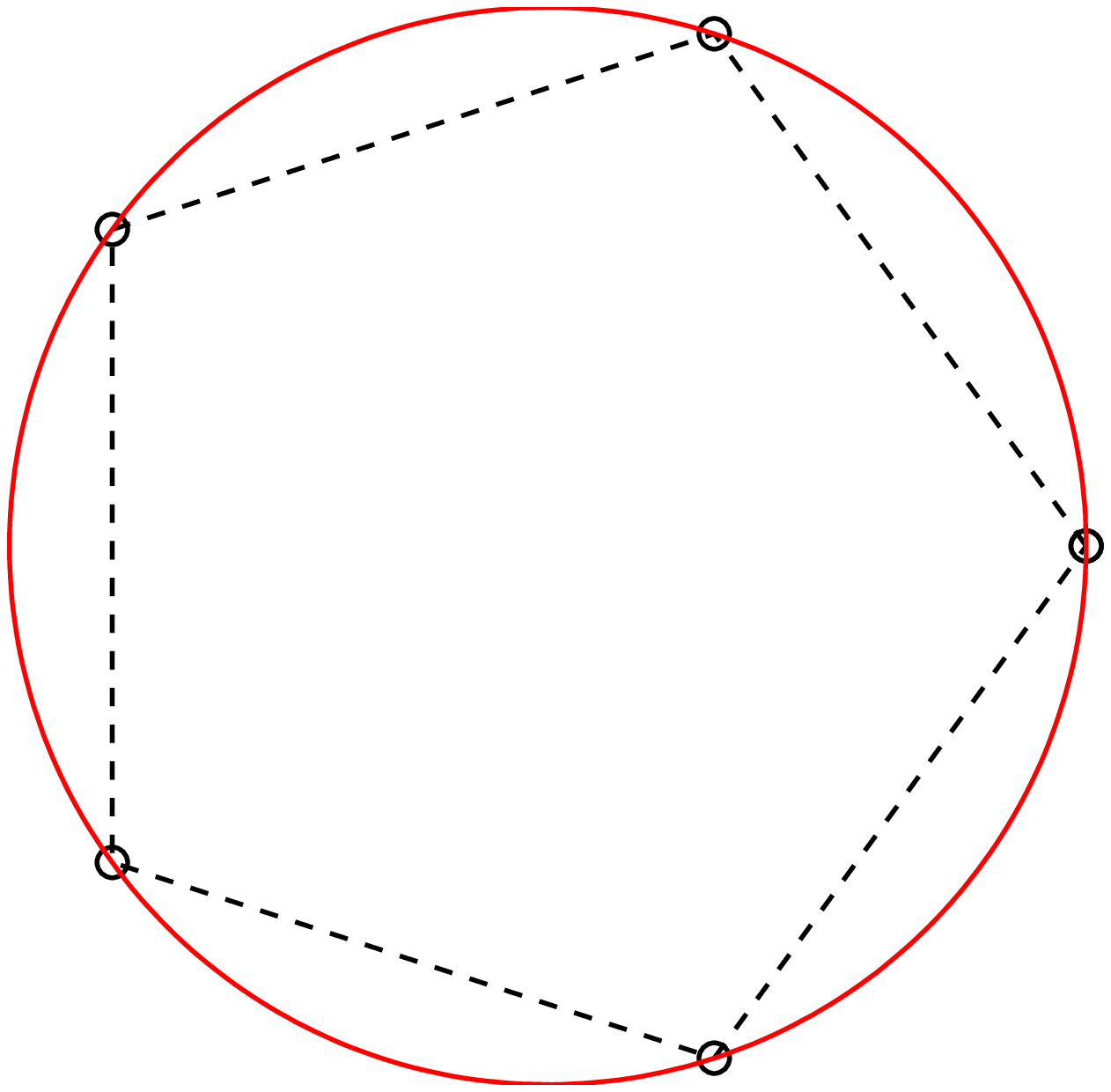}}\vspace{0.3cm}\\
{\includegraphics[trim= 1mm 6mm 1mm 6mm, clip, width=3.6cm]{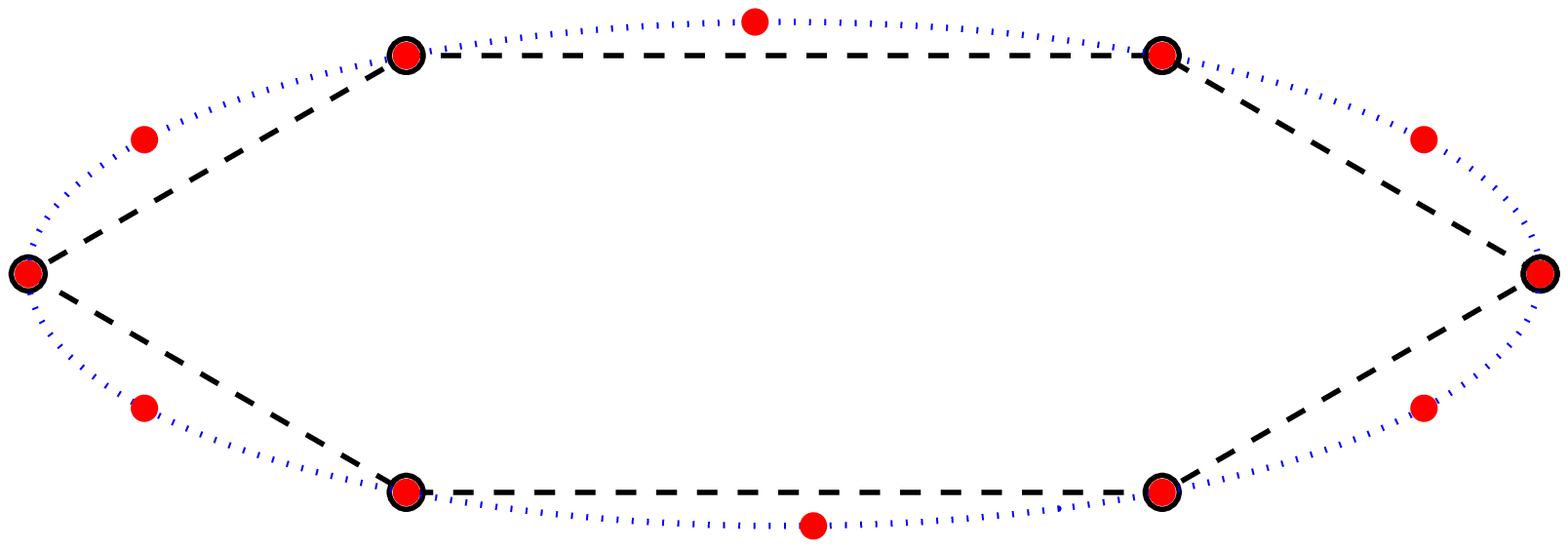}}\hspace{0.3cm}
{\includegraphics[trim= 1mm 6mm 1mm 6mm, clip, width=3.6cm]{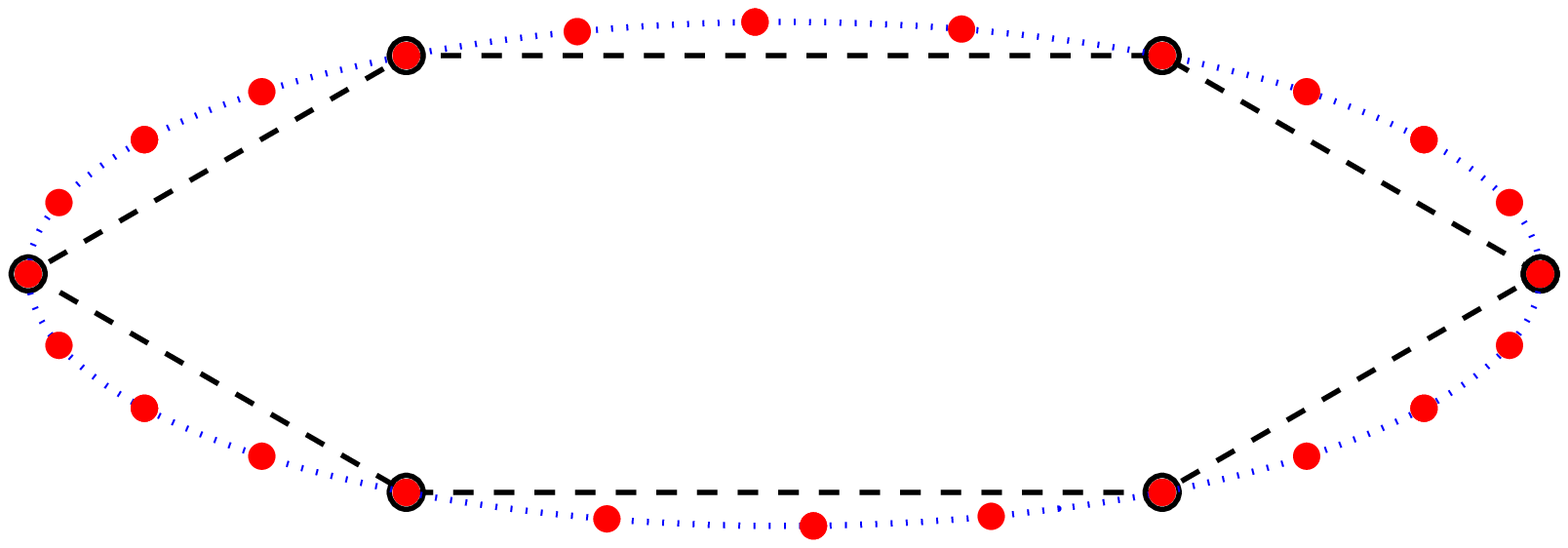}}\hspace{0.3cm}
{\includegraphics[trim= 0mm -3mm 1mm -3mm, clip, width=3.6cm]{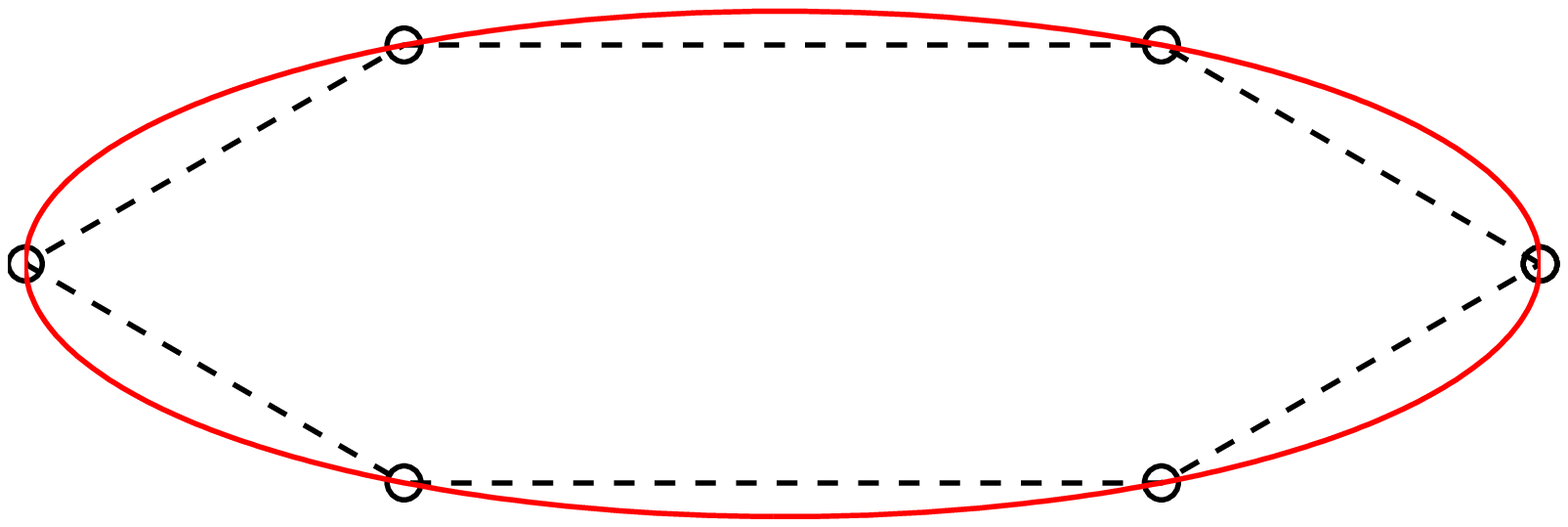}}\vspace{0.3cm}\\
{\includegraphics[trim= 6mm 1mm 6mm 1mm, clip, width=3.3cm]{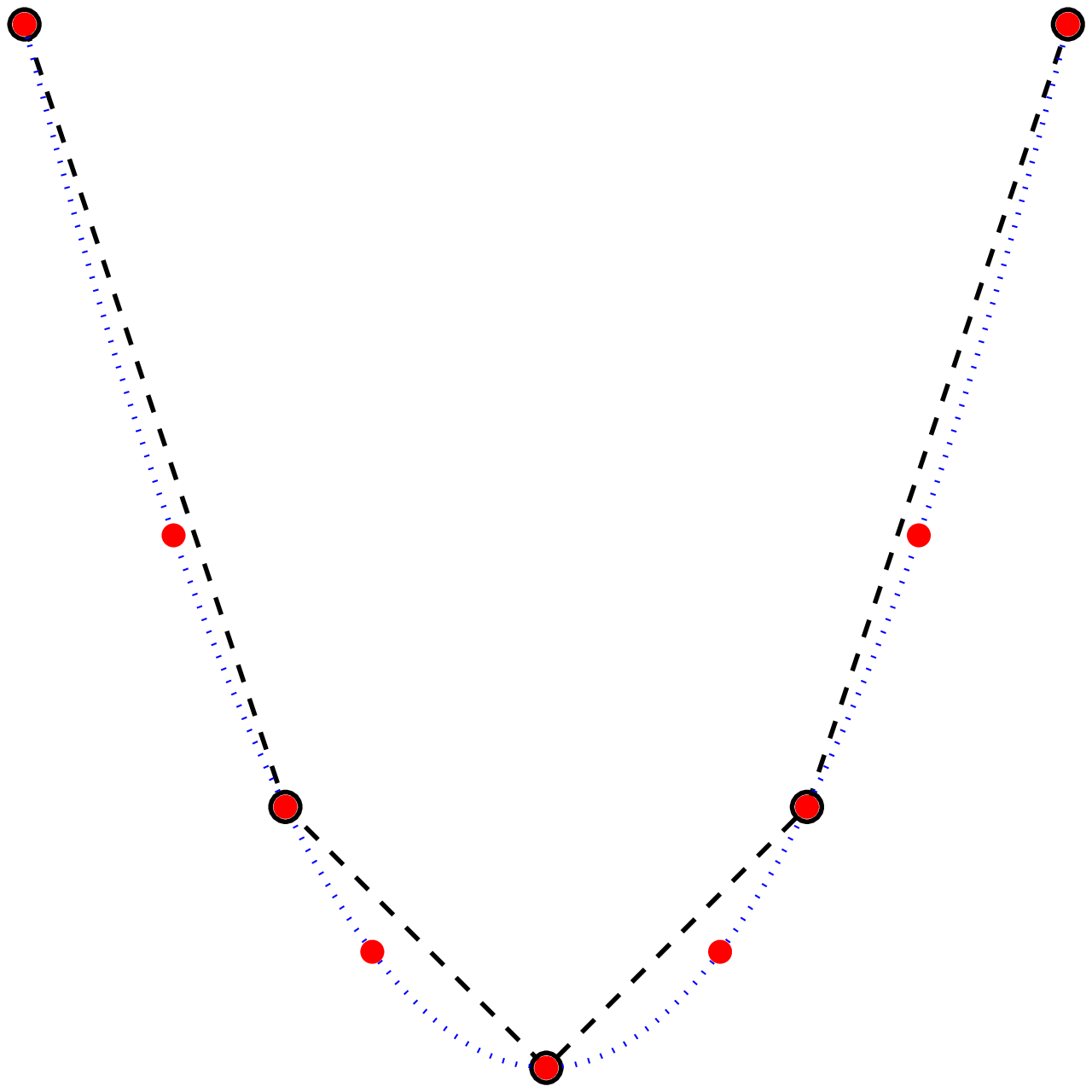}}\hspace{0.4cm}
{\includegraphics[trim= 6mm 1mm 6mm 1mm, clip, width=3.3cm]{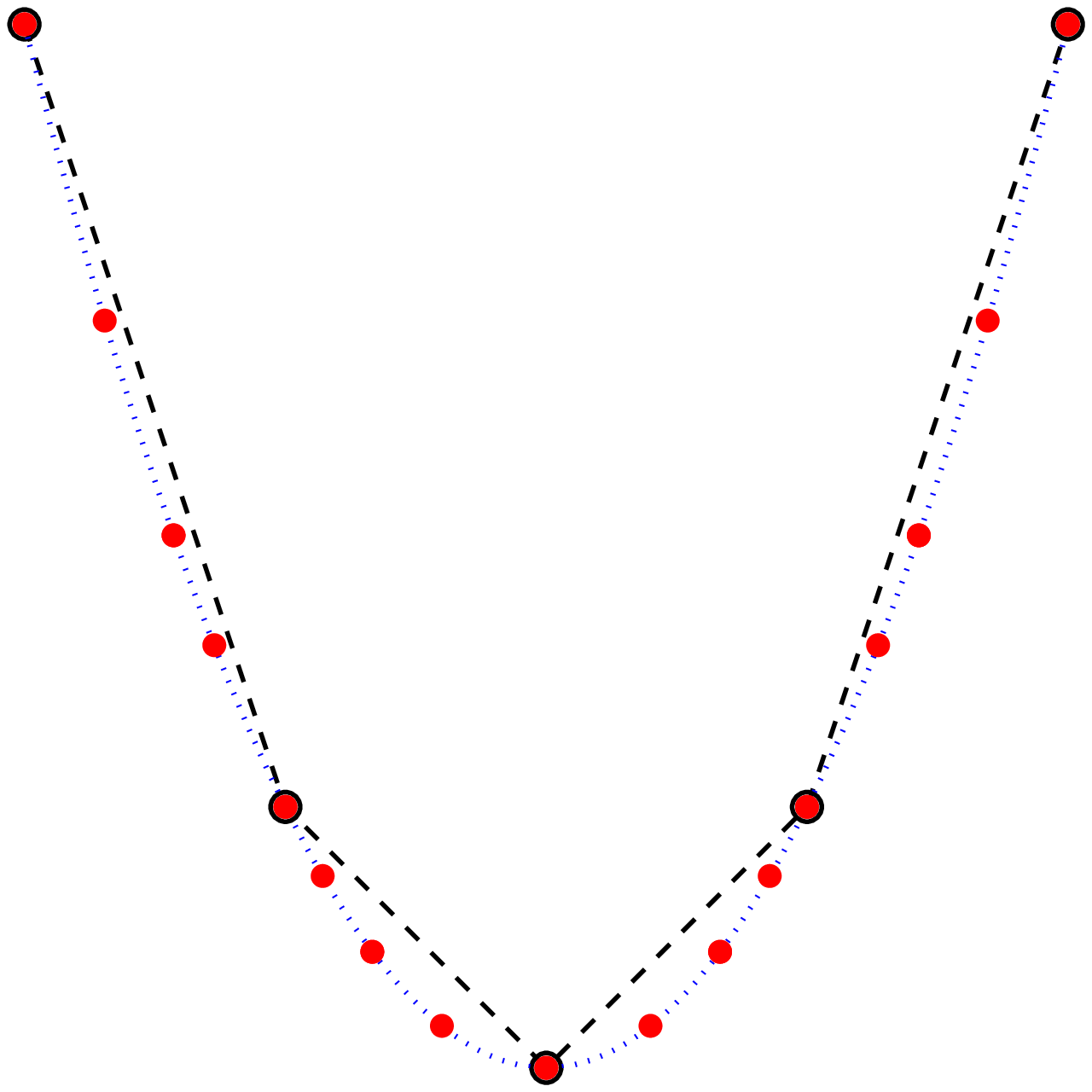}}\hspace{0.4cm}
{\includegraphics[trim= 1mm -2mm 1mm 1mm, clip, width=3.3cm]{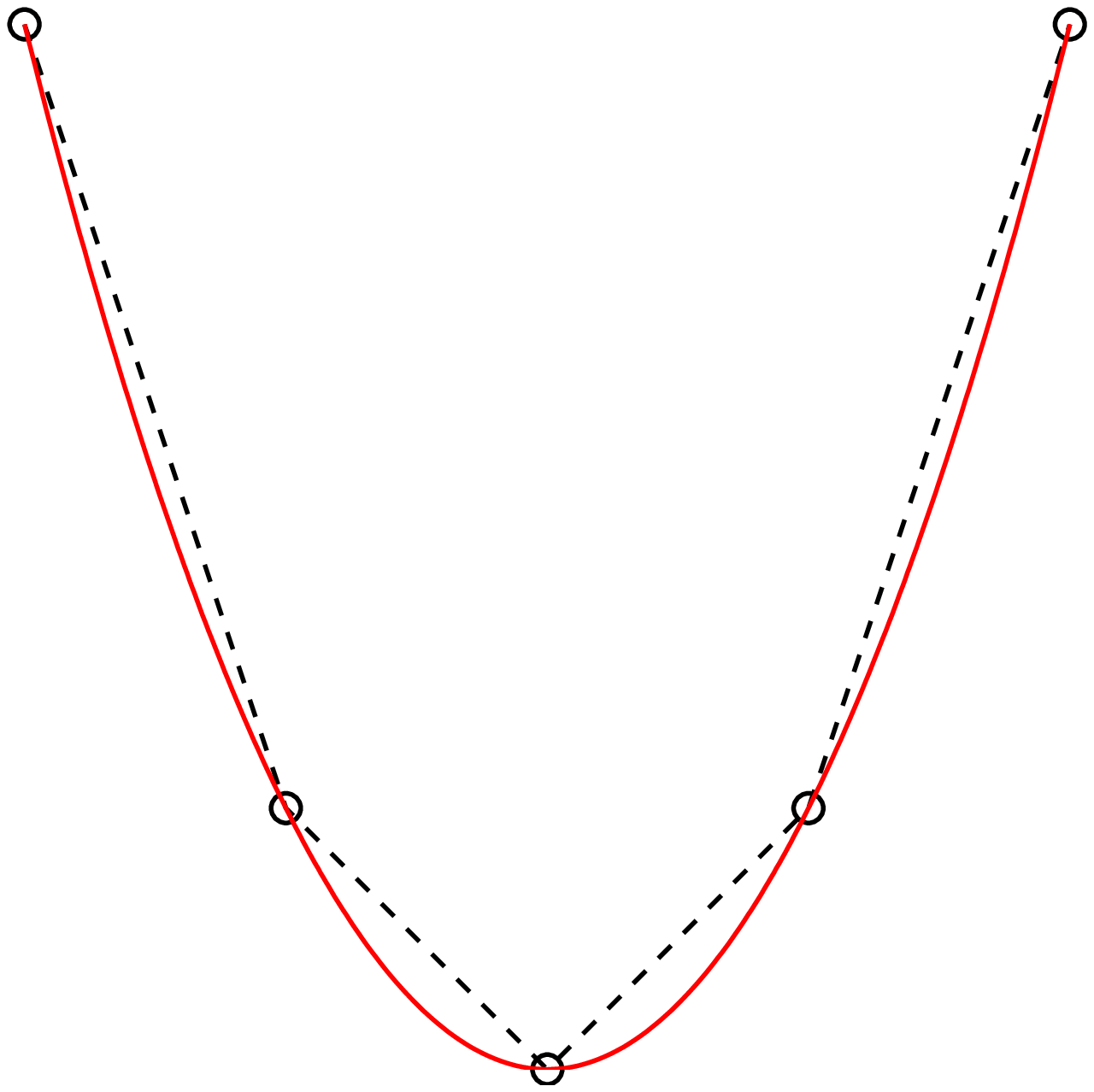}}\vspace{0.3cm}\\
\hspace{-0.5cm}
{\includegraphics[trim= 6mm 2mm 6mm 2mm, clip, width=1.5cm]{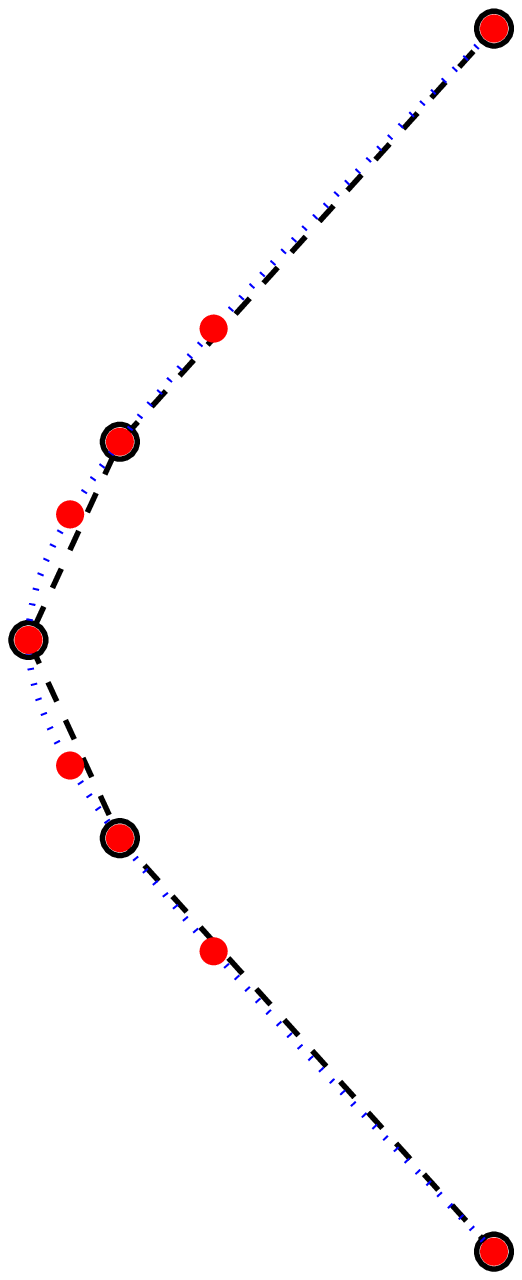}}\hspace{1.8cm}
{\includegraphics[trim= 6mm 2mm 6mm 2mm, clip, width=1.5cm]{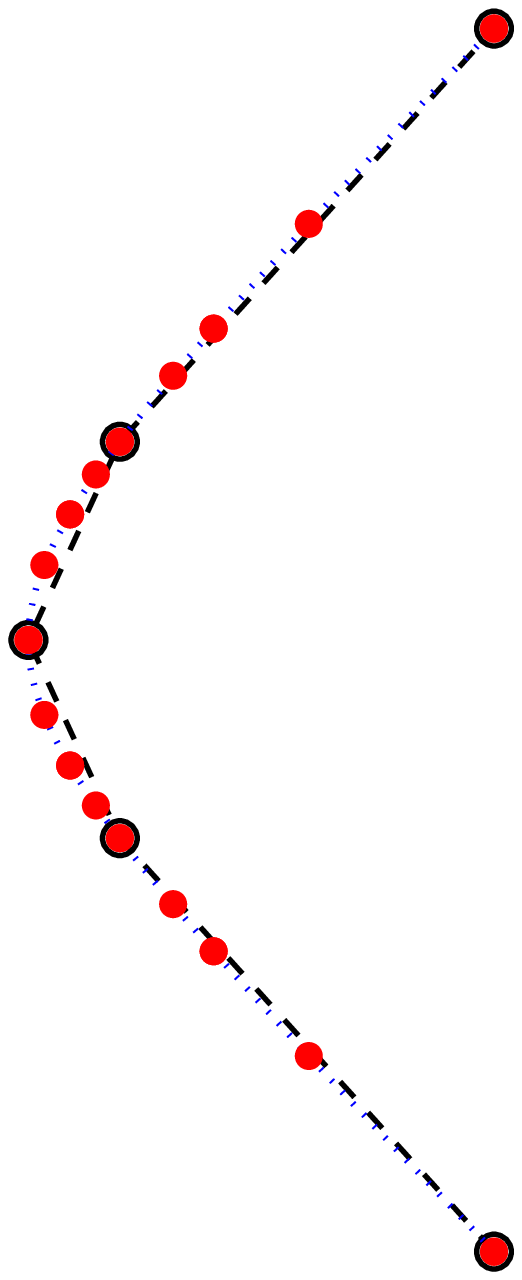}}\hspace{1.8cm}
{\includegraphics[trim= 2mm -1mm 2mm 1mm, clip, width=1.5cm]{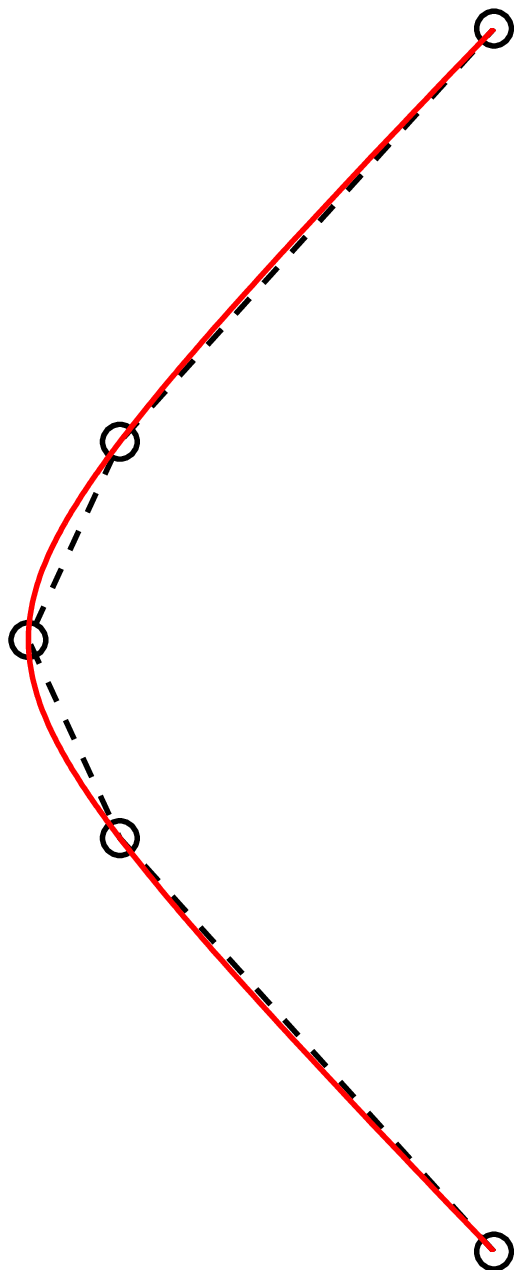}}
\caption{Uniform case: reproduction of conic sections from uniform
samples by applying the subdivision algorithm of Subsection
\ref{total_conv_algo}. The dotted blue line is the conic section to
be reconstructed. From left to right: points at 1st and 2nd level of
refinement; refined polyline after 6 steps of the algorithm.}
\label{fig5a}
\end{figure}

\begin{figure}[h!]
\centering
{\includegraphics[trim= 20mm 4mm 20mm 4mm, clip, width=2.6cm]{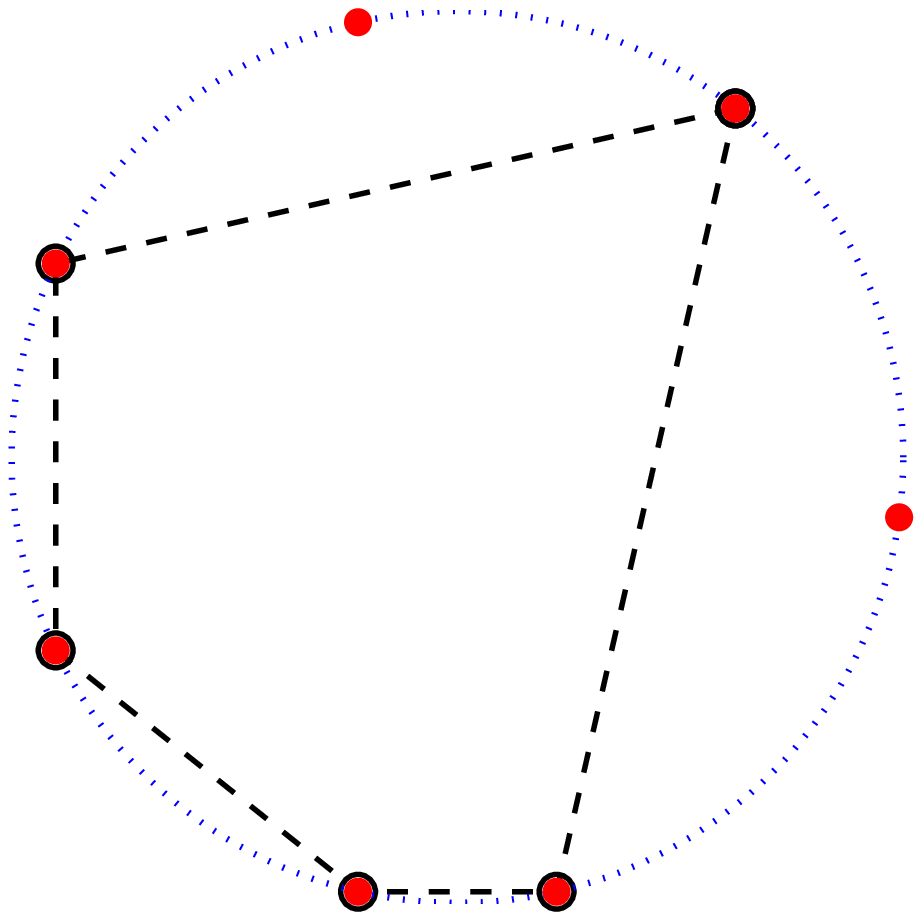}}\hspace{0.4cm}
{\includegraphics[trim= 20mm 4mm 20mm 4mm, clip, width=2.6cm]{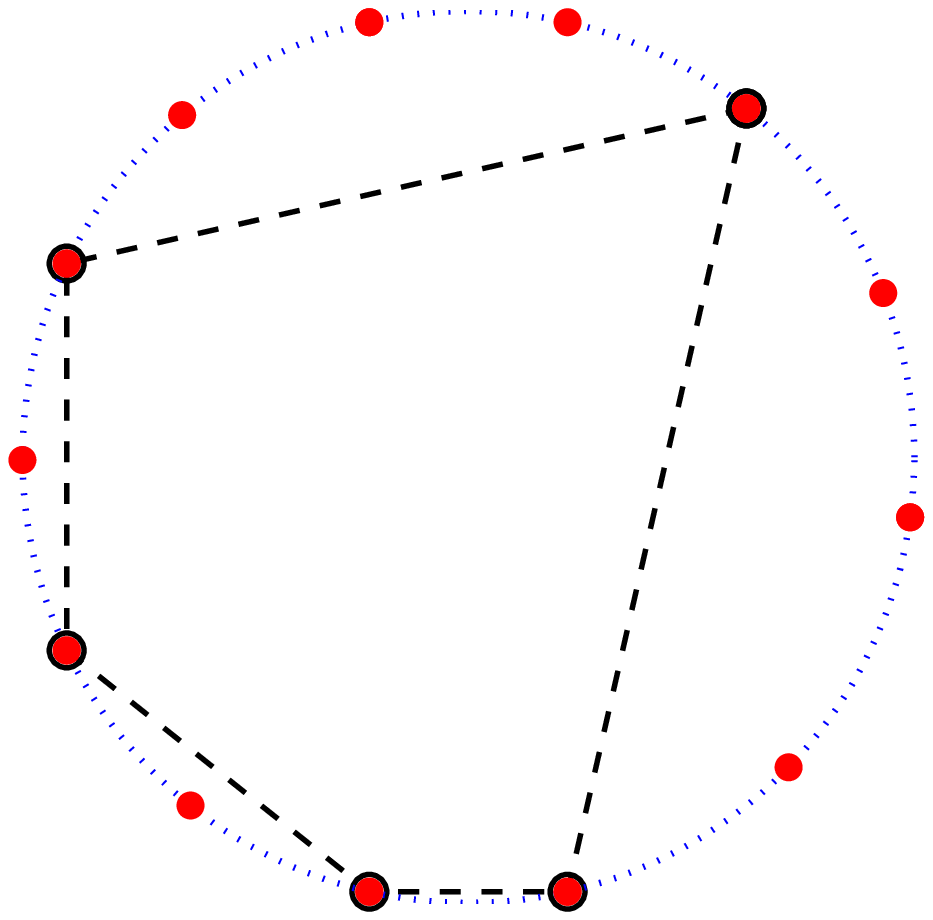}}\hspace{0.4cm}
{\includegraphics[trim= 20mm 4mm 20mm 4mm, clip, width=2.6cm]{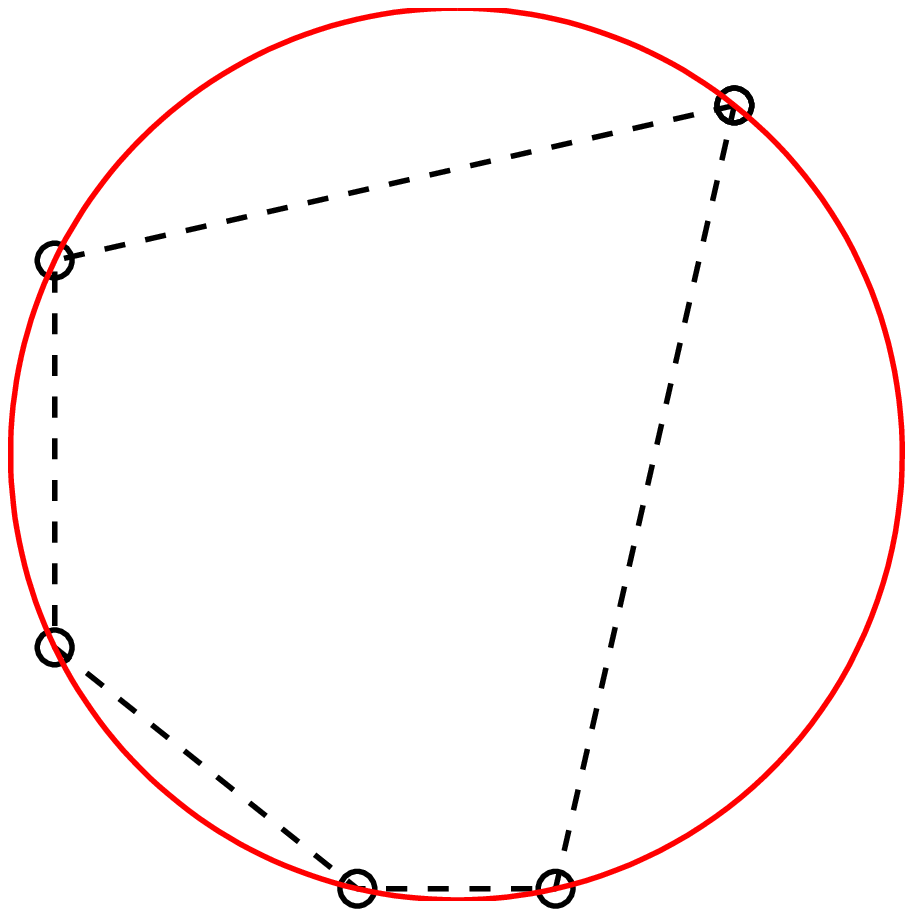}}\vspace{-0.2cm}\\
{\includegraphics[trim= 10mm 6mm 10mm 6mm, clip, width=2.9cm]{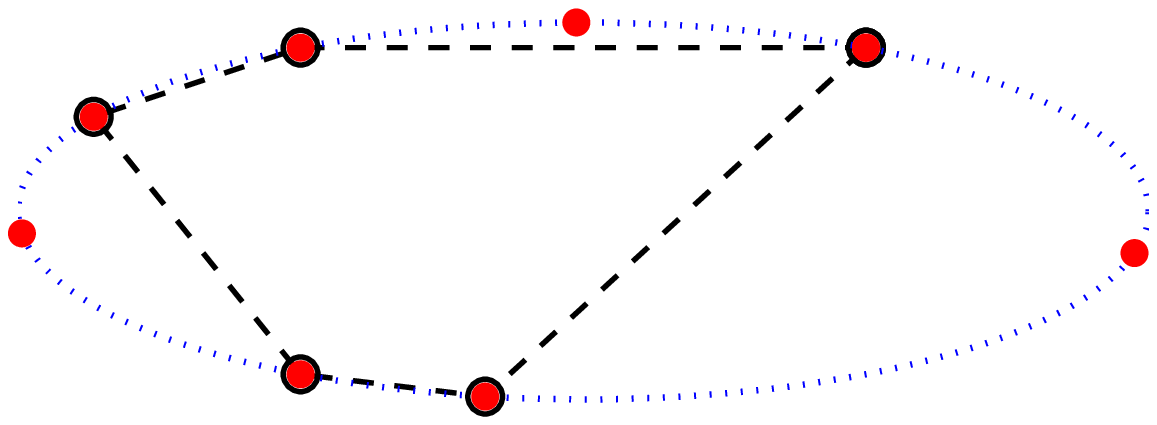}}\hspace{0.3cm}
{\includegraphics[trim= 10mm 6mm 10mm 6mm, clip, width=2.9cm]{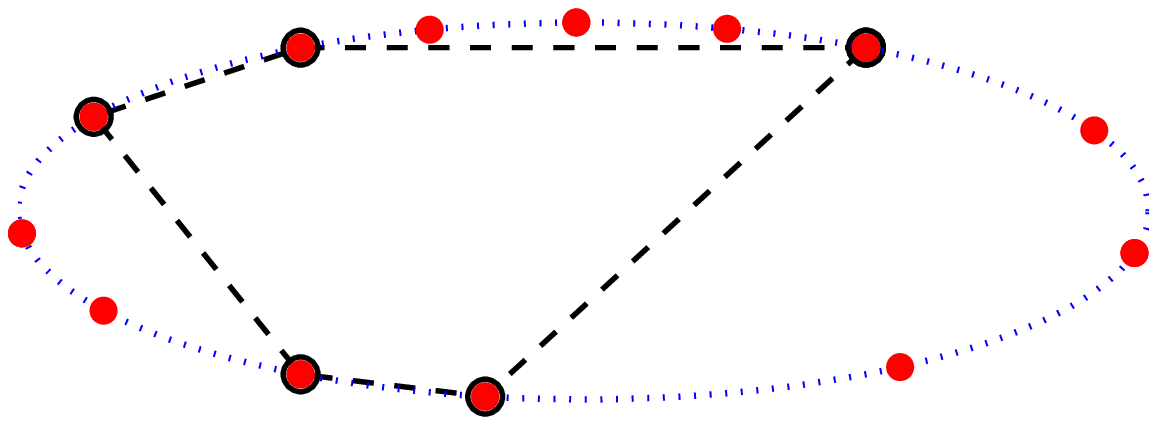}}\hspace{0.3cm}
{\includegraphics[trim= 10mm 6mm 10mm 6mm, clip, width=2.9cm]{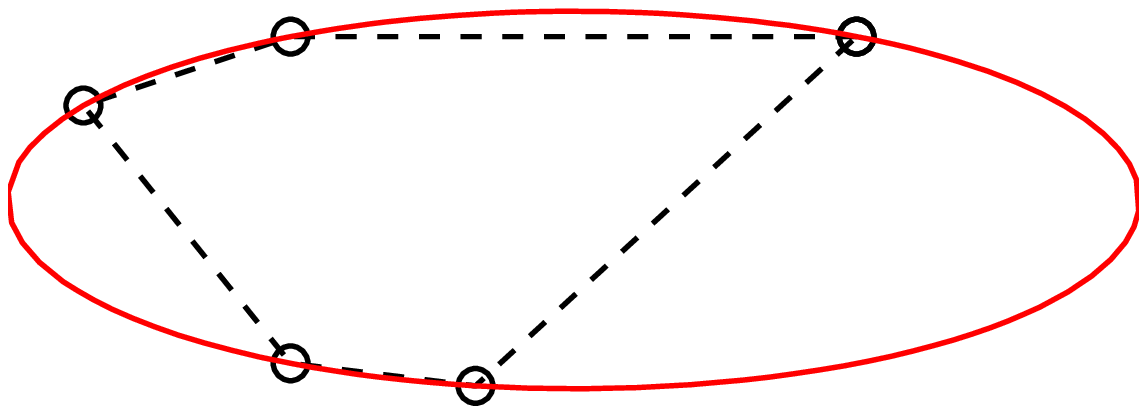}}\vspace{-0.4cm}\\
{\includegraphics[trim= 20mm 6mm 20mm 6mm, clip, width=2.8cm]{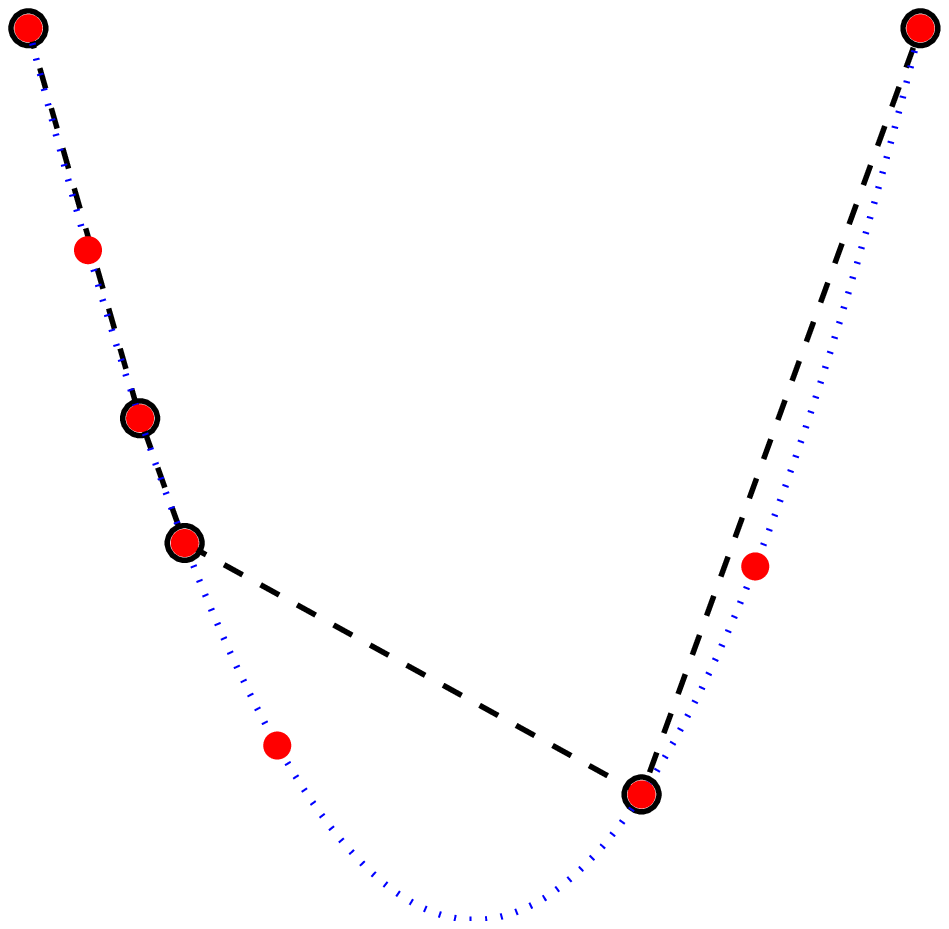}}\hspace{0.4cm}
{\includegraphics[trim= 20mm 6mm 20mm 6mm, clip, width=2.8cm]{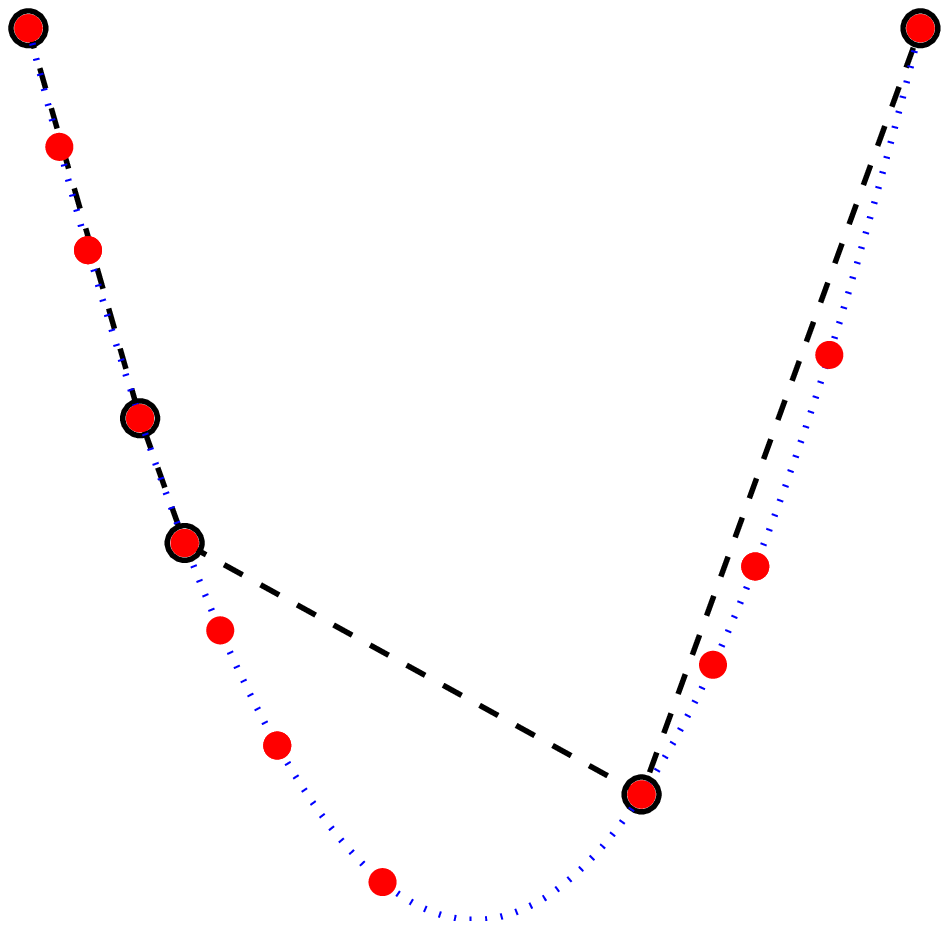}}\hspace{0.4cm}
{\includegraphics[trim= 20mm 8mm 20mm 6mm, clip, width=3.0cm]{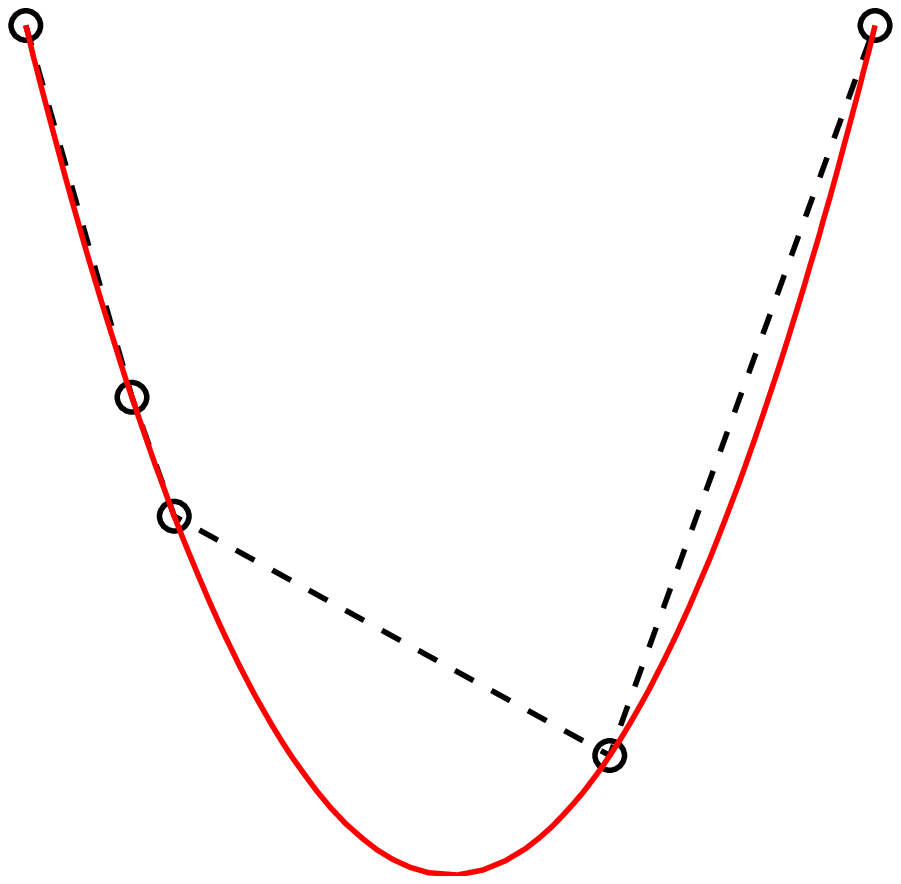}}\vspace{-0.0cm}\\
\hspace{-0.1cm}
{\includegraphics[trim= 25mm 6mm 20mm 6mm, clip, width=2.8cm]{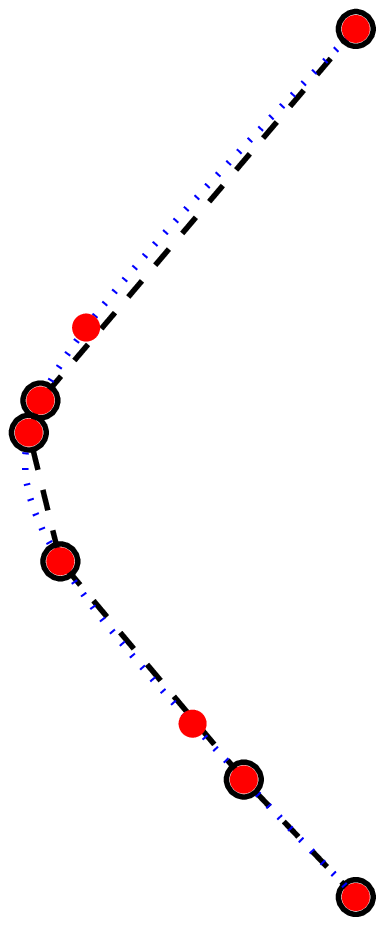}}\hspace{0.5cm}
{\includegraphics[trim= 25mm 6mm 20mm 6mm, clip, width=2.8cm]{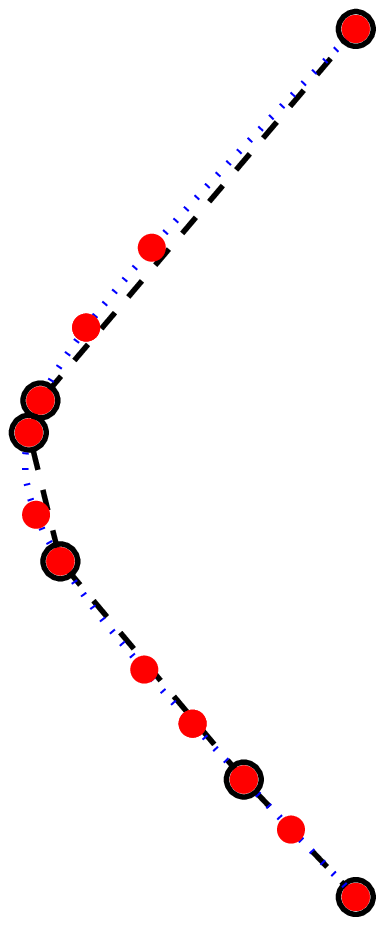}}\hspace{0.4cm}
{\includegraphics[trim= 25mm 7mm 22mm 6mm, clip, width=2.9cm]{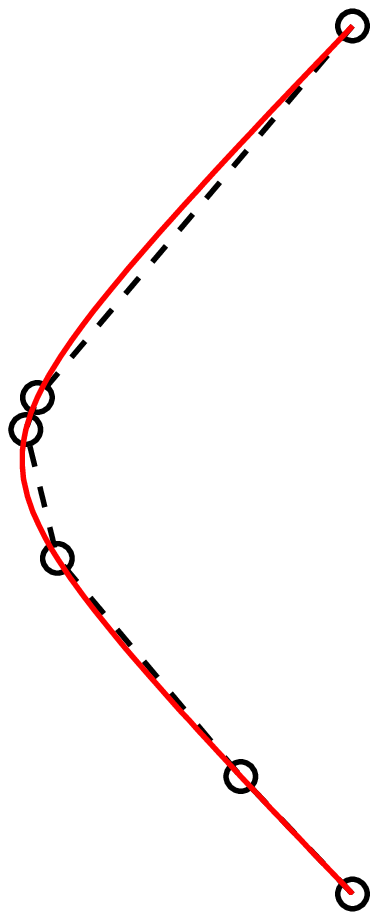}}
\caption{Non-uniform case: reproduction of conic sections from
non-equispaced samples by applying the adaptive subdivision algorithm of
Section \ref{secadapt}. The dotted blue line is the conic
section to be reconstructed. From left to right: points at 1st and
2nd level of refinement; refined polyline after 6 steps of the
algorithm.} \label{fig5b}
\end{figure}

\begin{figure}[h!]
\centering \hspace{-0.3cm}
{\includegraphics[trim
= 30mm 6mm 30mm 6mm, clip, width=2.8cm]{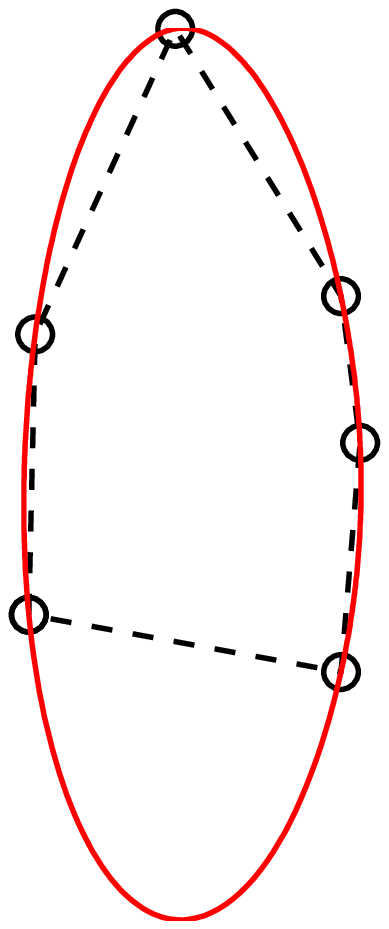}}
\hspace{0.4cm}
{\includegraphics[trim = 10mm 5mm 10mm 6mm, clip,width=3.8cm]{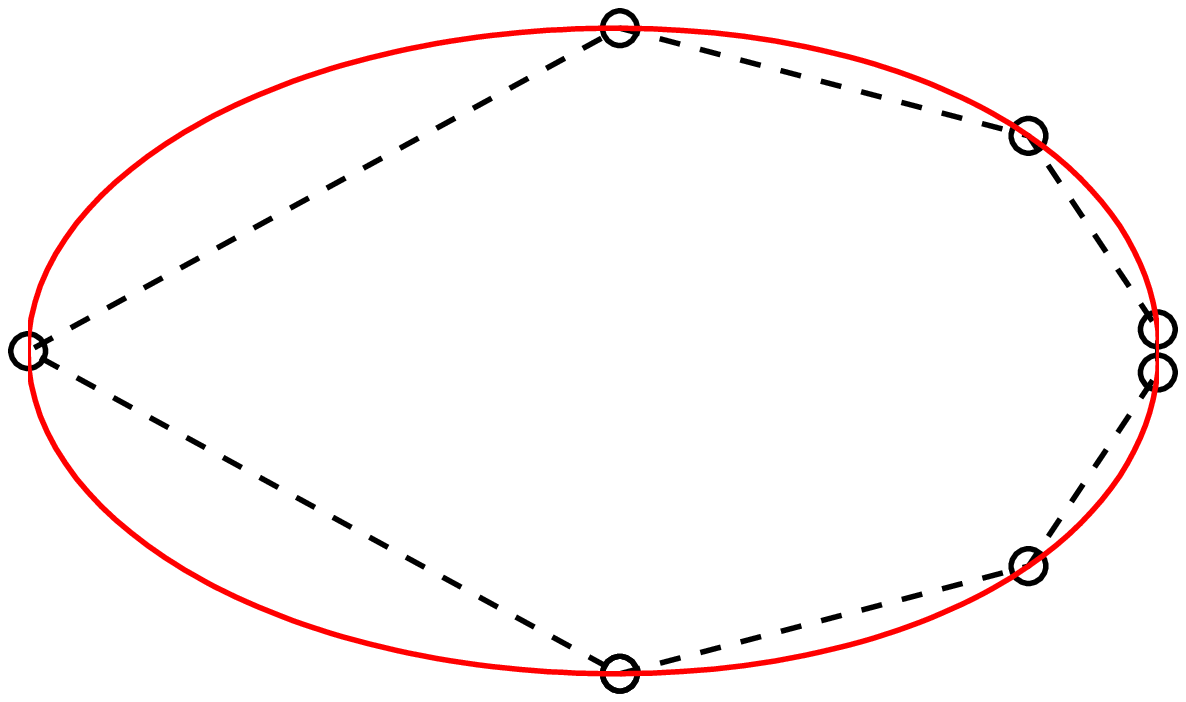}}\hspace{1.2cm}
{\includegraphics[trim = 30mm 6mm 30mm 6mm, clip, width=2.8cm]{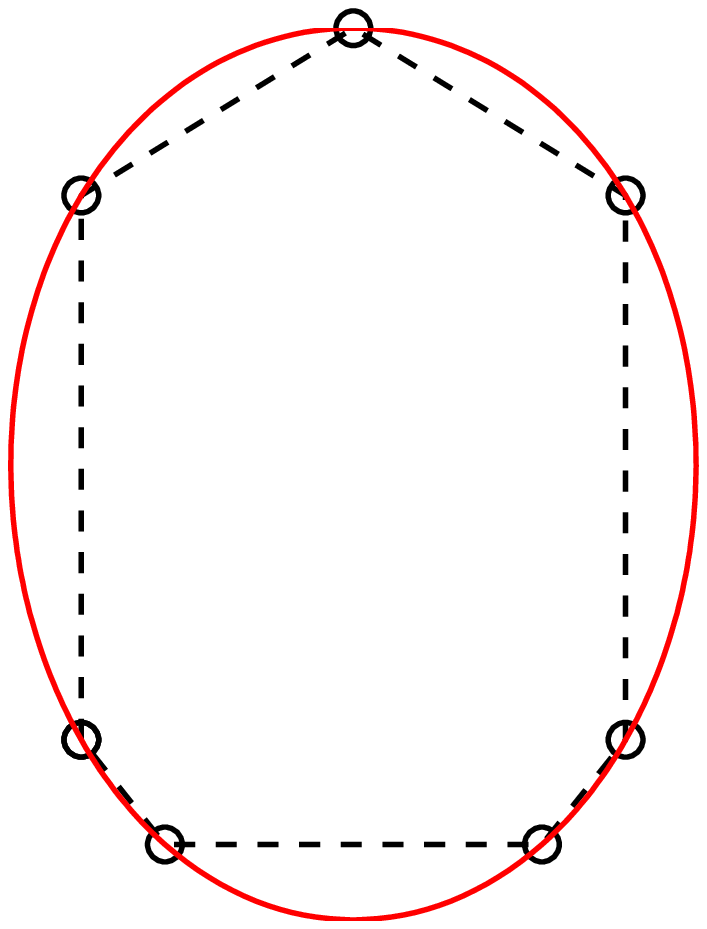}}\\
\vspace{0.3cm}
\hspace{0.1cm}
{\includegraphics
[trim = 45mm 6mm 45mm 6mm, clip,width=2.5cm]
{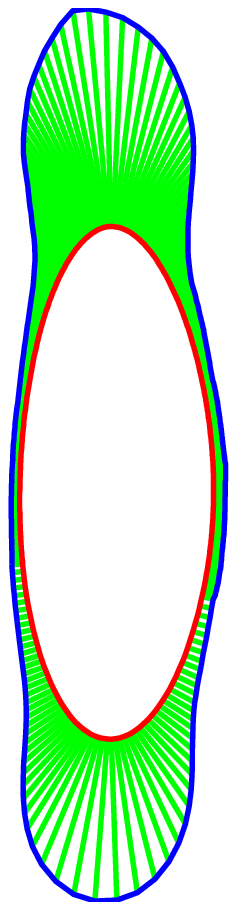}}\hspace{0.1cm}
{\includegraphics[trim = 10mm 5mm 10mm 6mm, clip,width=5.2cm]{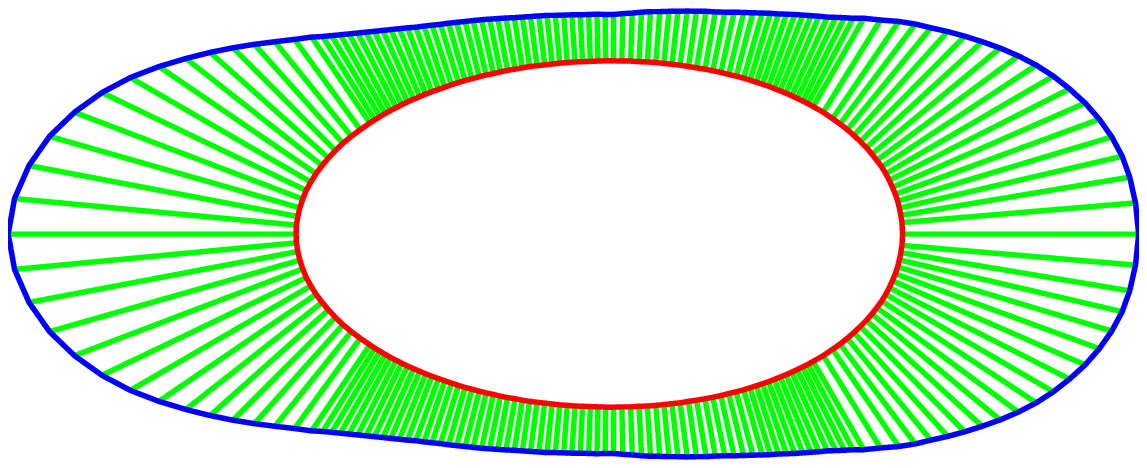}}\hspace{0.3cm}
{\includegraphics[trim = 33mm 6mm 33mm 6mm, clip, width=3.2cm]{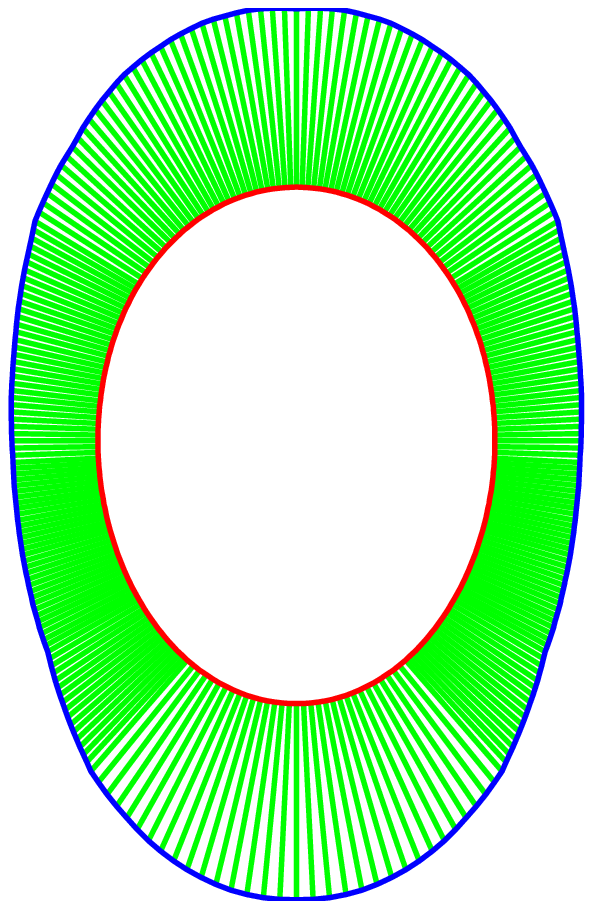}}
\caption{First row: application
examples of the subdivision algorithms of Subsection
\ref{total_conv_algo} (left) and Section \ref{secadapt} (center and right) to totally convex closed polylines that do not come from a conic section. Second row: curvature combs of the considered examples.} \label{fig3}
\end{figure}

\begin{figure}[h!]
\centering
\hspace{-0.3cm}
\subfigure[]{\includegraphics[trim = 30mm 6mm 30mm 6mm, clip, width=2.8cm]{fig77_a.eps}}\hspace{0.3cm}
\subfigure[]{\includegraphics[trim = 30mm 6mm 30mm 6mm, clip, width=2.8cm]{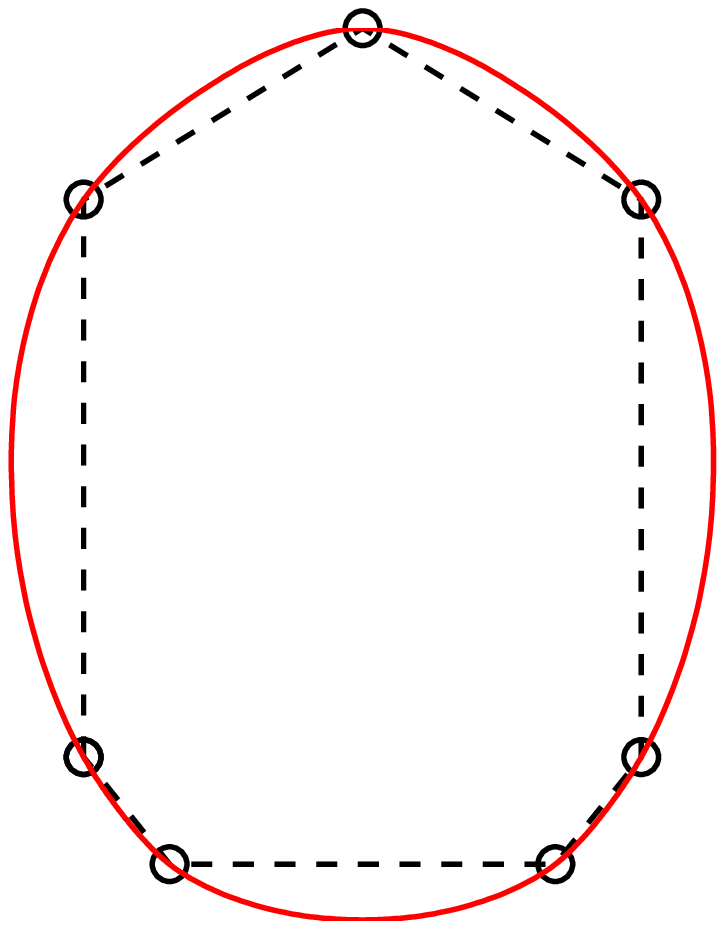}}\hspace{0.3cm}
\subfigure[]{\includegraphics[trim = 30mm 6mm 30mm 6mm, clip, width=2.8cm]{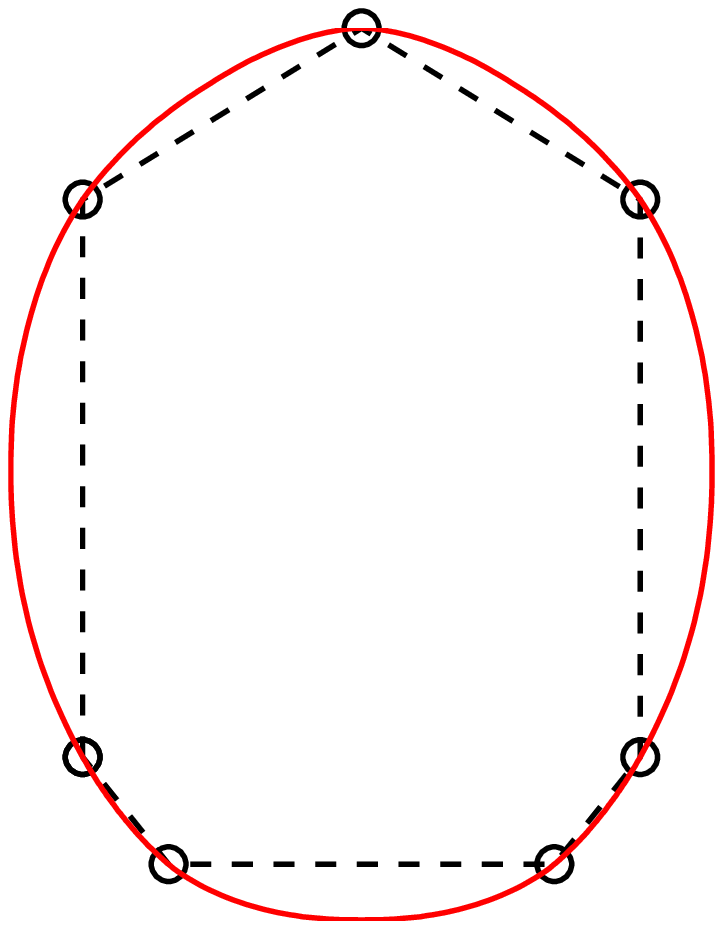}}\hspace{0.3cm}
\subfigure[]{\includegraphics[trim = 30mm 6mm 30mm 6mm, clip, width=2.8cm]{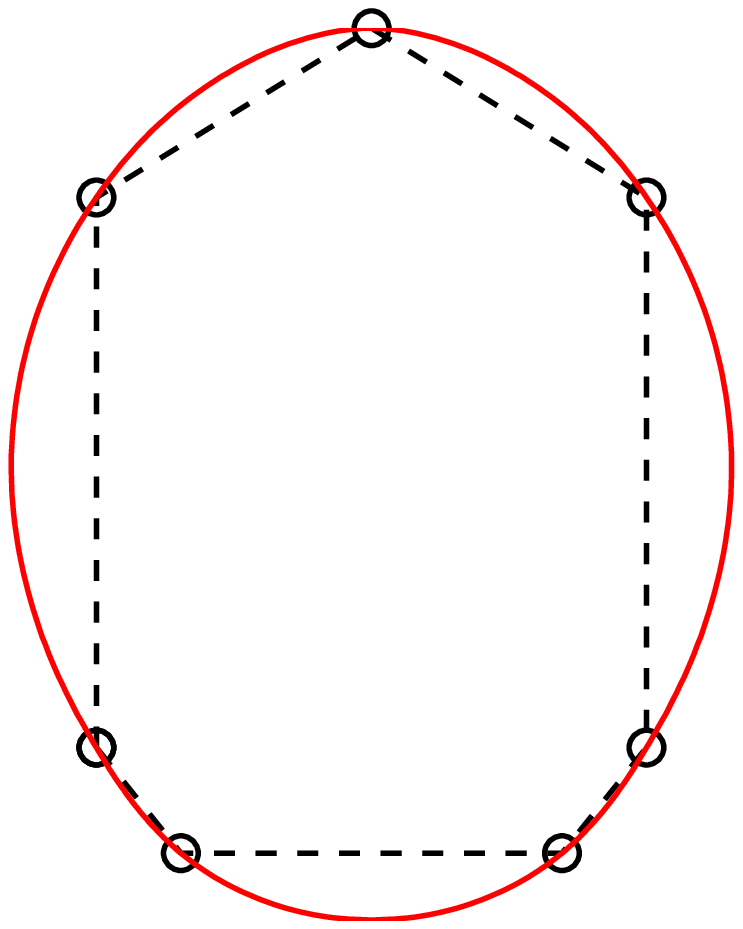}}\vspace{0.3cm}\\
\hspace{-0.3cm}
\subfigure[]{\includegraphics[trim = 33mm -6mm 33mm 6mm, clip, width=2.4cm]{fig77_a_comb.eps}}
\hspace{-0.2cm}
\subfigure[]{\includegraphics[trim = 38mm 10mm 38mm 10mm, clip, width=3.4cm]{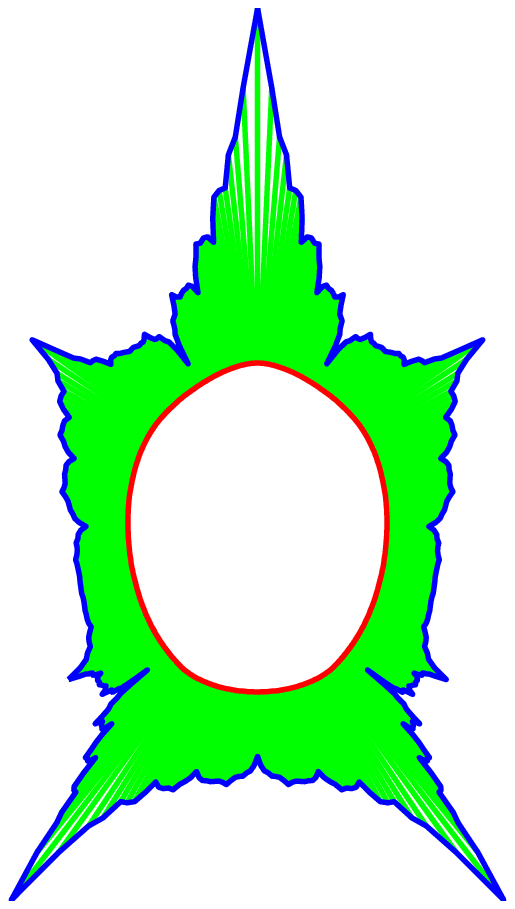}}\hspace{-0.1cm}
\subfigure[]{\includegraphics[trim = 38mm 7mm 38mm 10mm, clip, width=3.4cm]{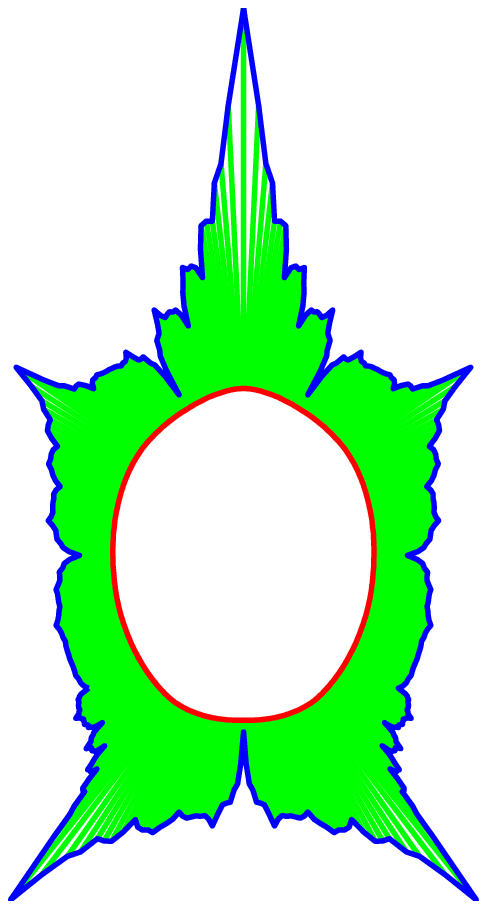}}\hspace{-0.1cm}
\subfigure[]{\includegraphics[trim = 33mm -3mm 33mm 6mm, clip, width=2.9cm]{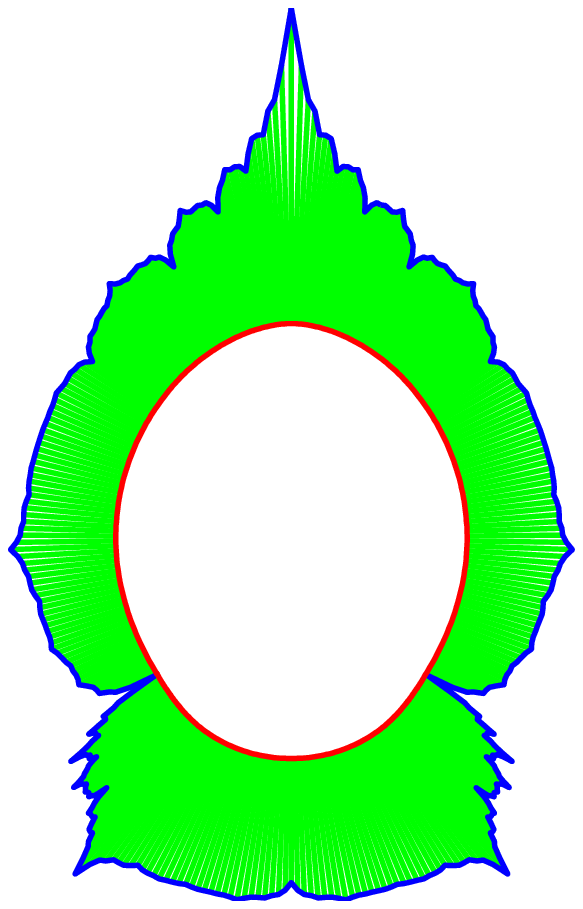}}
\caption{Comparison with the
linear, non-uniform scheme in \cite{BCR09b} and with the non-linear
schemes in \cite{DFH09} and \cite{SD05}. First row: refined
polylines obtained after 6 steps of (a) our adaptive algorithm of Section
\ref{secadapt}; (b) algorithm in \cite{BCR09b} with chord
length parameterization; (c) algorithm in \cite{DFH09} with chord
length parameterization; (d) algorithm in \cite{SD05}. Second row:
corresponding curvature combs.} \label{fig4}
\end{figure}

\begin{figure}[h!]
\centering \hspace{-0.7cm}
\resizebox{5.2cm}{!}{\includegraphics{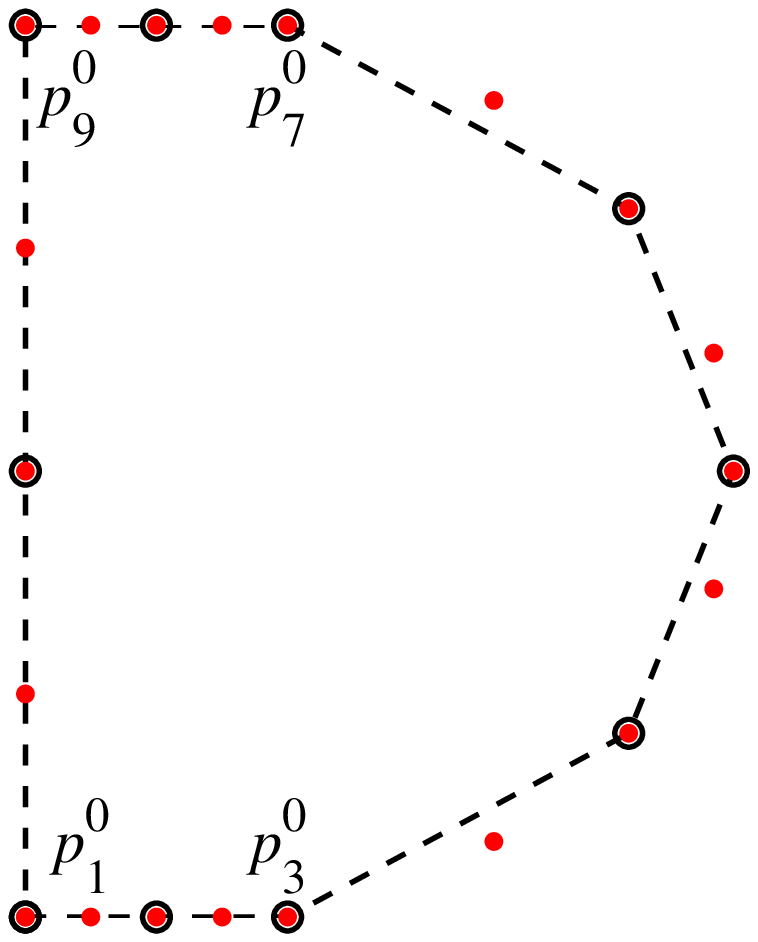}}\hspace{-1.1cm}
\resizebox{5.2cm}{!}{\includegraphics{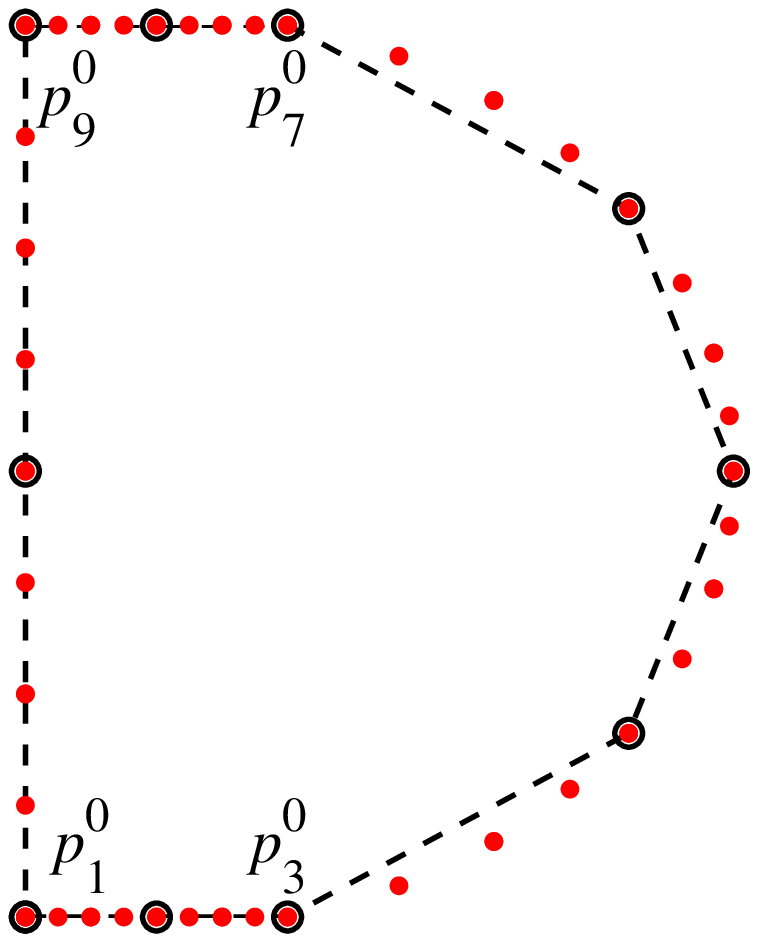}}\hspace{-1.1cm}
\resizebox{5.2cm}{!}{\includegraphics{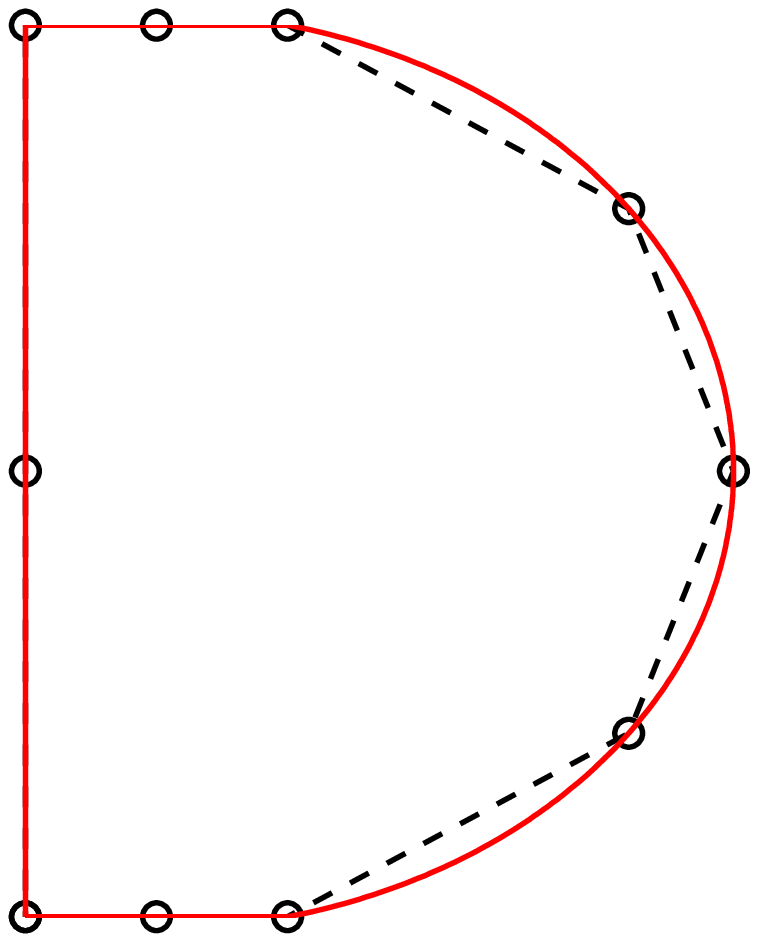}}
\caption{Application example of the subdivision algorithm of
Subsection \ref{non_conv_algo} to a closed sequence of data
containing collinear vertices. The labels in the figure denote the
end points of consecutive subpolygons. From left to right: points at
1st and 2nd level of refinement; refined polyline after 6 steps of
the algorithm.} \label{fig11}
\end{figure}

\begin{figure}[h!]
\centering \hspace{-0.2cm}
\resizebox{4.3cm}{!}{\includegraphics{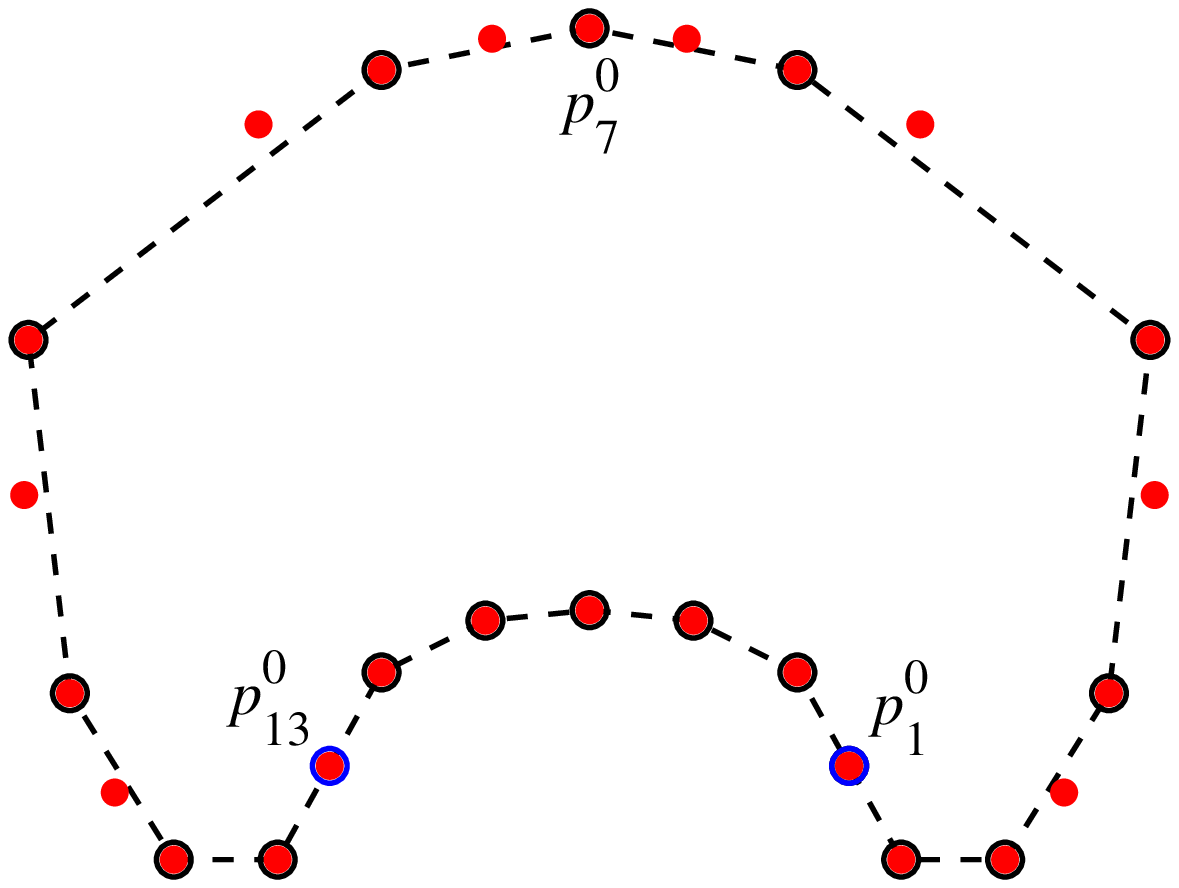}}\hspace{-0.1cm}
\resizebox{4.3cm}{!}{\includegraphics{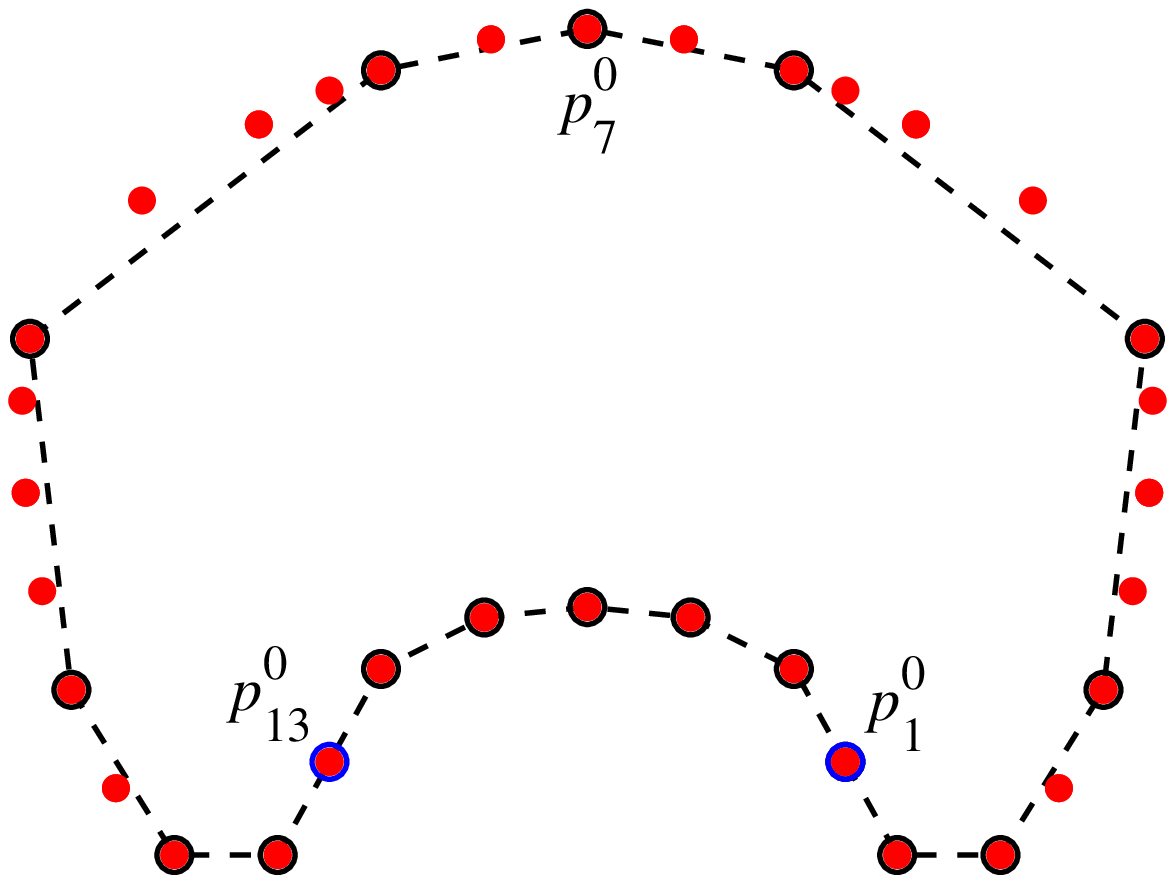}}\hspace{-0.1cm}
\resizebox{4.3cm}{!}{\includegraphics{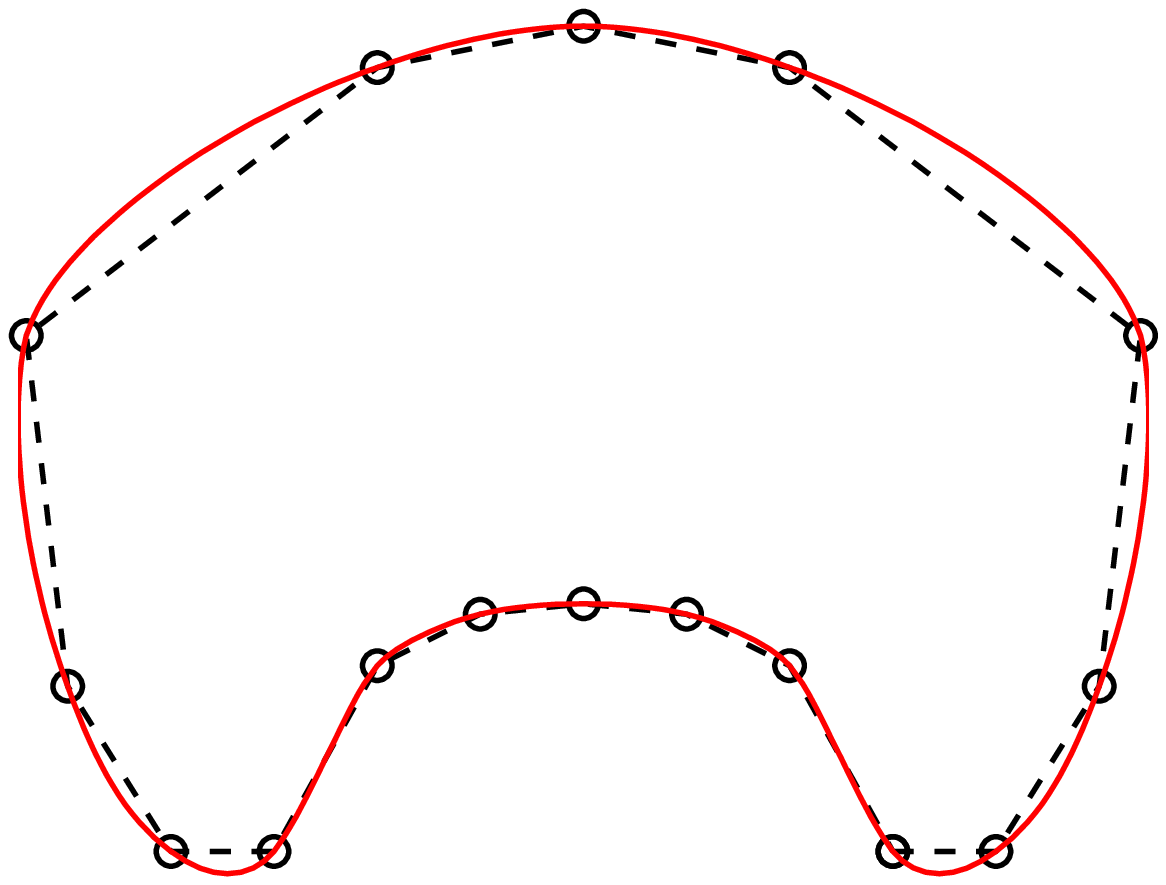}}
\caption{Application example of the adaptive subdivision algorithm of
Section \ref{secadapt} to a closed sequence of non-convex
data. The labels in the figure denote the end points of consecutive
subpolygons. From left to right: points at 1st and 2nd level of
refinement; refined polyline after 6 steps of the algorithm.}
\label{fig12}
\end{figure}

\begin{figure}[h!]
\centering \hspace{-0.8cm}
\resizebox{5.5cm}{!}{\includegraphics{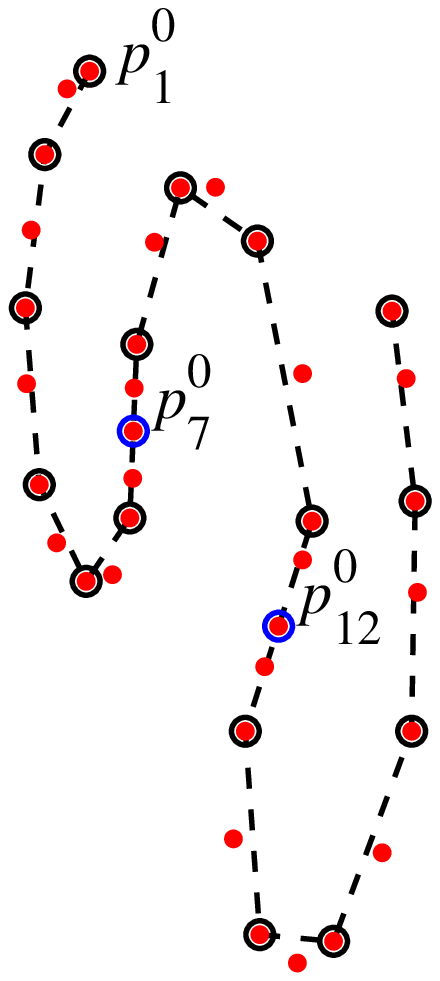}}\hspace{-1.8cm}
\resizebox{5.5cm}{!}{\includegraphics{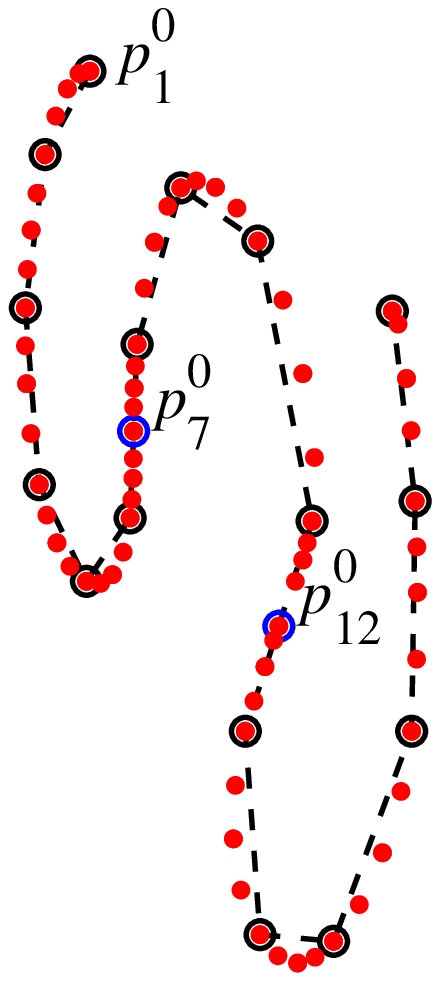}}\hspace{-1.8cm}
\resizebox{5.5cm}{!}{\includegraphics{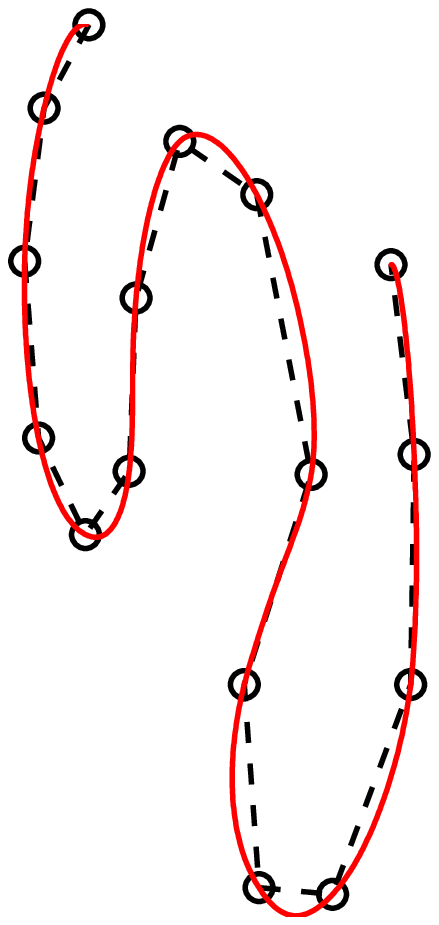}}
\caption{Application example of the subdivision algorithm of
Subsection \ref{non_conv_algo} to an open sequence of non-convex
data. The labels in the figure denote the end points of consecutive
subpolygons. From left to right: points at 1st and 2nd level of
refinement; refined polyline after 6 steps of the algorithm.}
\label{fig13}
\end{figure}

\begin{figure}[h!]
\centering
\hspace{-0.2cm}
{\includegraphics[trim = 50mm 10mm 50mm 10mm, clip, width=3.5cm]{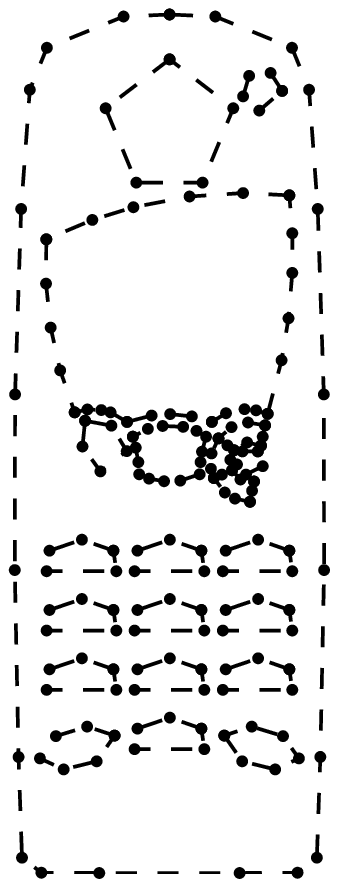}}\hspace{0.4cm}
{\includegraphics[trim = 50mm 10mm 50mm 10mm, clip, width=3.5cm]{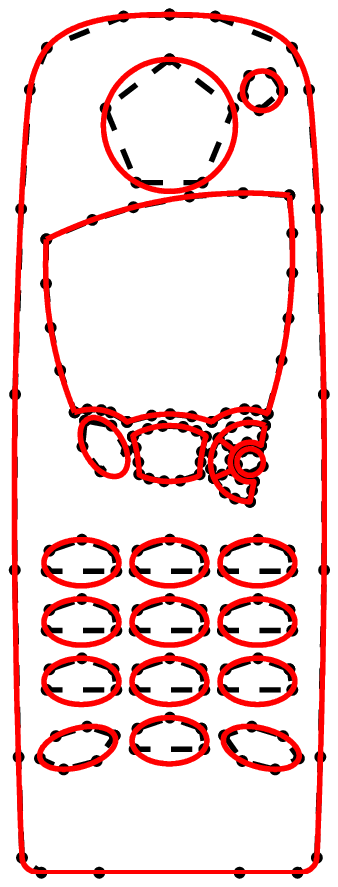}}\hspace{0.4cm}
{\includegraphics[trim = 50mm 10mm 50mm 10mm, clip, width=3.5cm]{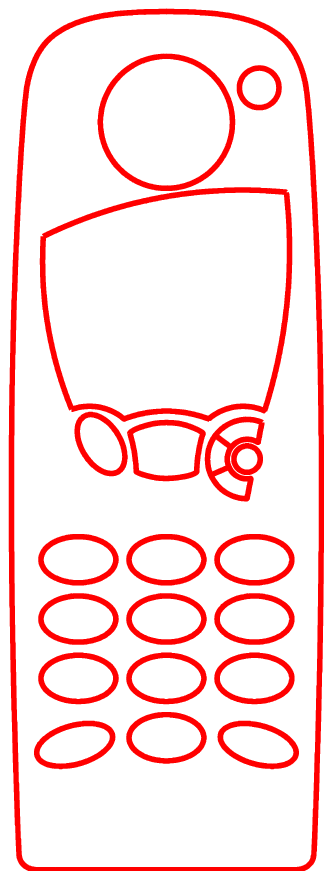}}\vspace{0.5cm}\\
\hspace{-0.3cm} {\includegraphics[trim = 25mm 10mm 25mm 10mm, clip,
width=4.3cm]{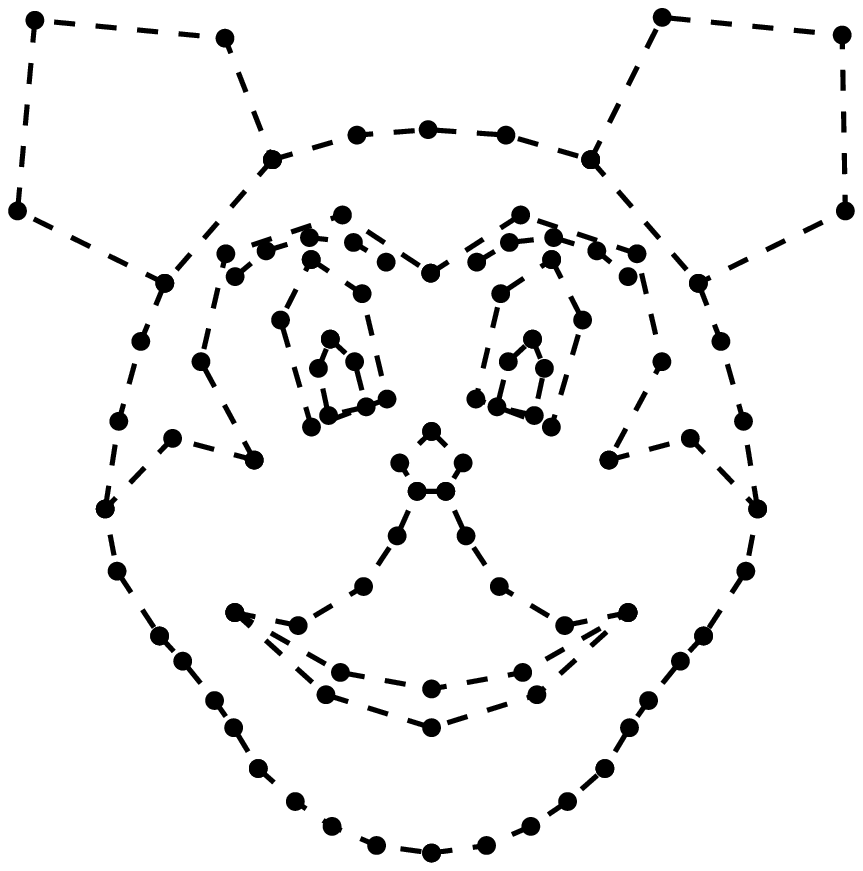}}\hspace{-0.1cm}
{\includegraphics[trim = 25mm 10mm 25mm 10mm, clip,
width=4.3cm]{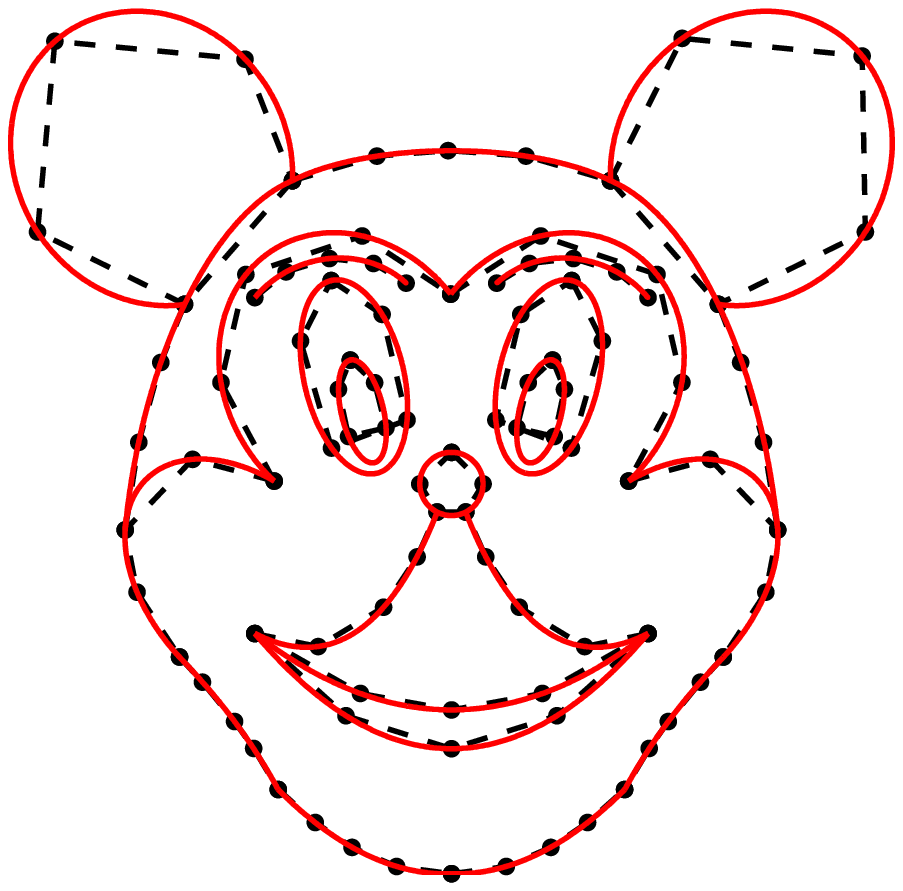}}\hspace{-0.1cm}
{\includegraphics[trim = 25mm 10mm 25mm 10mm, clip,
width=4.3cm]{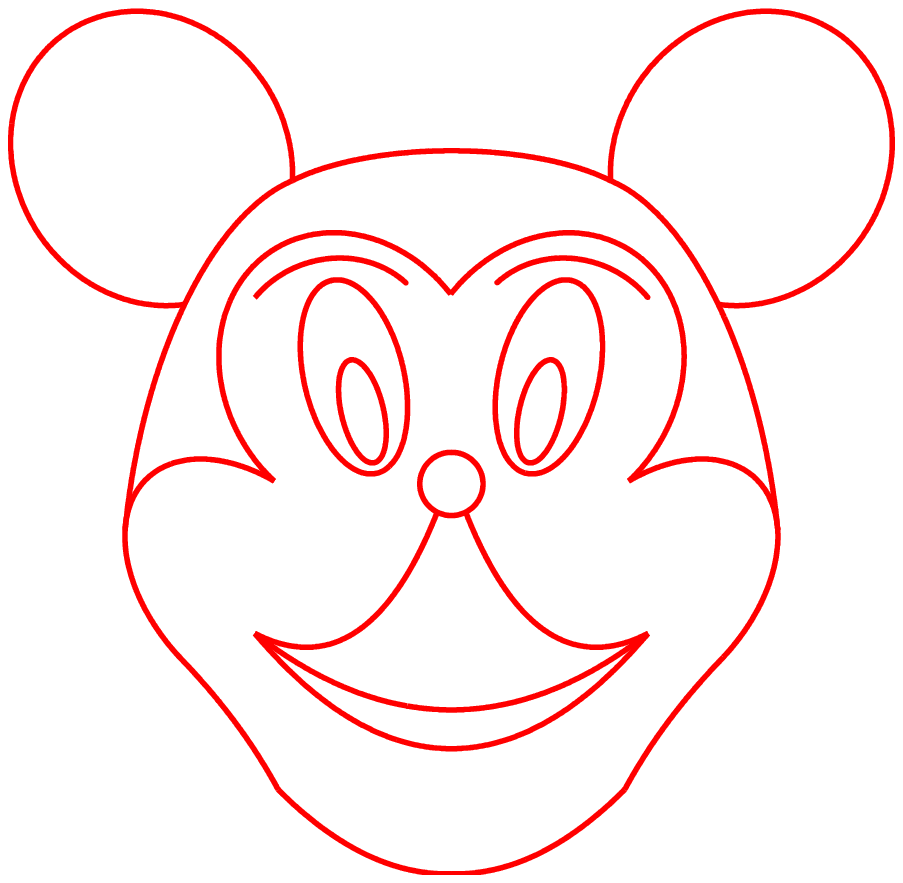}}
\caption{Application examples of the subdivision algorithm of
Subsection \ref{non_conv_algo} to complex data. Left: starting
polylines. Center and Right: refined polylines obtained after 7
steps of our algorithm. Data of first row: courtesy of think3.}
\label{fig_complex1}
\end{figure}

\begin{figure}[h!]
\centering
\hspace{-0.4cm} {\includegraphics[trim = 25mm 10mm 20mm 10mm, clip,
width=4.0cm]{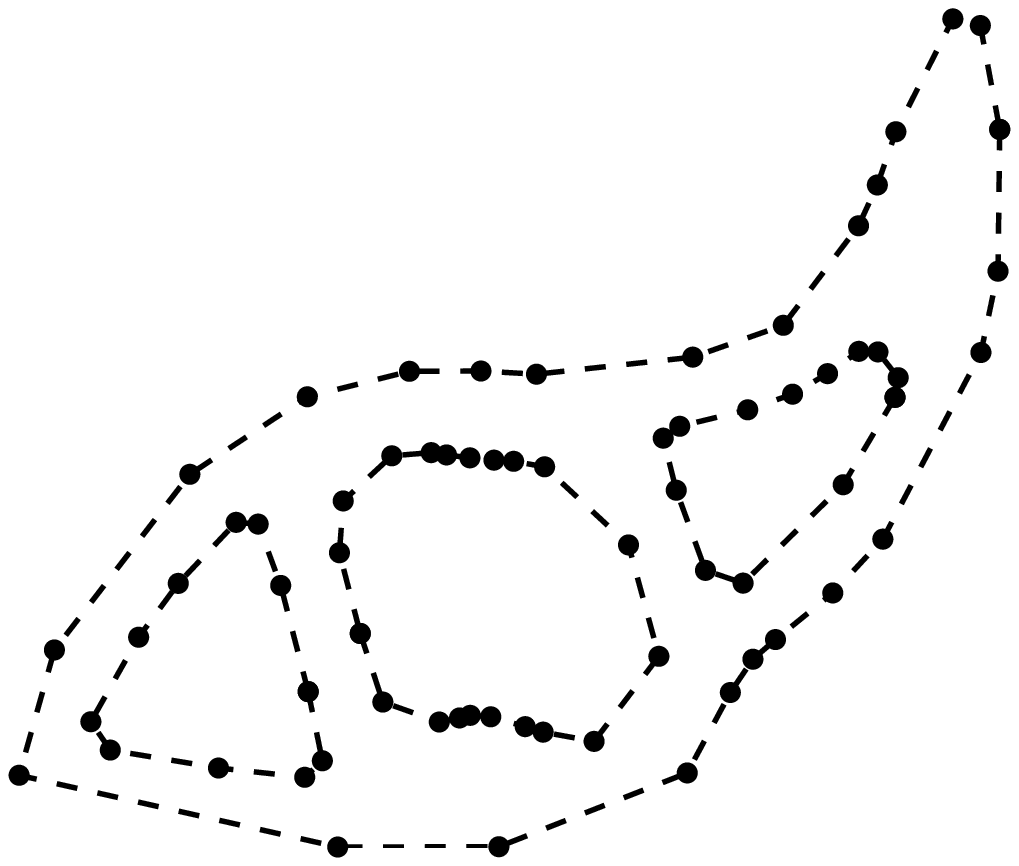}}
\hspace{-0.1cm} {\includegraphics[trim = 25mm 10mm 20mm 10mm, clip,
width=4.0cm]{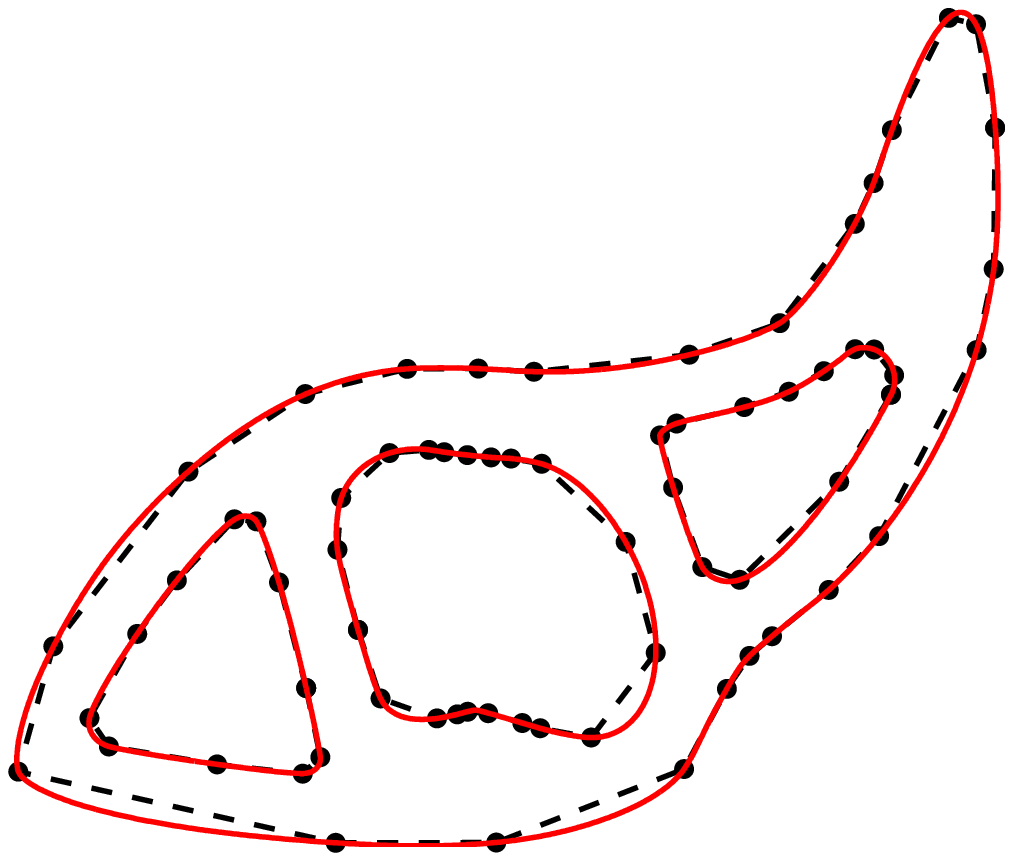}}
\hspace{-0.1cm} {\includegraphics[trim = 25mm 10mm 20mm 10mm, clip,
width=4.0cm]{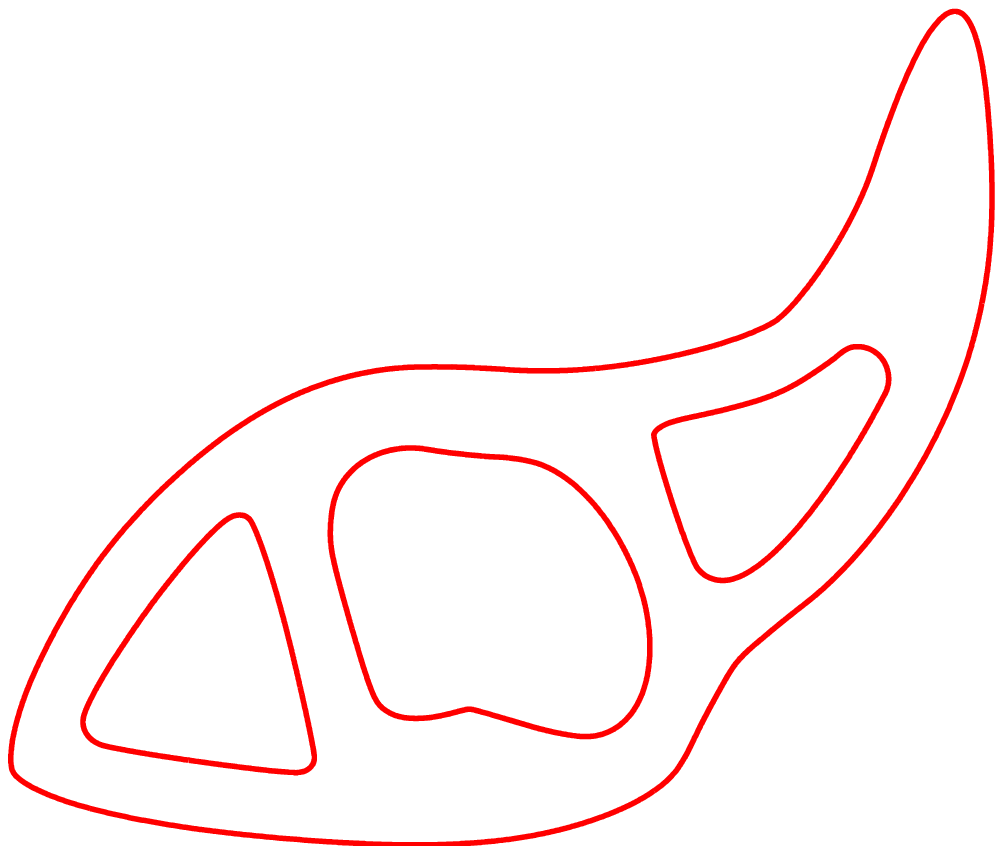}}
\caption{Application example of the subdivision algorithm of
Subsection \ref{non_conv_algo} to complex data. Left: starting
polylines. Center and Right: refined polylines obtained after 7
steps of our algorithm. Data are courtesy of think3.}
\label{fig_complex2}
\end{figure}

As concerns the forthcoming examples, we start by applying the subdivision algorithm of Subsection \ref{total_conv_algo} to totally convex closed and open polylines with nearly uniform edges (Figures \ref{fig7},
\ref{fig9}), and successively we exploit the adaptive version of the scheme in the case of polylines with highly non-uniform edges (Figures \ref{fig6}, \ref{fig8}). As it appears, the generated curves are
always convex and visually pleasing.\\
Then, concerning Figures \ref{fig5a} and \ref{fig5b}, the goal is to illustrate the conic
precision property of the proposed algorithms in both the uniform and
non-uniform cases. Moving from top to bottom, the four plots in the
figures have been generated by repeated application of the
subdivision scheme to points sampled from a circle, an ellipse, a
parabola and a hyperbola.\\
We then continue by showing that the curves computed through
the algorithms presented in Subsection \ref{total_conv_algo} and Section \ref{secadapt} are really
artifact free. In fact, although a limit curve can be apparently
artifact free, it is hard to tell from the display on the screen if
it is acceptable or not. Two curves may look very similar on the
screen, but their curvature plots may reveal important differences.
The most commonly used tool for revealing significant shape
differences is provided by the curvature comb of the curve. In
pictures \ref{fig3}, \ref{fig4} we have used the graphs of the
discrete curvature combs of the refined polylines to show that the
limit curves generated by our algorithms are indeed artifact free.\\
In particular, if we compare the results we get by
refining the polyline in Figure \ref{fig4}(a) for data, that do not
come from a conic section, through our adaptive algorithm and through the
subdivision algorithms in \cite{BCR09b}, \cite{DFH09} and
\cite{SD05} (Figures \ref{fig4}(b), (c) and (d) respectively), they
are only apparently very similar. Yet their curvature combs reveal
substantial differences showing that neither every non-linear nor
every non-uniform subdivision scheme is indeed artifact free (see
Figures \ref{fig4} (e)-(f)-(g)-(h)).\\
We close this section by illustrating the results of the subdivision
algorithm of Subsection \ref{non_conv_algo}. In the first example we
take the D-shape polyline of \cite{cinesi} in order to show the
ability of the scheme to reproduce collinear vertices (see Figure
\ref{fig11}). In the second and third examples we apply the
generalized subdivision scheme to
a closed, respectively open, sequence of non-convex data to illustrate its $G^1$ continuity and shape-preserving interpolation properties (Figures \ref{fig12}, \ref{fig13}).
The last three examples in Figures \ref{fig_complex1} and \ref{fig_complex2} deal with more complex polylines that
respectively represent the cover of a mobile phone, Micky Mouse face
and a bottle opener.\\
The shapes in Figure \ref{fig_complex1} are
designed by a collection of rectilinear and convex segments where
the most of the latter ones have been sampled from conic sections.
The shape in Figure \ref{fig_complex2} is made of 4 independent closed
polygons, some of which are non-convex. The data of the first
example in Figure \ref{fig_complex1} and that of Figure
\ref{fig_complex2} are courtesy of the CAD company think3
(www.think3.com). In all the considered experiments the proposed
subdivision scheme turns out to work very well and clearly manifests
all its characteristic features described in the previous sections.

\section{Conclusions} \label{concl}

Even though in the last years important steps forward have been
taken both in the construction and analysis of interpolatory
subdivision schemes \cite{S05}, several problems are still open and
need to be tackled in order to increase the strength and popularity
of subdivision in more and more fields of application.

First of all, unlike the non-interpolatory subdivision schemes, the
interpolatory ones usually generate shapes of inferior quality
because, if applied to points with an irregular distribution, they
provide a limit curve with more convexity changes than the starting
polygon. Since, in several applications it is important to guarantee
shape preservation, in this paper we have described a new
interpolating subdivision algorithm enjoying this important
property.

Because in CAGD it is also often necessary to have schemes able to
generate classical geometric shapes, we have enriched our
interpolating subdivision scheme with the capability of including
the exact representation of all conic sections. The different
methods of the literature \cite{BCR07,BCR09a,R09a} give solutions
only when the assigned points have a regular distribution. These
linear subdivision schemes use non-stationary refinement rules
associated with functional spaces defined via the union between
polynomial and exponential functions with a free parameter.

The idea we have explored in this paper is to provide a subdivision
scheme in which the conic section reproduction is obtained by
adaptive geometric constructions on the given points. The advantage
of doing this is that the presented non--linear scheme is able to
adapt itself to any data configuration, i.e., to arbitrary
irregularly distributed point sequences. Due to the underlying
construction, the properties of shape preservation, conic
reproduction as well as the proof of $G^1$ continuity of the limit
curve, follow straightforwardly. The method has been illustrated by
several significant examples.


\end{document}